\documentclass[]{amsart}
\usepackage[francais]{babel}
\usepackage{amsfonts,amssymb,epsfig}
\usepackage[latin1]{inputenc}
\usepackage[all]{xy}
\usepackage{amsmath}
\usepackage{delarray}
\usepackage{amssymb}
\usepackage{amscd}
\usepackage{diagrams}
\theoremstyle{definition}

\renewcommand{\epsilon}{\varepsilon}

\theoremstyle{remark}

\numberwithin{equation}{section}

\newcommand{\IE}{\text\textquestiondown}

\newcommand{\go}{\mathfrak}
\newcommand{\gr}{\mathcal}
\newcommand{\mbb}{\mathbb}
\newcommand{\mbf}{\mathbf}

\def\cad{c'est-à-dire }
\newcommand{\St}{$\check{\textrm{S}}$tefan }
\let\til=~

\newcommand{\lab}{labelstyle}

\newcommand{\sst}{\scriptstyle}
\newcommand{\ssst}{\scriptscriptstyle}
\newcommand{\tnat}{{\overset{\bullet}{}}}

\newarrow{ClopTo}{sqhook}---{vee}
\newarrow{OpTo}{hooka}---{vee}
\newarrow{iTo}{triangle}--->
\newarrow{ITo}{blacktriangle}--->
\newarrow{sTo}----{triangle}
\newarrow{STo}----{blacktriangle}

\newarrow{iLine}{triangle}----
\newarrow{ILine}{blacktriangle}----
\newarrow{Equal}=====%
\newarrow{DotTo}....>
\newarrow{DotsTo}....{triangle}
\newarrow{DotSTo}....{blacktriangle}
\newarrow{DotiTo}{triangle}...>
\newarrow{DotITo}{blacktriangle}...>
\newarrow{DashTo}{dash}{dash}{dash}{dash}>
\newarrow{DashsTo}{dash}{dash}{dash}{dash}{triangle}
\newarrow{DashSTo}{dash}{dash}{dash}{dash}{blacktriangle}
\newarrow{DashiTo}{triangle}{dash}{}{dash}>
\newarrow{DashITo}{blacktriangle}{dash}{dash}{dash}>
\newarrow{DashLine}{dash}{dash}{dash}{dash}{dash}
\newarrow{DashiLine}{triangle}{dash}{dash}{dash}{dash}
\newarrow{DashILine}{blacktriangle}{dash}{dash}{dash}{dash}
\newarrow{DotLine}.....

\def\rsum{\raisebox{2pt}{$\:\scriptscriptstyle+\:$}}
\def\lsum{\raisebox{2pt}{$\:\scriptscriptstyle+\:$}}
\def\dsum{\vphantom b\scriptscriptstyle+}
\def\usum{\vphantom b\scriptscriptstyle+}

\newarrowfiller{+}\rsum\lsum\dsum\usum
\newarrow{SumSTo}{+}{+}{+}{+}{blacktriangle}
\newarrow{SumSTo}{+}{+}{+}{+}{blacktriangle}
\newarrow{SumLine}{+}{+}{+}{+}{+}
\newarrow{SumTo}{+}{+}{+}{+}>
\newarrow{SumiTo}{triangle}{+}{+}{+}>

\def\rrnd{\raisebox{1.5pt}{$\:\scriptscriptstyle{\circ}\:$}}
\def\lrnd{\raisebox{1.5pt}{$\:\scriptscriptstyle{\circ}\:$}}
\def\drnd{\vphantom b\scriptscriptstyle{\circ}}
\def\urnd{\vphantom b\scriptscriptstyle{\circ}}

\newarrowfiller{rnd}{\rrnd}{\lrnd}{\drnd}{\urnd}
\newarrow{RndSTo}{rnd}{rnd}{rnd}{rnd}{blacktriangle}
\newarrow{RndLine}{rnd}{rnd}{rnd}{rnd}{rnd}
\def\merod{$\widetilde{\textbf{Grpd}}(\mbf{D})$ }
\def\holod{$[\textbf{Grpd}](\mbf{D})$ }

\def\rbul{\raisebox{1.5pt}{$\:\scriptscriptstyle{\bullet}\:$}}
\def\lbul{\raisebox{1.5pt}{$\:\scriptscriptstyle{\bullet}\:$}}
\def\dbul{\vphantom b\scriptscriptstyle{\bullet}}
\def\ubul{\vphantom b\scriptscriptstyle{\bullet}}

\newarrowfiller{bul}{\rbul}{\lbul}{\dbul}{\ubul}
\newarrow{BulSTo}{bul}{bul}{bul}{bul}{blacktriangle}
\newarrow{BulLine}{bul}{bul}{bul}{bul}{bul}
\def\imm{\text{\scriptsize{\textcircled{i}}}}
\def\rimm{\imm}
\def\limm{\imm}
\def\uimm{\imm}
\def\dimm{\imm}
\newarrowmiddle{imm}{\rimm}{\limm}{\uimm}{\dimm}
\newarrow{imTo}--{imm}->
\newarrow{plfTo}{LaTeX}-{imm}->
\newarrow{fidTo}{curlyvee}-{imm}->
\newarrow{istTo}{boldlittlevee}-{imm}->

\newarrow{qimTo}{o}--->
\def\etl{\text{\scriptsize{\textcircled{e}}}}
\def\retl{\etl}
\def\letl{\etl}
\def\uetl{\etl}
\def\detl{\etl}
\newarrowmiddle{etl}{\retl}{\letl}{\uetl}{\detl}
\newarrow{etaTo}--{etl}->
\newarrow{etasTo}--{etl}-{triangle}
\newarrow{ploTo}{triangle}-{etl}->
\newarrow{pofTo}{littleblack}-{etl}->
\def\sub{\text{\scriptsize{\textcircled{s}}}}
\def\rsub{\sub}
\def\lsub{\sub}
\def\usub{\sub}
\def\dsub{\sub}
\newarrowmiddle{sub}{\rsub}{\lsub}{\usub}{\dsub}
\newarrow{subTo}--{sub}->

\def\trl{\pitchfork}
\def\tre{\top}

\def\act{\text{{\fboxsep1pt\framebox{\tiny{\sffamily act}}}}}

\def\exa{\text{{\fboxsep1pt\framebox{\tiny{\sffamily exa}}}}}

\def\pb{\text{\raisebox{1.6pt}{{\fboxsep0.7pt\framebox{\tiny{\sffamily pb}}}}}}
\def\po{\text{\raisebox{1.6pt}{{\fboxsep0.7pt\framebox{\tiny{\sffamily po}}}}}}

\def\pbk{\SEpbk}

\newcommand{\gpd}[2]{#1\rightrightarrows #2}

\begin{document}
\title[]{Groupoïdes de Lie et Feuilletages}
\author{Jean PRADINES}
\address{26, rue Alexandre Ducos, F31500 Toulouse, France}
\email{jpradines@wanadoo.fr} 


\thanks{\emph{Mots clefs}: groupoïdes, Lie, feuilletages, \St.}

\thanks{-- L'auteur est redevable à Paul Taylor pour l'utilisation de son extension \emph{diagrams} dans l'exécution de nombreux schémas et diagrammes commutatifs.}

\thanks{-- Le point de vue présenté dans le résumé ci-dessus, que l'on s'efforce ici de vulgariser, s'est dégagé clairement à la suite des discussions entre A. Connes (et ses élèves: M. Hilsum, G. Skandalis), A. Haefliger, et W.T. van Est, lors du Colloque sur la \emph{Structure Transverse des Feuilletages}, tenu à Toulouse en 1982, à l'initiative de P. Cartier (voir le volume 116 de \emph{Astérisque}, SMF, Paris, 1984).}

\date{}

\begin{abstract} 
On étudie la relation entre les groupoïdes de Lie (groupoïdes différen\-tiables d'Ehresmann) et les feuilletages (avec singularités possibles) au sens de Suss\-mann-\St\!. On souligne la façon dont l'interaction entre une structure al\-gé\-bri\-que\-ment triviale (celle de graphe de relation d'équivalence) et une structure différentiable fait naître des groupes (groupes d'holonomie), dont la variation, discontinue, se lit dans l'objet unificateur (algébrico-différentiable) constitué par le groupoïde d'holonomie.\\
La classe d'équivalence de Morita de ce dernier (notion algébrico-dif\-fé\-ren\-tiable) munit l'espace des feuilles d'une structure, en un sens généralisé (structure transverse du feuilletage); celle-ci, notamment, attache, aux points de cet espace, des groupes, dont les sauts ne se comprennent bien que dans le groupoïde de Lie.
\end{abstract}

\maketitle 

\tableofcontents

\section{Introduction.}
\label{int}
Citant la conférence de Jean-Marie Souriau donnée lors des présentes \emph{Journées}, nous invoquerons la Géométrie comme source de quantification, et, selon F. Klein, les groupes comme source de la Géométrie\,\footnote{Et des états quantiques dans la conférence de J.-M. Souriau.}, groupes qu'il convient le plus souvent, à la suite de S. Lie, de munir d'une structure de variété\,\footnote{Voire, d'après J.-M. Souriau, d'une difféologie.}.

Toutefois nous voudrions ici plaider pour un élargissement de cette dernière source, consistant à remplacer les groupes de Lie par les \emph{groupoïdes de Lie}, \emph{alias} les groupoïdes différentiables, introduits, sous leur forme générale, par Charles Ehresmann en 1958 (voir l'annexe \ref{nh} pour des précisions historiques et terminologiques). Ces derniers sont en outre, nous le verrons, étroitement liés au second thème de ces \emph{Journées}, à savoir les feuilletages, puisque les orbites d'un groupoïde de Lie définis\-sent un \emph{feuilletage avec singularités} au sens de Sussmann-$\check{\textrm{S}}$tefan.

Nous nous attacherons ici seulement à ce deuxième aspect, faute de compétence suffisante pour explorer les groupoïdes de Lie en tant que sources possibles de quantification, ce que nous pensons pourtant être une voie prometteuse\,\footnote{Il conviendrait dans ce sens (ce que nous ne ferons pas ici) de considérer simultanément les objets infinitésimaux correspondant aux groupoïdes de Lie, à savoir les algébroïdes de Lie, que nous avions introduits en 1967 en vue d'étendre la correspondance groupes de Lie/algèbres de Lie.}. Les groupoïdes (et algébroïdes) de Lie sont des objets purement géométriques de dimension finie. Leurs espaces de sections ou de fonctions donnent naissance à des objets de dimension infinie, dont beaucoup de propriétés se lisent déjà de façon très naturelle au niveau géométrique. Mais nous fermerons ici cette parenthèse.

\subsection{Au carrefour de trois structures.}
Avant d'entrer dans plus de détails, annonçons d'abord ici informellement quelques idées directrices.

\smallskip
De même que les groupes de Lie, les \emph{groupoïdes de Lie} font interagir une structure \emph{algébrique} et une structure de \emph{variété}. Mais la structure algébrique de groupoïde unifie à son tour deux structures, qui constituent deux formes possibles de sa dé\-gé\-né\-res\-cence: la structure de \emph{groupe}, et celle de \emph{graphe} (de relation d'équivalence), laquelle passe (regrettablement) le plus souvent totalement inaperçue, du fait de son apparence trop \og triviale \fg\,d'un point de vue purement algébrique.

La notion naturelle de morphisme pour cette structure mixte est celle de foncteur différentiable entre groupoïdes de Lie, que nous appellerons ici brièvement \emph{morphismes}. Il peut être aussi utile de considérer des \emph{transformations naturelles} [McL] (différentiables ) entre de tels foncteurs.

La simple évocation de quelques formules familières (parmi beaucoup d'autres), que nous laissons au lecteur le soin de réinterpréter, selon le cas, en termes de tels morphismes ou de transformations naturelles, nous permettra dès maintenant de souligner l'intérêt de pouvoir notamment considérer comme cas particuliers d'une \emph{même} structure un graphe d'équivalence du côté source et un groupe du côté but:
$$\overrightarrow{CA}=\overrightarrow{CB}+\overrightarrow{BA}\,;$$
$$g_{k,i} (x)=g_{k,j} (x)\cdot g_{j,i} (x)\quad(x\in {U_i}\cap{U_j}\cap{U_k}),$$
$$g'_{j,i}={h_j}(x)^{-1}\cdot g_{j,i}(x)\cdot{h_i}(x)\quad(x\in{U_i}\cap{U_j})\quad\text{(cocycles cohomologues)}\,\footnote{Il faudrait ici considérer un recouvrement $(U_i)_{i\in I}$ d'une variété $B$, la variété $U$ somme des $U_i$, et le graphe $R$ de la relation d'équivalence sur $U$ définie par la projection canonique $U\rightarrow B$\,. Les \emph{cocycles} $g:R\rightarrow G,\,(j,i;x)\mapsto g_{j,i} (x)$ sont alors des cas particuliers de morphismes, donc de foncteurs, et la condition, pour deux cocycles, d'être \emph{cohomologues} apparaît comme un cas particulier d'isomorphisme (fonctoriel) de deux foncteurs (voir Mac Lane [McL]). Les cocycles de Haefliger s'obtiennent en prenant pour $G$ un \emph{pseudogroupe}, qui peut être vu comme un groupoïde (étale).};$$
$$\delta(z,x)=\delta(z,y)\cdot\delta(y,x)\,,$$
$$\text{où }\delta:G\times G\rightarrow G, (y,x)\mapsto y\cdot x^{-1}\,\text{est la \og division \fg\, du groupe }G\,\footnote{On pourrait voir que la version infinitésimale de cette dernière formule, laquelle donne un morphisme d'algébroïdes de Lie, n'est autre que l'équation dite de Maurer-Cartan.}.$$

\bigskip

\begin{figure}[!htb]
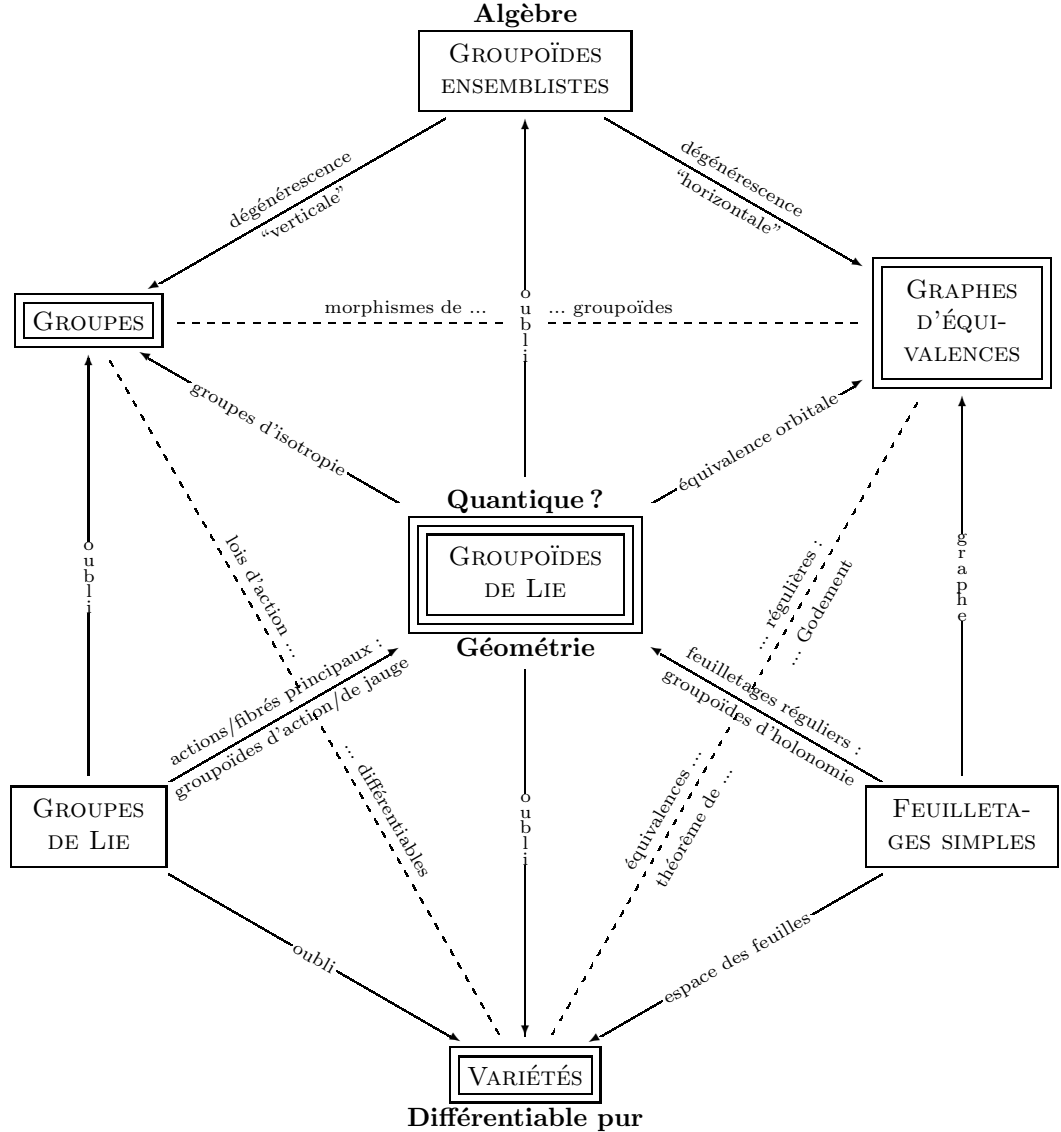

\def\oub{{\begin{diagram}[size=0.6em,objectstyle=\scriptstyle]\text o\\ \text u\\ \text b\\ \text l\\ \text i\\ \end{diagram}}}
\def\gr{{\begin{diagram}[size=0.6em,objectstyle=\scriptstyle]\text g\\ \text r\\ \text a\\ \text p\\ \text h\\ \text e\\ \end{diagram}}}
\def\AA{\overset{\textstyle\textbf{Quantique ?}}{\underset{\textstyle\textbf{Géométrie}}{\boxed{\boxed{\boxed{\begin{array}{c}\textsc{Groupoïdes}\\ \textsc{de Lie}\\ \end{array}}}}}}}
\def\BB{\underset{\textstyle\textbf{Différentiable pur}}{\boxed{\boxed{\text{\textsc{Variétés}}}}}}
\def\CC{\boxed{\boxed{\text{\textsc{Groupes}}}}}
\def\DD{\boxed{\boxed{\begin{array}{c}{\text{\textsc{Graphes}}}\\{\text{\textsc{d'équi-}}}\\{\text{\textsc{valences}}}\\ \end{array}}}}
\def\EE{\overset{\textstyle\textbf{Algèbre}}{\boxed{\begin{array}{c}{\text{\textsc{Groupoïdes}}}\\{\text{\textsc{ensemblistes}}}\\ \end{array}}}}
\def\FF{\boxed{\begin{array}{c}{\text{\textsc{Groupes}}}\\{\text{\textsc{de Lie}}}\\ \end{array}}}
\def\GG{\boxed{\begin{array}{c}{\text{\textsc{Feuilleta-}}}\\{\text{\textsc{ges simples}}}\\ \end{array}}}
\begin{diagram}[w=8.258em,h=4.768em,tight,labelstyle=\scriptstyle]
&&\EE&&\\
&\ldTo^{\text{dégénérescence}}_{\text{``verticale''}}
&&\rdTo^{\text{dégénérescence}}_{\text{``horizontale''}}&\\
\CC&\rDashLine^{\quad\quad\quad\quad\quad\text{morphismes de ...}}&\uTo\til{\oub}&\lDashLine^{\text{... groupoïdes}\quad\quad\quad\quad\quad\quad\quad}&\DD\\
&\rdDashLine(1,3)^{\quad\quad\quad\quad\quad\quad\quad\text{lois d'action ...}}
\luTo\til{\text{groupes d'isotropie}}&&\ruTo\til{\text{équivalence orbitale}}
\ldDashLine(1,3)^{\text{... régulières:}\quad\quad\quad\quad\quad}_{\text{... Godement}\quad\quad\quad\quad\quad}&\\
\uTo\til{\oub}&&\AA&&\uTo\til{\gr}\\
&\ruTo^{\text{actions/fibrés principaux:}}_{\text{groupoïdes d'action/de jauge}}&&\luTo^{\text{feuilletages réguliers:}}_{\text{groupoïdes d'holonomie}}&\\
\FF&&\rdDashLine(1,3)^{\text{... différentiables}\quad\quad\quad\quad\quad}\ldDashLine(1,3)^{\quad\quad\quad\quad\quad\text{équivalences ...}}_{\quad\quad\quad\quad\quad\text{théorême de ...}}\dTo\til{\oub}&&\GG\\
&\rdTo\til{\text{oubli}}&&\ldTo\til{\text{espace des feuilles}}&\\
&&\BB&&\\
\end{diagram}
\caption{\label{ggv}Triple interaction: groupes/graphes/variétés.}
\end{figure}

La figure \ref{ggv} ci-après résume et schématise cette triple interaction entre groupes, graphes, et variétés, qui sera développée dans la suite.

\medskip
Les groupoïdes de Lie figurent au centre (dans un triple cadre), et unifient les trois structures, disposées (dans des cadres doubles) aux sommets d'un triangle, dont les côtés, en traits discontinus, schématisent les interactions de deux de ces structures. Celles-ci engendrent les trois structures mixtes inscrites dans les cadres simples.

\smallskip
Par oubli de la structure algébrique, on obtient, vers le bas, les structures de variétés, dans le domaine du différentiable pur, tandis que l'oubli de la structure de variété donne, vers le haut, les structures de (petits) groupoïdes (ensemblistes ou algébriques), dans le domaine de l'Algèbre pure.

\smallskip
Ces derniers peuvent dégénérer, soit en groupes, soit en graphes (d'équivalence), et, de façon générale, font apparaître:

	-- des groupes, les groupes d'isotropie (qu'il est commode de penser comme la partie \og verticale \fg\,du groupoïde);

	-- une relation d'équivalence, tant sur la base (ensemble des objets) que sur l'ensemble des flèches, dont les classes sont les composantes transitives, ou orbites, du groupoïde et de sa base (décomposition orbitale ou \og horizontale \fg\, du groupoïde).

\smallskip

L'interaction des groupes et des variétés donne naissance aux groupes de Lie et aux lois d'action différentiables. Ces dernières sont, nous le verrons, décrites par certains morphismes de groupoïdes de Lie (que nous appellerons, de façon assez naturelle, \og acteurs \fg\,\footnote{De façon générale nous adopterons le suffixe \emph{eur} pour les morphismes de groupoïdes, \emph{alias} fonct\emph{eur}s.}.). D'autre part, la théorie des fibrés principaux (c'est-à-dire du cas particulier des actions que nous appellerons \og principales\fg\,\footnote{C'est-à-dire, si l'on se restreint au point de vue purement ensembliste, abstraction faite de la différentiabilité, libres.}), telle que la concevait C. Ehresmann, conduit à la considération des \og groupoïdes de jauge \fg\,, qui furent son premier exemple de groupoïde de Lie (voir l'annexe \ref{nh}).

\smallskip

L'interaction des variétés et des graphes d'équivalence se manifeste par le \og thé\-orème de Godement \fg, lequel caractérise les relations d'équivalence dites régulières, précisément par la propriété que leur graphe est un groupoïde de Lie (dégénéré); ceci s'applique notamment aux feuilletages simples (classes connexes).

Les feuilletages réguliers non simples donnent naissance à un groupoïde de Lie non dégénéré, à savoir le groupoïde d'holonomie, qui est souvent appelé \og graphe \fg\, du feuilletage (voir l'annexe historique \ref{nh}), dont les groupes d'isotropie sont les groupes d'holonomie introduits antérieurement dans la thèse de G. Reeb (1952). Nous verrons que, plus généralement, les groupoïdes de Lie sont étroitement liés aux feuilletages avec singularités dits de Sussmann-\St.

\subsection{Orientation.}
Respectant ce qui nous a semblé être la vocation de ces \emph{Jour\-nées}, nous chercherons à suggérer quelques pistes de réflexion, tout en restant la plupart du temps largement informel; un traitement formel et détaillé entraînerait d'ailleurs beaucoup trop loin.

Notre point de vue est très voisin de celui exposé par A. Weinstein au cours de ces \emph{Journées}. Nous signalerons les différences de terminologie et de notations pour faciliter le passage d'un exposé à l'autre.

Dans une première section, nous esquissons un cadre plus général (\emph{\og diptyques \fg}) permettant de resituer les groupoïdes de Lie dans un contexte de grou\-poïdes structurés analogue à celui qu'avait développé C. Eheresmann, mais qui nous paraît plus propre à tenir compte de la spécificité de la catégorie des variétés tout en couvrant un large éventail d'autres exemples utiles\,\footnote{Il serait éminemment souhaitable de construire un cadre susceptible d'accueillir d'autres exemples intéressants (dont certains ont été évoqués au cours de ces \emph{Journées}, et pour lesquels des définitions spécifiques ont parfois été données) de groupoïdes porteurs d'une structure supplémentaire, \og compatible \fg\, en un sens à préciser, avec la structure algébrique, mais que l'on ne sait pas actuellement décrire comme groupoides structurés au sens d'Ehresmann, faute de savoir définir la bonne catégorie : groupoïdes symplectiques, de Poisson, de contact, de Jacobi, riemanniens, mesurés \ldots

Il se peut qu'il soit nécessaire d'introduire quelques retouches à nos axiomes, et notamment de considérer des catégories munies d'une involution (ici triviale).}.

Ce cadre fournit une méthode générale pour transférer de façon \og automatique \fg\, certaines constructions d'apparence purement ensembliste dans des catégories plus riches, et singulièrement celle des variétés, lorsque l'on a réussi à traduire ces constructions par des \emph{diagrammes commutatifs} appropriés, où l'on met en évidence \emph{injections} et \emph{surjections}, ainsi que certaines propriétés des \emph{carrés commutatifs}.

Cela nous permettra fréquemment dans la suite de nous limiter à une description de la structure des groupoïdes et des morphismes dans le cadre ensembliste, mais d'une manière permettant de suggérer ce qui s'étend aux groupoïdes différentiables, et plus gé\-né\-ra\-le\-ment \og structurés \fg\, par une catégorie convenable. Il faut d'ailleurs noter que beaucoup des propriétés que nous \og rappellerons \fg\, sont souvent mal connues, pour ne pas dire méconnues, même dans le cadre purement algébrique, bien que très élémentaires dans ce cas.

Nous pourrons alors donner une définition de l'\emph{équivalence de Morita dif\-fé\-ren\-tiable}, susceptible de nombreuses variantes, et de la structure (en un sens généralisé) d'un \emph{espace d'orbites} comme un classe d'équivalence (différentiable) de Morita de groupoïdes de Lie.

Nous introduirons enfin la notion de \emph{feuilletage singulier} au sens de \St\, sous une forme plus simple que la définition originale, puis étudierons la relation, à la fois étroite et riche, entre feuilletages (réguliers ou singuliers) et groupoïdes de Lie, selon ce qu'annonce le titre de l'exposé.

\section{Un cadre formel pour passer des ensembles aux variétés.}
\subsection{Motivations.}
La plupart des notions classiques de la théorie générale des catégories sont, à notre sens, mal adaptées à la catégorie des variétés, qui, de ce fait, est générale\-ment considérée comme une \og mauvaise \fg\, catégorie. Notamment:

-- elle contient des monomorphismes et épimorphismes pathologiques; inversement les plongements et les surmersions ne peuvent être caractérisés par des propriétés purement catégoriques;

-- elle n'est pas complète, et, notamment, ne possède pas suffisamment de produits fibrés et de quotients.

Il est souvent commode de la plonger dans une catégorie jugée \og meilleure \fg, par exemple la catégorie des espaces difféologiques de J.-M. Souriau, dans laquelle on a notamment des quotients quelconques; on a en particulier une difféologie canonique sur l'espace des feuilles d'un feuilletage.

Toutefois l'inconvénient d'une telle procédure est d'introduire une multitude d'objets auxiliaires hautement pathologiques jouant un rôle de parasites, et de conduire à des théorèmes qui ne distinguent pas toujours suffisamment les vraies variétés parmi de tels objets.

\smallskip
Nous suivrons ici une approche très différente, reposant sur quelques observations \og expérimentales \fg\, qui nous conduisent à considérer la catégorie des variétés comme au contraire \og excellente \fg\,:

-- les \emph{plongements} et les \emph{surmersions} constituent deux sous-catégories privilégiées possédant de remarquables propriétés de stabilité;

-- un nombre considérable de catégories d'usage courant dans des domaines très variés des Mathématiques sont munies de façon naturelle d'une ou plusieurs paires de sous-catégories possédant ces mêmes propriétés formelles;

-- cependant, dans une catégorie arbitraire, ces propriétés ne sont pas vérifiées dans leur totalité par les deux sous-catégories formées de \emph{tous} les monomorphismes et \emph{tous} les épimorphismes, pas plus, à notre connaissance, que par celles que l'on chercherait à définir en imposant à ces derniers des propriétés catégoriques sup\-plé\-men\-taires de type général, telles que d'être \og scindés \fg\,, \og stricts \fg\,, etc., de sorte qu'il s'agit bien d'une propriété des catégories envisagées, qui nous semble liée précisément à leur utilité;

-- bien qu'individuellement ces propriétés soient d'apparence anodine et qu'elles soient en petit nombre, il s'avère que leur conjonction est très puissante, et permet, comme on le verra, d'aller beaucoup plus loin que l'on ne pourrait s'y attendre \emph{a priori}.

\smallskip
Du fait que la plupart des définitions et propriétés envisagées ici (dans la mesure où nous ne considérons pas ici l'aspect infinitésimal) sont valables dans ce cadre, nous nous placerons dans ce contexte général, esquissé ci-après sous le nom de \og \emph{diptyques} \fg\,, pour notre définition générale des groupoïdes structurés.

\subsection{Diptyques.}\label{dip}\footnote{On trouvera des indications beaucoup plus détaillées dans l'exposé que nous avons donné à Krynica [P 2003].}

Une \emph{structure de diptyque} $\textbf{D}=(D;D_i,D_s)$ sur une catégorie $D$\,\footnote{On écrira souvent $D$ pour $\textbf{D}.$} consiste en les \emph{données} suivantes:

-- des \emph{produits} pour les familles d'objets finies non vides\,\footnote{Dans le cas fréquent où le produit de la famille vide existe aussi, c'est un objet final, noté $\bullet$.};

-- deux \emph{sous-catégories} $D_i/D_s$, appelées à jouer le rôle des injections/surjections dans \textbf{Ens}, ou des plongements/surmersions dans \textbf{Dif}\,\footnote{Lorsqu'il ne crée pas trop d'ambiguïté, le symbole / est utilisé d'une façon quelque peu laxiste mais suffisamment transparente, ici et dans la suite, pour éviter l'emploi lourd et répété de la formule \og respectivement\fg\,, plus correcte.}.

On notera $D_\ast$ le sous-groupoïde des \emph{isomorphismes} (flèches inversibles), et $D_r=D_iD_s$ la sous-classe des flèches dites \og \emph{régulières}\,\fg, c'est-à-dire décomposables selon $D_i$ et $D_s$. En vue des applications \emph{il est essentiel de ne pas imposer que $D_r$ soit une sous-catégorie, et \emph{a fortiori} que toute flèche soit régulière}.

Les flèches de $D_i/D_s$ seront appelées les \og \emph{bons monos/épis}\fg\, et notées en général avec une queue/tête triangulaire, pleine ou non, selon l'un des modèles 
\begin{center}
$\riTo$/$\rsTo$, ou $\rITo$/$\rSTo$.
\end{center}

Ces données sont assujetties à satisfaire à quelques \emph{axiomes de stabilité} d'apparence très simple, que, pour ne pas alourdir, nous préfèrerons énoncer dans l'\emph{annexe} B, en indiquant seulement ci-après (Tab. \ref{ed}) quelques \emph{exemples} fondamentaux (parmi une multitude d'autres).
\begin{table}[!hbt]
\caption{\label{ed}\large\scshape Quelques exemples de Diptyques}
\def\DD{
\begin{center}
\begin{diagram}[size=2.2em,tight,inline,abut]
{\cdot}
&\riTo&{\cdot}\\
\dTo&&\dTo\\
{\cdot}&\riTo&{\cdot}\\
\end{diagram}
\end{center}
}
\def\EE{
\begin{center}
\begin{diagram}[size=2.2em,tight,inline,abut]
{\cdot}&\rsTo&{\cdot}\\
\dTo&&\dTo\\
{\cdot}&\rsTo&{\cdot}\\
\end{diagram}
\end{center}
}
\def\GG{
\begin{center}
\begin{diagram}[size=2.2em,tight,inline,abut]
{\cdot}
&\rsTo&{\cdot}\\
\dTo&\pb&\dTo\\
{\cdot}&\rsTo&{\cdot}\\
\end{diagram}
\end{center}
}
\def\FF{
\begin{diagram}[size=1.1em,tight,inline]
\cdot&\rsTo&&&\cdot\\
&\cdot\rdsTo(1,1)[abut]&&\rusTo(3,1)[abut]&\\
\dTo\ldTo(1,3)[abut]&&\pb&&\dTo[abut]\\
&&&&\\
\cdot&&\rsTo[abut]&&\cdot\\
\end{diagram}
}
\begin{tabular}[]{||p{4em}||p{5em}|p{3em}|p{5.5em}|p{5em}|p{4.5em}||}
\hline
{}&Catégories\newline abéliennes&Topos&\bfseries Top\newline [TopHaus]&\bfseries Dif\newline \bfseries [DifHaus]&
\begin{center}
$\square D$
\end{center}
\\
\hline
\hline
\begin{center}
\scshape bons monos\newline$\riTo$
\end{center}
&
\begin{center}
tous
\end{center}
&
\begin{center}
tous
\end{center}
&
\begin{center}
plongements\newline [propres]
\end{center}
&
\begin{center}
idem $\ldots$
\end{center}
&\DD \\
\hline
\hline
\begin{center}
\scshape bons épis\newline$\rsTo$
\end{center}
&
\begin{center}
tous
\end{center}
&
\begin{center}
tous
\end{center}
&
\begin{center}
surjections\newline \emph{ouvertes}
\end{center}
&
\begin{center}
surmersions
\end{center}
&\EE \\
\cline{2-6}
{}&\multicolumn{2}{c|}{\scshape{variantes} :}&
\begin{center}
-- propres\newline -- étales\newline -- rétro-\newline connexes
\end{center}
&
\begin{center}
idem $\ldots$
\end{center}
&
\begin{center}
\GG ou\FF
\end{center}
\\
\hline
\end{tabular}
\end{table}

Ces axiomes permettent notamment de construire certains\,\footnote{De façon précise, les \emph{bons} \emph{carrés cartésiens} (ou \emph{bons pull backs}) sont ceux pour lesquels la flèche canonique $B'\riTo B\times A'$ est un bon mono. 

Dans \textbf{Dif}, la condition classique de transversalité est suffisante pour qu'il en soit ainsi, mais non nécessaire, et laisse échapper beaucoup de carrés cartésiens indispensables (notamment les intersections régulières de sous-variétés de basse dimension). En revanche il y a dans \textbf{Dif} des carrés cartésiens pathologiques, de sorte que, là aussi, nous pensons que la notion purement catégorique générale n'est pas la bonne.} carrés cartésiens (\emph{pull back squares}) (c'est-à-dire, rappelons-le, possédant la propriété universelle des produits fibrés), qui seront notés :

\begin{center}
\begin{diagram}[size=2em,tight,labelstyle=\scriptstyle,inline]
B'\SEpbk&\rTo\til{f'}&A'\\
\dTo\til v&&\dTo\til u\\
B&\rTo\til f&A
\end{diagram}
ou
\begin{diagram}[size=2em,tight,labelstyle=\scriptstyle,inline]
B'&\rTo\til {f'}&A'\\
\dTo\til v&\pb&\dTo\til u\\
B&\rTo\til f&A
\end{diagram}
\end{center}

Un rôle fondamental sera également joué dans la suite par les carrés \emph{pleins}, c'est-à-dire les carrés commutatifs $B'A'BA$ admettant une factorisation de la forme :
\begin{diagram}[w=2em,h=2em,tight,midshaft,labelstyle=\scriptstyle]
B'&\rTo\til{f'}&&&A'\\
&R\rdsTo(1,1)\til{r}[abut] \SEpbk&&\ruTo(3,1)\til{\pi_1}&\\
\dTo\til{v}\ldTo\til{\pi_2}(1,3)&&&&\dTo\til{u}\\
&&&&\\
B&&\rTo\til{f}&&A\\
\end{diagram}
où le carré $RA'BA$ est cartésien, et $r:B'\rsTo R$ un bon épi.

Dans le tableau \ref{ed} ci-joint annoncé, \textbf{Haus} signifie \emph{séparé}.

$\square D$ (notation d'Ehresmann) désigne la catégorie des \emph{carrés commutatifs} (quatuors d'Ehresmann) de $D$, munie de la loi de composition horizontale (par juxtaposition).

Une application continue est dite \emph{rétro-connexe} si l'image réciproque de tout point (donc de tout connexe) est connexe. (Observer que \og rétro-compacte \fg\, serait synonyme de \og propre \fg).

Parmi beaucoup d'autres exemples utiles non mentionnés dans le tableau, figurent notamment (avec des choix appropriés) les espaces de Banach, les fibrés vectoriels\ldots

Nous pensons que la prédilection des catégoriciens purs pour les catégories abé\-lien\-nes et les topos est liée au fait (plutôt rare dans les exemples concrets) que ces catégories partagent avec la catégorie \textbf{Ens} (non rappelée dans le tableau) la propriété que l'on peut y prendre pour bons monos/épis \emph{tous} les mono\-morphismes/ épi\-mor\-phismes.

\section{Groupoïdes algébriques.}
Comme annoncé, on considère dans cette section les groupoïdes d'un point de vue purement ensembliste ou algébrique, sans chercher à énoncer des définitions formelles complètes, mais en vue de préciser nos notations (très variables d'un auteur à l'autre), et en adoptant des points de vue et des notations qui facilitent le passage au cas structuré (notamment en mettant en évidence certaines injections et surjections, destinées à devenir de bons monos/épis dans le cas structuré).

\subsection{Notations.}\label{not}
Le schéma ci-après (fig. \ref{rsg}) est destiné à faciliter la mémorisation des axiomes des groupoïdes (ensemblistes), lesquels peuvent être brièvement définis comme les (petites\,\footnote{Cela veut dire que les flèches et les objets constituent des ensembles.}) catégories dont toutes les flèches sont inversibles. Les \emph{morphismes} de groupoïdes ne sont autres que les foncteurs entre groupoïdes.

\begin{figure}[!hbt]
\fbox{
\begin{diagram}[w=3.3em,h=1.7em,tight,abut,labelstyle=\scriptstyle,midshaft,thick]
{}&&&&\rLine&&&&{}&&\\
&&&&&&&&&&\\
&&&&\mathop{\bullet}^{yx}_{\swarrow\searrow}&&&&\dLine&&\\
&&&\mathop{\bullet}^y_{\swarrow\searrow}&&\mathop{\bullet}^x_{\swarrow\searrow}&&&G&&G=G^{(1)}\\
&&&&&&&&\dLine&&{}\\
\dLine&\rLine&\overset{\omega(\beta(y))}\bullet&\rLine&\overset{\omega(\alpha(y))=\omega(\beta(x))}\bullet&\rLine&\overset{\omega(\alpha(x))}\bullet&\rLine&\omega(B)&&\\
&&&&&&&&\dLine&&\\
{}&&&&\rLine&&&&{}&&\dsTo^\beta \uiTo\til{\omega}\dsTo_\alpha\\
&&&&&&&&&&\\
&&&{}_{\swarrow}&\lLine[thin]_{yx=\gamma(y,x)}&{}_\diagdown&&&&&\\
&&&\underset{y}\curvearrowleft&&\underset{x}\curvearrowleft&&&{}&&{}\\
{}&\rLine&\overset{\beta(y)}\bullet&\rLine&\overset{\alpha(y)=\beta(x)}\bullet&\rLine&\overset{\alpha(x)}\bullet&\rLine&B&&B=G^{(o)}\\
\end{diagram}}

\caption{\label{rsg}Représentation schématique d'un groupoïde.}
\end{figure}

Les \emph{données} d'une structure de groupoïde \textbf{G} sont donc :\\
-- l'ensemble des \emph{flèches}, noté $G$ ou $G^{(1)}$\,\footnote{Notation liée au point de vue simplicial, que nous ne développerons pas ici: voir [P 2003].};\\
-- l'ensemble des \emph{objets} ou \emph{base} $B$ ou $G^{(0)}$ (de façon laxiste, on écrira souvent $G$ ou $\gpd{G}{B}$ pour $\textbf{G}$);\\
-- la loi \emph{unité}, qui associe à tout objet une unité, et définit une injection $\omega=\omega_G:B\riTo G$\,; on \emph{identifiera} (le plus souvent implicitement) $B$ avec $\omega(B)$;\\
-- les surjections \emph{source} et \emph{but} $\alpha=\alpha_G,\ \beta=\beta_G:G\rsTo B$\,\footnote{Notations d'Ehresmann. On rencontre souvent la notation $s$ ou $d$ pour \emph{source} ou \emph{domain} et $r$ pour \emph{range} (but, nommé aussi \emph{codomain}). Nous écrirons de préférence la source à droite et le but à gauche; signalons que C.M. Marle adopte cette convention, mais intervertit les notations $\alpha,\ \beta$, ce qui rétablit l'ordre\ldots alphabétique. Signalons aussi la notation simpliciale $d_0,d_1$ pour $\beta,\ \alpha$.} (admettant $\omega$ pour section commune);\\
-- la loi de \emph{composition} $\gamma=\gamma_G:\Gamma G\rsTo G$ (qui se trouve être une surjection), définie sur le produit fibré $\Gamma G$ de $\alpha$ et $\beta$.

\smallskip
L'\emph{inversion} (ou \emph{symétrie}) n'est pas une donnée supplémentaire, puisque son existence et son unicité sont garanties par les axiomes; elle détermine une bijection, que l'on notera: $$\varsigma=\varsigma_G:G\rightarrow G,\ x\mapsto x^{-1}.$$

\smallskip
En fait il nous sera commode pour la suite de considérer la structure \textbf{G} comme définie plutôt par la donnée (équivalente) du quintuplet:
$$\textbf{G}=(G,B,\omega,\alpha,\delta)$$
où $\delta=\delta_G:(y,x)\mapsto yx^{-1}$ est la \emph{division} (ou différence); c'est une surjection $\delta_G:\Delta G\rsTo G$, définie sur le graphe $G^{(2)}=\,\wedge G=\Delta G$ de la relation d'équivalence sur $G$ définie par la surjection $\alpha$, autrement dit sur le produit fibré de la flèche $\alpha$ par elle-même.

De ces données, on déduit immédiatement $\varsigma(x)=x^{-1}=\delta(\omega(\alpha(x)),x))$, puis $\beta=\alpha\circ\varsigma= \delta\circ\text{diag}_G$, et $\gamma(y,x)=\delta(y,\varsigma(x))$.

Les notations $\wedge G=\Delta G$ évoquent le diagramme triangulaire suivant, dans 
\begin{diagram}[h=1.5em,w=1em,tight,objectstyle=\scriptstyle,labelstyle=\scriptstyle,textflow,abut]
&&a&&\\
&\ldTo(2,2)\til y&&\rdTo(2,2)\til x&\\
c&&\lTo\til{yx^{-1}}&&b\\
\end{diagram}
lequel on a $\alpha(x)=\alpha(y)=a,\ \beta(x)=b,\ \beta(y)=c$. On voit que 
l'application $\delta$ permet de compléter le triangle commutatif dans \textbf{G}.

Sur la figure \ref{rsg} ci-jointe, qui schématise la loi de composition $\gamma$, on voit qu'il est commode de se représenter les éléments $x,\,y$ de $G$, tantôt comme des points de $G$ munis de deux projections sur $B$ ou $\omega(B)$, tantôt comme des flèches reliant deux points de la base $B$.

\smallskip
 Pour tout objet $a\in B$, nous noterons:
\begin{itemize}
	\item ${G_a}=\alpha^{-1}(a)$ la $\alpha$-fibre de $a$;
	\item ${_aG}=\beta^{-1}(a)$ sa $\beta$-fibre;
	\item ${_a{G}_a}={_aG}\cap{G_a}$ le \emph{groupe d'isotropie}\,\footnote{Groupe de cohérence pour W.T. van Est.} de $a$.
\end{itemize}

\subsection{Le transiteur (anchor map).}\label{am}
Nous poserons: $$\tau_G=\tau=(\beta,\alpha):G\rightarrow B\times B.$$ C'est un morphisme (foncteur) lorsqu'on regarde $B\times B$ comme un groupoïde \emph{banal}, graphe de la relation d'équivalence grossière; il traduit les propriétés d'(in)\-tran\-si\-ti\-vité du groupoïde \textbf{G}. Conformément à nos principes généraux de terminologie (usage du suffixe \emph{eur}), nous l'appellerons le \emph{transiteur} de \textbf{G}. La terminologie usuelle est \emph{ancre}\,\footnote{Peu parlante pour les peuples sans tradition maritime suffisante.}, traduction de l'expression \emph{anchor map}, introduite par K. Mackenzie [McK]. 

Les propriétés d'injectivité ou de surjectivité de $\tau_G$ jouent, nous le verrons, un rôle crucial, encore plus apparent dans le cas structuré, du fait qu'il existe alors de multiples variantes et raffinements de ces deux notions.

\bigskip
Du point de vue \emph{ensembliste}, $\tau_G$ admet toujours une factorisation en surjection et injection\,\footnote{La situation sera beaucoup plus riche et complexe dans le cas structuré. Les pointillés dans le diagramme attirent l'attention sur les flèches qui n'existeront pas toujours dans le cas structuré.}:
\begin{diagram}[h=1em,w=3em,labelstyle=\scriptstyle]
&{}&\rLine[abut]\til{\tau_G}&{}&\\
G\ruLine(1,1)[abut]&\rDotsTo\til s&R&\rDotiTo\til i&B\times B\rdTo(1,1)[abut]\\
\end{diagram}

$R$ est le \emph{graphe d'une relation d'équivalence} (qui sera dite \emph{orbitale}) sur la base $B$, dont les classes sont les composantes transitives ou \emph{orbites} du groupoïde.

Le \emph{noyau}\,\footnote{Le noyau d'un morphisme de groupoïdes est évidemment défini comme image réciproque de la base ou de l'ensemble des unités.} $K$ de $\tau_g$ sera appelé \emph{plurigroupe d'isotropie}: c'est la somme des groupes d'isotropie.

La situation est résumée par le diagramme suivant, dans lequel le carré du haut 
\begin{diagram}[w=3.5em,h=2.5em,textflow]
K\pbk&\rDotsTo\til k&B\\
\dDotiTo\til j&&\diTo\til{\text{diag}_B}\\
G&\rTo\til{\tau_G}&B\times B\\
\dDotsTo\til s&&\dEqual\\
R&\rDotiTo\til i&B\times B
\end{diagram}
est cartésien. Les deux bords verticaux sont des exemples (très dégénérés) de ce qu'il convient d'appeler des \emph{suites exactes} (courtes) de groupoïdes\,\footnote{La notion générale résulte de celle d'extenseurs, introduite plus loin.}. La suite de gauche fait apparaître tout groupoïde (considéré d'un point de vue purement ensembliste et algébrique) $\textbf{G}$ comme \emph{extension du graphe de l'équivalence orbitale $R$ (\og partie horizontale \fg) par le plurigroupe d'isotropie $K$ (\og partie verticale \fg)}.

\subsection{Dégénérescences.}\label{degen}
Le diagramme précédent est susceptible de différents types de \emph{dégénérescence}:
\begin{itemize}
	\item les conditions suivantes sont équivalentes:
\begin{itemize}
	\item $\tau_G$ est surjective;
	\item $i$ est bijective;\\
	$G$ est dit \emph{transitif}\,\footnote{Cohérent dans la terminologie de W.T. van Est.};
\end{itemize}
	\item les conditions suivantes sont équivalentes:
\begin{itemize}
	\item $\tau_G$ est injective;
	\item $s$ est bijective;
	\item $k$ est bijective;\\
	$G$ est alors (isomorphe à) un graphe d'équivalence (\emph{dégénérescence horizontale}); on dira que $G$ est \emph{principal}, ou \emph{de Godement};
\end{itemize}
\item $\tau_G$ est bijective: on dira que $G$ est \emph{banal}\,\footnote{La terminologie usuelle \og grossier \fg\, est à exclure, car elle devient trop ambiguë pour les groupoïdes topologiques.};
	\item les conditions suivantes sont équivalentes:
\begin{itemize}
	\item $\alpha=\beta$;
	\item $j$ est bijective;\\
	on dira que $G$ est un \emph{plurigroupe}\,\footnote{Une terminologie telle que \og groupe multiple \fg\, entrerait en conflit avec la notion de catégorie n-uple. Signalons que Douady-Lazard ont utilisé le terme de \og groupe \fg\, dans le sens où nous dirons \og plurigroupe \fg.}: c'est une somme de groupes, qui se réduit à un \emph{groupe} si de plus $B$ est un singleton (\emph{dégénérescence verticale});
\end{itemize}
\item $\omega$ est bijective: $G$ sera dit \emph{nul}\,\footnote{La terminologie répandue \og discret \fg\, est à rejeter pour la même raison que \og grossier \fg.}.
\end{itemize}

\subsection{Décomposition orbitale ensembliste.}\label{doe}
\begin{figure}[!hbt]
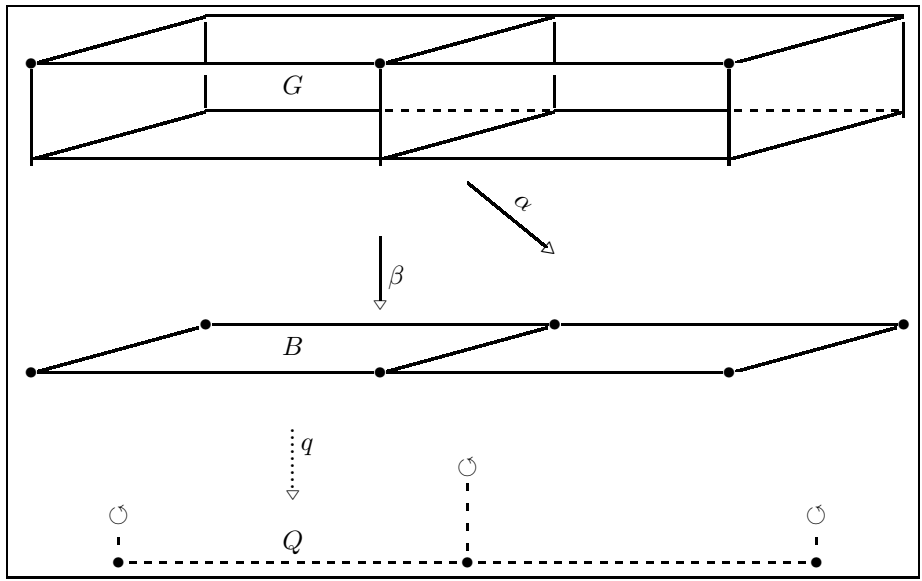

\fbox{
\begin{diagram}[w=3.3em,h=.9em,tight,labelstyle=\textstyle,midshaft,thick]
&&{}&\rLine[abut]&&&{}&&\rLine[abut]&&{}\\
&\ldLine[abut]&\dLine[abut]&&&\ldLine[abut]&\dLine[abut]&&&\ldLine[abut]&\dLine[abut]\\
\bullet&\rLine[abut]&\HonV&&\bullet&\rLine[abut]&\HonV&&\bullet&&\\
\dLine[abut]&&\dLine&G&\dLine[abut]&&\dLine&&\dLine[abut]&&\\
&&{}&\rLine[abut]&\VonH&\rDashLine[abut]&{}&\rLine[abut]&\VonH&\rDashLine[abut]&{}\\
&\ldLine[abut]&&&&\ldLine[abut]&&&&\ldLine[abut]&\\
{}&\rLine[abut]&&&{}&\rLine[abut]&&&{}&&\\
&&&&&{}&&&&&\\
&&&&&&\rdsTo[abut](1,3)^\alpha&&&&\\
&&&&{}&&&&&&\\
&&&&&&{}&&&&\\
&&&&\dsTo[abut]_\beta&&&&&&\\
&&&&{}&&&&&&{}\\
&&{\bullet}&\rLine[abut]&&&{\bullet}&&\rLine[abut]&&{\bullet}\\
&\ldLine[abut]&&B&&\ldLine[abut]&&&&\ldLine[abut]&\\
{\bullet}&\rLine[abut]&&&{\bullet}&\rLine[abut]&&&{\bullet}&&\\
&&&&&&&&&&\\
&&&{}&&&&&&\\
&&&\dDotsTo[abut]_q&&&&&&&\\
&&&&&{\circlearrowleft}&&&&&\\
&&&{}&&&&&&&\\
&{\circlearrowleft}&&&&\dDashLine[abut]&&&&{\circlearrowleft}&\\
&\dDashLine[abut]&&&&&&&&\dDashLine[abut]&\\
&\bullet&&\rDashLine[abut]^Q&&\bullet&\rDashLine[abut]&&&\bullet&\\
\end{diagram}}
\caption{\label{do}Décomposition orbitale d'un groupoïde.}
\end{figure}
\bigskip

Toujours d'un point de vue purement ensembliste, on peut introduire \emph{l'ensemble orbital} $Q$, quotient de la base $B$ par l'\emph{équivalence orbitale} définie par le graphe $R$. On notera 
$$q:B\rDotsTo Q$$
la projection canonique (de nouveau les pointillés soulignent le fait qu'en général $q$ n'existera pas dans le cas structuré).

La situation est résumée par le diagramme suivant. On note que le carré cartésien 
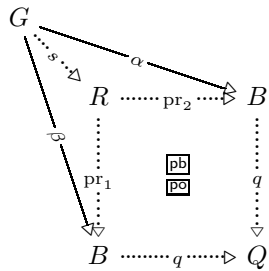
\begin{diagram}[size=3em,tight,labelstyle=\scriptstyle,textflow]
G&&&\\
&R\rdDotsTo(1,1)\til s\rdsTo(1,3)\til{\beta}\rdsTo(3,1)\til{\alpha}&\rDotsTo\til{\text{pr}_2}&B\\
&\dDotsTo\til{\text{pr}_1}&\genfrac{}{}{0pt}{1}{\pb}{\po}&\dDotsTo\til q\\
&B&\rDotsTo\til q&Q
\end{diagram}
est aussi \emph{cocartésien} (\emph{pushout square}); autrement dit la flèche $q$ peut aussi être vue comme le co-égalisateur [McL] des deux projections du graphe $R$. C'est aussi le co-égalisateur de $\alpha$ et $\beta$; le carré $GBBQ$ est aussi cocartésien.

On a aussi une surjection $r:R\rDotsTo Q$, que l'on peut même considérer comme un morphisme de groupoïdes nul, \cad ayant pour but un groupoïde nul.

La flèche composée $r\circ s:G\rDotsTo Q$ détermine la \emph{décomposition de $G$ en somme de ses composantes transitives}, que nous appellerons \emph{décomposition orbitale du groupoïde $G$}.\footnote{Bien qu'en général, lorsque $G$ sera un groupoïde de Lie, $Q$ ne soit pas défini en tant que variété, nous verrons que la décomposition orbitale détermine un feuilletage au sens de \St.}.

La structure algébrique d'un groupoïde se ramène donc à l'étude d'une composante transitive, laquelle fait l'objet du paragraphe suivant. On verra qu'on peut lui attacher un groupe, défini à isomorphisme près.

La situation est schématisée par la figure \ref{do} ci-jointe, où les classes d'équivalence sont dans des plans verticaux.

L'espace orbital est dessiné en tirets, pour souligner son caractère \og virtuel \fg \,(dans le cas structuré). Les lignes verticales en tirets surmontées d'un flèche fermée sont censées évoquer les groupes attachés à ses points, par lesquels ces derniers se trouvent munis d'une sorte de structure, en un sens généralisé, qui est l'ombre (platonicienne ?) ou le reflet de la situation décrite en traits pleins aux étages supérieurs.

\subsection{Trivialisation ensembliste des composantes transitives.}\label{tect}
Nous nous inté\-ressons donc ici à la structure d'un groupoïde \emph{transitif}.

\begin{figure}[!hbt]
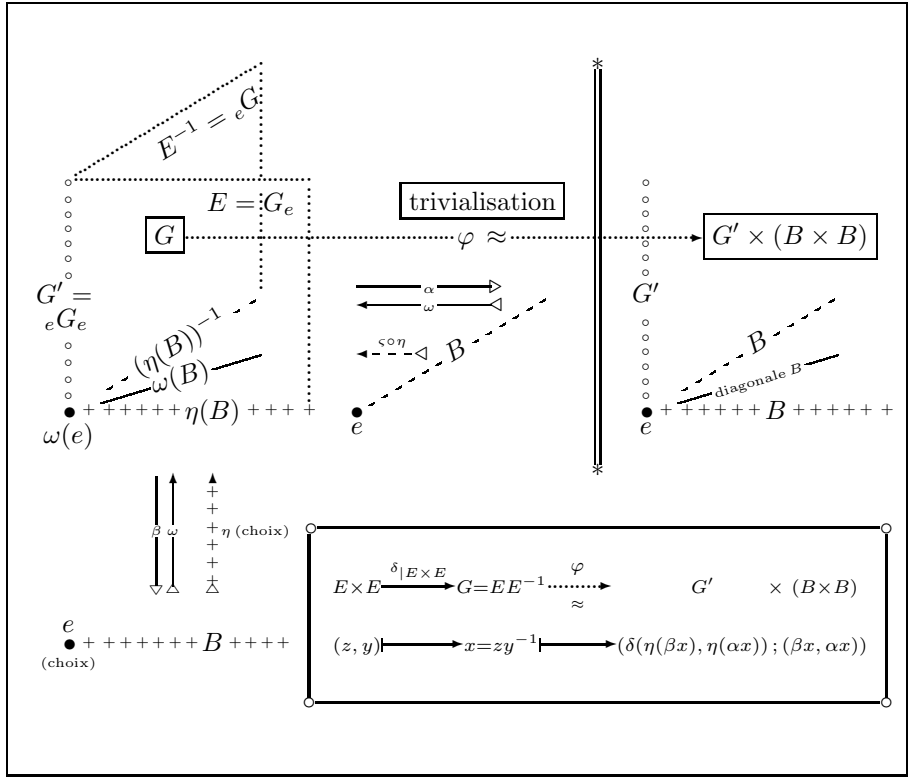

\fbox{
\def\diagG'{\begin{diagram}[size=.9em,tight,center]
\drnd\\G'=\,\,\\{_eG_e\,}
\end{diagram}}
\begin{diagram}[w=1.82em,h=2.2em,p=.2em,tight]
&&&&&&&&&&&&&&&&&&\\
&&&&{\cdot}&&&&&&&{\ast}&&&&&&&\\ 
&&&\ruDotLine(4,2)[abut]_{\quad\qquad E^{-1}={_eG}\,\,}&\dDotLine[abut]&&&&&&&& &&&&&\\
\urnd&&\rDotLine[abut]_{\,\,\,\,\qquad\qquad E=G_e}&&\HonV&{\cdot}&&&&&&&\urnd&&&&&\\
&&\boxed G&\rDotLine[abut]&\dDotLine[thin]&&&&&&\rDotTo[abut,midshaft]^{\boxed{\text{trivialisation}}}\til{\varphi\,\approx}& &&&&\boxed{G'\times(B\times B)}&&\\
\dRndLine[abut]\til{\diagG'}&&&&{}&&{}&\pile{\rsTo[abut,thick,midshaft,\lab=\ssst]\til{\alpha}\\ \liTo[abut,thin,midshaft,\lab=\ssst]\til{\omega}}&&{}&{}&\dLine[abut]\dLine[abut]& \dRndLine[abut]\til{\overset{}{G'}}&&&&{}&\\
&&&\ldDashLine(4,2)[abut,midshaft]\til{(\eta(B))^{-1}}&{}&\dDotLine[abut]&{}&\lDashiTo[abut,\lab=\ssst]^{\varsigma\circ\eta}&{\qquad}&\ruDashLine(4,2)[abut]\til{B}&&&&&&\ldDashLine(4,2)[abut]\til{B}&{}&\\
\overset{\drnd}{\underset{\textstyle{\omega(e)}}{\quad\bullet\rsum}}&&&\rSumLine[abut]\til{\eta(B)}\ruLine[abut,thin](4,1)\til{\omega(B)}&&\rsum&\underset{\textstyle{e}}\bullet&&&&&&
\overset{\drnd}{\underset{\textstyle{e}}{\quad\bullet\rsum}}&&&\ruLine[abut,\lab=\ssst](4,1)\til{\text{diagonale}\,B}\rSumLine[abut,midshaft]\til{B}&&\rsum\\
&&{}&{}&&&&&&&&{{\ast}}&&&&&&\\
&&\dsTo[abut,thick,\lab=\ssst]\til{\ssst\beta}\uiTo[abut,thin,\lab=\ssst]\til{\omega}&\uSumiTo[abut,\lab=\ssst]_{\eta\,(\text{choix})}&&\,{\circ}&&&&&&\rLine[thick,abut]&&&&&&{\circ}\\
&&{}&{}&&&\sst{E\times E}&\rTo[abut,\lab=\ssst]^{\delta_{|E\times E}}&&\sst{G=EE^{-1}}&&\rDotTo[\lab=\ssst,midshaft,abut]^{\sst\varphi\,}_{\approx\,}&&&\quad\qquad\sst{G'\qquad\times\,\,(B\times B)\,\,\,\,}&&&\\
\overset{\textstyle{e}}{\underset{\ssst{(\text{choix})}}{\quad\bullet\rsum}}&&&\rSumLine[abut]\til{B}&&{\,}\dLine[thick,abut]^{\,}&\sst{(z,\,y)}&&\rMapsto[abut]&\sst{x=zy^{-1}}&&\rMapsto[abut]&&&\sst{(\delta(\eta(\beta x),\,\eta(\alpha x))\,;\,(\beta x,\,\alpha x))}&&&{\,}\dLine[thick,abut]\\
&&&&&{\circ}&&&&&\rLine[thick,abut]&&&&&&&{\circ}\\
&&&&&&&&&&&&&&&&&&\\
\end{diagram}}
\caption{\label{tct}Trivialisation \emph{ensembliste} d'une composante transitive.}
\end{figure}

Ici encore la construction sera schématisée par une figure (fig\ref{tct}).

Sur cette figure on a \emph{dissocié deux copies} de la base $B$ du groupoïde $G$, figurées, l'une avec des $+$, l'autre avec des $-$, considérées, l'une comme le but de $\beta$, l'autre de $\alpha$, ces deux projections étant munies chacune de la section commune $\omega$.

La construction repose sur le \emph{choix} d'un objet $e\in B$. On pose $E=G_e$, $E^{-1}=\varsigma(E)={_eG}$, et l'on note $G'={_eG_e}$ le \emph{groupe d'isotropie} de $e$ (schématisé sur la figure par une ligne verticale formée de $\circ$). On peut observer que l'on a $G=\delta(E,E)=EE^{-1}$\,\footnote{Notation fréquemment utilisée par C. Ehresmann dans sa construction du groupoïde de jauge d'un fibré principal (voir ci-après).}.

\emph{Par l'hypothèse de transitivité}, la restriction de $\beta$ à $E$ est encore \emph{surjective}, ce qui permet d'en \emph{choisir} une \emph{section} $\eta$ (ceci n'étant en général plus possible dans le cas structuré). (On peut noter que la restriction de $\alpha$ à $E^{-1}$ admet $\varsigma\circ\eta$ pour section.)

Alors on vérifie facilement que les formules encadrées dans la figure \ref{tct} définissent un isomorphisme (algébrique) $\varphi$ entre $G$ et le groupoïde (\og \emph{trivial} \fg) produit du groupe $G'$ par le groupoïde banal $B\times B$ (trivialisation).

En composant $\varphi$ avec la projection canonique sur $G'$, on définit une \emph{équivalence de groupoïdes} (au sens algébrique [McL]) surjective entre $G$ et le groupe $G'$.
\bigskip

\emph{Remarque}: En raison de ce qui précède, la plupart des catégoriciens purs consi\-dèrent les groupoïdes comme des objets inintéressants, puisque équivalents à des groupes. Il s'agit là d'un double erreur d'appréciation:
\begin{itemize}
	\item d'une part on ne doit pas confondre \og trivialisable \fg\, et \og trivial \fg\, : nous verrons que, même d'un point de vue purement algébrique, l'étude des morphismes de groupoïdes (étrangement mal connue, bien que très élémentaire d'un point de vue purement algébrique) ne se réduit nullement à celle des morphismes de groupes\,\footnote{Ceci ne doit pas surprendre: ce n'est pas parce que les espaces vectoriels sont caractérisés à isomorphisme près par le cardinal de leur dimension que la théorie des application linéaires se réduirait à celle des cardinaux !};
	\item d'autre part la trivialisation précédente n'est plus (en général) possible dans le cadre différentiable ou structuré; il conviendra de remplacer l'équivalence précédente par une forme convenable de l'équivalence de Morita.
\end{itemize}

\smallskip
En fait, comme nous le verrons, la théorie des groupoïdes \emph{transitifs} dans le cadre différentiable est essentiellement équivalente à celle des \emph{fibrés principaux}, telle que l'avait conçue C.\til Ehres\-mann (voir l'annexe A).

\emph{Ceci permet de concevoir un groupoïde de Lie général comme décrivant une variation différentiable d'un fibré principal, laquelle autorise une variation discontinue du groupe structural}.

Ceci sera illustré par le cas particulier des feuilletages: la considération du \emph{groupoïde d'holonomie} décrit de façon différentiable la variation du \emph{revêtement d'holonomie} des feuilles.

\section{Groupoïdes D-structurés.}
On se place désormais dans un diptyque fixé \textbf{D}, l'exemple fondamental pour notre présent propos restant le diptyque $\textbf{D}=\textbf{Dif}$\,\footnote{De préférence ici à \textbf{DifHaus}, du moins pour l'étude des feuilletages, car celle-ci exige de pouvoir admettre des variétés non séparées et des plongement non fermés.}.

\subsection{Groupoïdes dans un diptyque.}
La présentation donnée ci-dessus des grou\-poïdes (algébriques) comme quintuplets $\textbf{G}=(G,B,\omega,\alpha,\delta)$, se laisse aisément transférer dans le cadre des diptyques, et aboutit à la définition des groupoïdes \textbf {D}-structurés, ou brièvement \textbf{D}-groupoïdes. Les groupoïdes de Lie sont les \textbf{Dif}-groupoïdes\,\footnote{\'A chaque variante proposée dans le tableau \ref{ed} correspondra une variante intéressante de cette notion, que nous n'examinerons pas ici.}.

On trouvera beaucoup de précisions et d'autres développements dans notre exposé [P 2003].

Disons simplement que les termes du quintuplet devront être des objets ou flèches de $D$, avec la \emph{condition supplémentaire} (annoncée par nos notations) que $\alpha$ soit un \emph{bon épi}, ce qui assure l'existence d'un (bon) produit fibré comme source de $\delta$: 
\begin{diagram}[size=3em,tight,labelstyle=\scriptstyle]
\Delta G&\riTo &G\times G&\pile{\rsTo^{\alpha\ \ }\\ \rsTo_{\alpha\ \ }}&B.\\
\end{diagram}

(Noter que les propriétés que $\omega$ soit un bon mono et $ \delta$ un bon épi n'ont pas à être imposées, car ce sont en réalité des conséquences faciles des axiomes des groupoïdes et des diptyques).

Nous allons maintenant examiner ce qui subsiste, dans le cadre des diptyques, des constructions précédentes, liées aux propriétés d'injectivité/surjectivité de $\tau_G$. Ceci nous conduira à une description diagrammatique de la théorie des fibrés différentiables principaux selon Ehresmann, laquelle étend cette théorie lorsqu'on remplace la catégorie \textbf{Dif} par un diptyque quelconque.

\subsection{D-groupoïdes réguliers.}
Les considérations des paragraphes \ref{am} et \ref{doe} s'é\-ten\-dent essentiellement au cas structuré dans le cas particulier important où la flèche $\tau_G$ est supposée \emph{régulière} (\ref{dip}), \cad où l'on a une factorisation:
\begin{diagram}[h=1em,w=3em,labelstyle=\scriptstyle]
&{}&\rLine[abut]\til{\tau_G}&{}&\\
G\ruLine(1,1)[abut]&\rsTo\til s&R&\riTo\til i&B\times {B.}\rdTo(1,1)[abut]\\
\end{diagram}

On dira que le \textbf{D}-groupoïde \textbf{G} est régulier.

\smallskip
Cette situation contient les deux cas particuliers suivants.

\subsubsection{Groupoïdes principaux et de Godement.}\label{god}

Nous dirons que \textbf{G} est:
\begin{itemize}
	\item \emph{de Godement} si $\tau_G:G\riTo B\times B$ est un bon mono ($s$ est un iso);
	\item \emph{principal} s'il existe un bon épi $q:B\rsTo Q$ tel que le carré suivant 
	\begin{diagram}[size=2.5em,tight,labelstyle=\scriptstyle,textflow]
G&\rsTo\til\alpha&B\\
\dsTo\til\beta&\genfrac{}{}{0pt}{1}{\pb}{\po}&\dsTo\til q\\
B&\rsTo\til q&Q
\end{diagram}
soit cartésien (donc parfait (\ref{qca}), \cad également cocartésien).
\end{itemize}

On a: principal $\Rightarrow$ de Godement.

\smallskip
Le \emph{théorème de Godement} (caractérisation des équivalences régulières par leur graphe), tel qu'il figure dans J.-P. Serre [LALG] (voir aussi N. Bourbaki [VAR]), équivaut à la \emph{réciproque}, lorsque $\textbf{D}=\textbf{Dif}$.

Le diptyque \textbf{D} sera dit \emph{de Godement} si tout groupoïde de Godement est principal.

Cette condition est vérifiée dans tous les exemples du tableau \ref{ed}.

\emph{Nous supposerons désormais, chaque fois que cela sera utile, que} \emph{D}\emph{est un diptyque de Godement.}

Alors, pour tout \textbf{D}-groupoïde régulier, on pourra écrire dans $D$ le diagramme \begin{diagram}[size=2.5em,tight,labelstyle=\scriptstyle,textflow]
G&&&\\
&R\rdsTo(1,1)\til s\rdsTo(1,3)\til{\beta}\rdsTo(3,1)\til{\alpha}&\rsTo\til{\text{pr}_2}&B\\
&\dsTo\til{\text{pr}_1}&\genfrac{}{}{0pt}{1}{\pb}{\po}&\dsTo\til q\\
&B&\rsTo\til q&Q
\end{diagram}
rencontré en \ref{doe} dans le cadre ensembliste. Le carré $GBBQ$ sera aussi cocartésien.

\subsubsection{Groupoïdes $s$-transitifs.}
C'est le cas \emph{symétrique du précédent}.

Nous dirons que \textbf{G} est \emph{$s$-transitif}\,\footnote{Dans cette terminologie la lettre $s$ se réfère à la donnée de diptyque $D_s$. Signalons (voir l'annexe historique A) que, pendant plusieurs décennies, le terme de \emph{groupoïde de Lie} a été réservé aux \textbf{Dif}-groupoïdes $s$-transitifs.} lorsque le transiteur $\tau_G:G\rsTo B\times B$ est un bon épi.

Ici la construction de \ref{tect} ne s'étend pas en général, sauf dans certains diptyques où les bons épis sont scindés (\emph{split}), \cad admettent une section. Nous reviendrons sur cette étude plus loin.

\section{Morphismes de D-groupoïdes.}
\subsection{Notations.}\label{notm}
On considère ici un deuxième \textbf{D}-groupoïde \textbf{H} de base $E=H^{(0)}$.

\smallskip
Un \textbf{D}-\emph{morphisme} (ou foncteur) $\textbf{f}:\textbf{H}\rTo \textbf{G}$ est défini par un couple de flèches de $D$, soit $f:H\rTo G,\ f^{(0)}=p:E\rTo B$ telles que: 
\begin{enumerate}
	\item $f$ et $f^{(0)}$ commutent avec $\omega_H$ et $\omega_G$ ainsi qu'avec $\alpha_H$ et $\alpha_G$, ce qui, par la propriété universelle des produits fibrés, détermine une flèche $\Delta f:\Delta H\rTo\Delta G$;
	\item $\Delta f$ et $f$ commutent avec $\delta_H$ et $\delta_G$.
	\end{enumerate}

On notera \textbf{Grpd(D)} la catégorie des \textbf{D}-morphismes de \textbf{D}-groupoïdes.

La définition est résumée par les carrés commutatifs de $D$ suivants:

\begin{center}
\begin{diagram}[size=2.5em,tight,inline,labelstyle=\scriptstyle]
E&\rTo\til p&B\\
\diTo\til{\omega_H}&\text{O}(f)&\diTo\til{\omega_G}\\
H&\rTo\til{f}&G\\
\end{diagram}
\quad
\begin{diagram}[size=2.5em,tight,inline,labelstyle=\scriptstyle]
H&\rTo\til{f}&G\\
\dsTo\til{\alpha_H}&\text{A}(f)&\dsTo\til{\alpha_G}\\
E&\rTo\til p&B
\end{diagram}
\quad
\begin{diagram}[h=2.5em,w=3.5em,tight,inline,labelstyle=\scriptstyle]
\Delta H&\rTo\til {\Delta f}&\Delta G\\
\dsTo\til{\delta_H}&\text{D}(f)&\dsTo\til{\delta_G}\\
H&\rTo\til f&G
\end{diagram}
\quad.
\end{center}

On aura alors aussi le carré commutatif:

\begin{center}
\begin{diagram}[h=2.5em,w=4em,tight,inline,labelstyle=\scriptstyle]
H&\rTo\til{f}&G\\
\dTo\til{\tau_H}&\text{T}(f)&\dTo\til{\tau_G}\\
E\times E&\rTo\til {p\times p}&B\times B
\end{diagram}
\quad.
\end{center}

Notons qu'à l'exception du carré \text{A}(f), ces carrés peuvent se lire aussi dans \textbf{Grpd(D)}, en considérant des groupoïdes nuls ou banaux.

\smallskip
On dira que \textbf{f} est un $i$-foncteur/un $s$-foncteur, lorsque $f$ (et, par suite, comme on peut le montrer, $f^{(0)}$ et $\Delta f$) sont de bons monos/épis\,\footnote{Contrairement à ce que l'on pourrait penser naïvement de prime abord, ce ne sont pas les bons candidats pour munir la catégorie \textbf{Grpd(D)} d'une structure (attendue!) de diptyque. Il y a en fait \emph{plusieurs} façons très utiles de munir \textbf{Grpd(D)} d'une structure de diptyque (de Godement), mais il convient pour ce faire d'imposer aux morphismes, candidats pour être de bons épis, diverses conditions supplémentaires, empruntées aux paragraphes suivants; ceci nous entraînerait trop loin de notre propos.}.

\smallskip
Il se trouve que les \emph{propriétés algébriques} importantes d'un morphisme de grou\-poïdes \textbf{f} sont reflétées par celles des \emph{carrés commutatifs} $\text{A}(f)$ (ou $\text{D}(f)$) et $\text{T}(f)$. L'intérêt de cette formulation diagrammatique est que, lorsqu'on la transfère au cadre structuré, conformément à notre philosophie générale, elle dit bien davantage que le contenu ensembliste sous-jacent (que nous laisserons le plus souvent au lecteur le soin d'énoncer), et conduit \og automatiquement \fg\, aux bonnes notions.

\bigskip Le morphisme \textbf{f} sera dit \emph{principal} si sa source est un groupoïde principal, \emph{nul} si sa source et/ou son but est nul (\ref{degen}).

\smallskip
La \emph{symétrie} détermine un \emph{isomorphisme} entre le groupoïde \textbf{G} et son \emph{opposé}:
$$\varsigma_G:\mathbf{G}\stackrel{\approx}{\rightarrow}\mathbf{G^{op}}.$$
Pour tout morphisme $\textbf{f}:\textbf{H}\rTo \textbf{G}$, on a le carré commutatif:

\begin{center}
\begin{diagram}[size=3em,tight,labelstyle=\scriptstyle,inline]
H&\rTo\til{f}&G\\
\dsTo\til{\varsigma_H}&\text{S}(f)&\dsTo\til{\varsigma_G}\\
H^{\text{op}}&\rTo\til{f^{\text{op}}}&G^{\text{op}}\\
\end{diagram}\quad.
\end{center}
(Noter que les morphismes $f$ et $f^{\text{op}}$ ont la \emph{même} flèche sous-jacente dans $D$.)

\smallskip
Ainsi que nos notations le soulignent, les carrés commutatifs que l'on vient d'écrire décrivent des \emph{transformations naturelles} [McL] \emph{entre endofoncteurs} de la catégorie \textbf{Grpd(D)} (lesquelles associent, à tout morphisme de groupoïdes, un carré commutatif de groupoïdes):
$$\omega:\Omega\rTo\til\tnat\text{Id},\quad\delta:\Delta\rTo\til\tnat\text{Id},\quad\tau:\text{Id}\rTo\til\tnat\Theta,\quad\varsigma:\text{Id}\rTo\til\tnat (.^{\text{op}})\,,$$
où $\Omega$ est le foncteur \og \emph{base} \fg\, (qui, à tout groupoïde, associe le groupoïde \emph{nul} déterminé par sa base), et $\Theta$ est le foncteur \og \emph{banal} \fg\, (qui, à tout groupoïde, associe le groupoïde banal déterminé par sa base). On note que $\varsigma$ est involutive.

\bigskip
\emph{On va s'intéresser successivement aux propriétés de} $\text{T}(f)$, \emph{puis de} $\text{A}(f)$

\subsection{\'Equivalences, inducteurs.}\emph{On considère ici le carré} $\text{T}(f)$. Il est facile de voir que, d'un point de vue purement algébrique, il reflète les propriétés de \emph{fidélité/ plénitude} du foncteur \textbf{f} [McL], lesquelles se lisent sur les propriétés d'injectivité/ surjectivité de la flèche canonique vers le produit fibré (ensembliste).

Pour simplifier, nous nous limiterons d'abord au cas où $\textbf{f}:\textbf{H}\rsTo\textbf{G}$ est un $s$-foncteur ce qui assure l'existence d'un (bon) produit fibré.

On dira que \textbf{f} est une $s$-équivalence/un $s$-extenseur si le carré \text{T}(f) est bien car\-tésien/ plein (\ref{dip}), ce que résume le tableau ci-dessous.

\begin{center}
\begin{tabular}[]{|c|c|}
\hline
\textbf{T(f)}&f\\
\hline
p.b.&$s$-équivalence\\
\hline
plein&$s$-extenseur\\
\hline
\end{tabular}
\end{center}

\subsubsection{$s$-équivalences, inducteurs.}\label{seq}
La définition des $s$-équivalences se traduit donc par le (bon) carré cartésien (dans $D$, mais aussi dans \textbf{Grpd(D)}) ci-dessous.

\begin{diagram}[h=3em,w=4em,tight,labelstyle=\scriptstyle]
H\SEpbk&\rsTo^f_\sim&G\\
\dTo\til{\tau_H}&&\dTo\til{\tau_G}\\
E\times E&\rsTo_{\sim}^{p\times p}&B\times B\\
\end{diagram}

On écrira $\textbf{H}=p^{\ast\ast}(\textbf{G})$, et on dira que \textbf{H} est \emph{induit} (\emph{pulled back}) par \textbf{G} le long de $p:E\rsTo B$.

On conserve cette terminologie quand on omet l'hypothèse $(E\rightarrow B)\in D_s$, mais que l'on exige seulement que le carré $\text{T}(f)$ soit bien cartésien (voir note de bas de page de \ref{dip}), et l'on dira alors que \textbf{f} est un \emph{inducteur} (la condition algébrique sous-jacente est que \textbf{f} soit \emph{pleinement fidèle}).

\smallskip
Un exemple très dégénéré, mais utile, de $s$-équivalence est $R\rsTo Q$, où $R$ est le graphe d'une relation d'équivalence régulière sur $B$, de quotient $Q$, considéré ici comme groupoïde \emph{nul}.

Un \textbf{D}-groupoïde \emph{principal} $R$ est ainsi caractérisé par la propriété d'être $s$-équivalent à un groupoïde \emph{nul}.
\smallskip
\paragraph{\emph{Inducteurs submersifs et ouverts.}} Dans le cas où $\mbf{D}=\mbf{Dif}$, le produit fibré existe toujours quand $p$ est une \emph{submersion} non nécessairement surjective, en particulier quand $E$ est un \emph{ouvert} $U$ de $B$. On peut noter $_U{G}_U$ le groupoïde induit sur $U$.

\subsubsection{\textbf{D}-équivalences, $i$-équivalences.}
\label{deq}
Ayant ainsi transféré au cadre structuré la notion algébrique de \og pleinement fidèle \fg, on saura faire de même pour la notion d'équivalence (catégorique [McL]), dès que l'on aura donné aussi une traduction diagrammatique de la notion de foncteur \emph{essentiellement (ou génériquement) surjectif}. C'est ce que réalise le diagramme (assez classique) suivant (dans \textbf{D}), 
\begin{diagram}[w=1.5em,h=1em,tight,labelstyle=\scriptstyle,textflow,midshaft]
&&&{}&\rLine[abut]\til f&{}&&  {}&\rLine[abut]\til b&{}&&&\\
H&&\ruLine(3,1)[abut]\rTo\til u&&p^{\ast}(G)\SEpbk&&\rdTo(3,1)[abut]\rTo\til q\ruLine(3,1)[abut]&&G  &&\rsTo\til{\beta_G}\rdsTo(3,1)[abut]&&B\\
&\rdsTo(4,4)\til{\alpha_H}&&&&&&&  &&&&\\
&&&&\dsTo\til a&&&&\dsTo\til{\alpha_G}  &&&&\\
&&&&&&&&  &&&&\\
&&&&E&&\rTo\til p&&B  &&&&\\
\end{diagram}
où l'on a factorisé $\text{A}(f)$ en introduisant un produit fibré, et où la surjectivité essentielle se traduit par la surjectivité (ensembliste) de la flèche composée $b=\beta_G\circ q$, que l'on remplacera donc, conformément à notre principe méta\-mathématique directeur, par la condition:
$$b\in D_s.$$

Lorsqu'il en est ainsi, on dira que $f$ est \emph{essentiellement} $s$.

Cette condition est trivialement vérifiée par le $i$-morphisme $\omega_G:B\riTo G$.

\smallskip

Toute scission (= section) d'un $s$-équivalence est une $i$-équivalence.

\smallskip

\subsubsection{\'{E}quivalences submersives et ouvertes.}\label{eqso}
 Lorsque $\mbf{D}=\mbf{Dif}$ et que $p$ est une submersion, la condition de $s$-essentialité se réduit à la condition ensembliste. En particulier, lorsque $H$ est le groupoïde induit par $G$ sur un ouvert $U$ de $B$, elle signifie seulement que $U$ rencontre toutes les orbites du groupoïde: nous parlerons alors d'\emph{équivalence ouverte}.

Un cas particulier important est celui où la variété de base $B$ est somme de deux sous-variétés ouvertes $B_2,\,B_1$, sur lesquels $G$ induit les sous-groupoïdes ouverts $G_2,\,G_1$ et où les restrictions $p_2,\,p_1$ de $\beta_G,\,\alpha_G$ à l'ouvert $$E=\tau_G^{-1}(B_2\times B_1)$$ sont (des submersions) \emph{surjectives}. Alors les injections canoniques $\iota_2,\,\iota_1$ de $G_2,\,G_1$ dans $G$ définissent des \emph{équivalences ouvertes} de $G_2,\,G_1$ vers $G$\,\footnote{C'est cette situation qui a été prise pour définition de la Morita-équivalence dans la conférence d'A. Weinstein. Celle-ci est aussi définie par ce que nous appellerons plus loin la bifibration principale (\ref{dipap}) $B_2\lsTo\til{p_2} E\rsTo\til{p_1}B_1$ et les $s$-équivalences qu'elle détermine entre $G$ et un sous-groupoïde ouvert $K$ de $\square\,G$, induit par $p_1$ ou $p_2$: $G_2\lsTo^{q_2}_\sim K\rsTo^{q_1}_\sim G_1$ (voir \ref{eqmor}). Nous privilégions ici les $s$-équivalences afin que nos constructions s'étendent au cadre diptyque, où l'on ne dispose que des \emph{sur}mersions. Nous verrons (\ref{fract}) qu'\emph{a posteriori} toutes les notions de \textbf{D}-équivalence conduisent à la même notion de Morita-équivalence dans le cadre de la théorie des fractions, qui rend toutes ces équivalences \emph{inversibles}.}. On observera que la variété $G$ est alors somme disjointe des quatre sous-variétés ouvertes
$$G =G_1 \cup G_2 \cup E \cup E^{-1}.$$ Avec les notations de \ref{genot}, on écrira:
$$G_2\rploTo^{\iota_2}_\sim G\lploTo^{\iota_1}_\sim G_1.$$

\subsection{$s$-extenseurs.}\label{sex}
En utilisant les notations de \ref{seq}, la définition ci-dessus se traduit maintenant par le diagramme suivant (que l'on peut lire aussi bien dans $D$ que dans \textbf{Grpd(D)}), où $h$ est par hypothèse un bon épi (et $f=g\circ h$),
\begin{diagram}[h=1.5em,w=2.5em,tight,labelstyle=\scriptstyle,textflow]
K\SEpbk&&\rsTo^ k&&R\SEpbk&&\rsTo^r_\sim&&B\\
&&&&&&&&\\
\diTo\til j&&&&\diTo\til i&&&&\diTo\til{\omega_G}\\
&{}&&\rLine[abut]\til f&&&&{}&\\
H\ruLine(1,1)[abut]&&\rsTo^h&&p^{\ast\ast}(G)\SEpbk&&\rsTo_\sim^g&&G\rdsTo(1,1)[abut]\\
&\rdTo(4,4)\til{\tau_H}&&&&&&&\\
&&&&\dTo\til t&&&&\dTo\til{\tau_G}\\
&&&&&&&&\\
&&&&E\times E&&\rsTo_{\sim}^{p\times p}&&B\times B\\
\end{diagram}
dans lequel on a fait apparaître les noyaux des morphimes, obtenus par des produits fibrés avec $\omega_G$ (noter que les carrés cartésiens se composent, horizontalement et verticalement).

\smallskip
On dira:
\begin{itemize}
	\item que $H$ est une \emph{extension} de $G$ par $K=\text{Ker}f$;
	\item que \begin{diagram}[size=2.5em,inline,labelstyle=\scriptstyle]
K&\riTo^ j&H&\rsTo^ f&G\\
\end{diagram}
est une \emph{suite exacte} (courte) de \textbf{D}-groupoïdes\,\footnote{Il n'y a pas ici de zéro, qu'il faudrait pouvoir écrire à droite et à gauche pour respecter la notation du cas des groupes.}, ainsi que \begin{diagram}[size=2.5em,inline,labelstyle=\scriptstyle]
R&\riTo^ i&p^{\ast\ast}(G)&\rsTo^ g_{\sim}&G\\
\end{diagram}.
\end{itemize}

On notera que le noyau d'une $s$-équivalence est de Godement, donc principal.

\subsection{Le contenu algébrique sous-jacent.}\label{casj}
Observons d'abord le cas particulier où $p$ est l'identité de $B$ (on dira alors que \textbf{f} est un $B$-morphisme, ou est \emph{unifère}), ou plus généralement un isomorphisme (ce qui est toujours le cas pour les morphismes de groupes).

Lorsqu'il en est ainsi les équivalences se réduisent à des isomorphismes, le noyau d'un $s$-extenseur est un \emph{plurigroupe} (\ref{degen}), et la situation généralise formellement mot à mot (\emph{mutatis mutandis}) la théorie des quotients de groupes par un sous-groupe invariant (il faut seulement prendre en compte le fait que la loi de composition n'est pas partout définie).

Dans le cas différentiable, la présentation donnée par J.-P. Serre dans [LALG], essentiellement diagrammatique et reposant sur le théorème de Godement, s'étend alors immédiatement dans notre cadre des diptyques de Godement (qui en est d'ailleurs inspiré).

Le cas général nécessite (pour traiter le cas structuré) des diagrammes un peu plus élaborés (à cause des classes bilatères), et utilise notre notion générale de bon carré cartésien (voir note de bas de page de \ref{dip}); il a été traité dans [P 1986].

Il faut noter toutefois que le traitement du contenu algébrique sous-jacent se réduit ici encore pratiquement à celui du cas des groupes\,\footnote{Il est d'autant plus curieux d'observer que cette situation ultra-élémentaire semble universellement et obstinément ignorée, même d'auteurs très au fait de situations beaucoup plus élaborées.}, la seule différence formelle importante étant qu'ici les classes à droite et à gauche ne coïncident plus, et doivent être systématiquement remplacées par les classes \emph{bilatères}. Pour cette raison nous pensons qu'il est bon d'écrire dans le cas général $G=K\backslash H/K$ ou $H//K$ (quotient bilatère).

Une fois de plus, il est agréable, même du point de vue algébrique, d'unifier ici, par le concept d'extenseur, les deux cas de dégénéréscence du noyau, en (pluri)groupe ou en graphe d'équivalence, lesquels conduisent respectivement aux extenseurs unifères et aux $s$-équivalences.

\smallskip En revanche le paragraphe suivant traite d'un cas de morphisme de groupoïdes fort important (déjà évoqué dans l'introduction) qui dégénère trivialement dans le cas des groupes (car le noyau de ces morphismes est nul!).

\subsection{Acteurs.}\label{acex}
Nous considérons maintenant le carré $\text{A}(f)$ (il serait en fait équi\-va\-lent, bien qu'un peu plus lourd, de considérer $\text{D}(f)$, qui a l'avantage théorique de pouvoir être regardé dans \textbf{Grpd(D)}, et pas seulement dans \textbf{D}).

Il est inutile ici d'imposer la condition $p\in D_s$ (du fait que la propriété $\alpha_G\in D_s$ assure déjà l'existence d'un (bon) produit fibré factorisant le carré), bien qu'elle soit souvent satisfaite dans les applications (nous rencontrerons plus loin des $i$-acteurs à propos du diagramme du papillon, mais aurons surtout affaire à des $s$-acteurs).

On dira que \textbf{f} est un \textbf{D}-\emph{acteur}\,\footnote{Ou \textbf{D}-\emph{foncteur d'action}.}/-\emph{exacteur} si le carré \text{A}(f) est cartésien/plein (\ref{dip}), ce que résume le tableau ci-dessous.

\begin{center}
\begin{tabular}[]{|c|c|}
\hline
\textbf{A(f)}&f\\
\hline
p.b.&acteur\\
\hline
plein&exacteur\\
\hline
\end{tabular}
\end{center}

Un acteur sera dit \emph{principal} lorsque sa source est un groupoïde principal. Cette terminologie sera justifiée dans une prochaine section.

\subsubsection{Acteurs et lois d'action.}\label{ala}
La définition est traduite par le carré cartésien suivant:

\begin{center}
\begin{diagram}[size=3em,tight,labelstyle=\scriptstyle,inline]
H\SEpbk&\rTo\til f&G\\
\dsTo\til{\alpha_H}&&\dsTo\til{\alpha_G}\\
E&\rTo\til p&B\\
\end{diagram}\quad.
\end{center}

Notre terminologie vient de ce que, du point de vue ensembliste sous-jacent, il y a correspondance biunivoque, à isomorphisme près, entre acteurs et lois d'action de groupoïdes\,\footnote{Cette notion a reçu dans la littérature des noms variés, dont aucun n'évoque cette propriété pourtant fondamentale, et dont aucun n'est adapté au cas structuré: foncteurs d'hypermorphismes (Ehresmann), fibrations discrètes (Grothendieck et autres catégoriciens), \emph{star bijection} (Ronnie Brown), déroulements (W.T. van Est); par ailleurs le terme de \og \emph{groupoïde d'action} \fg\, s'est aujourd'hui répandu dans la littérature pour désigner la \emph{source} d'un foncteur d'action, dont la structure de groupoïde ne suffit pas, à elle seule, à caractériser l'action.}.

Nous esquissons, dans le cadre ensembliste, cette correspondance, facile à vérifier (dont la traduction diagrammatique permet de lui attribuer automatiquement un sens dans le cadre structuré). Nous nous aiderons de la figure \ref{act} ci-jointe, où seules les bases sont dessinées, qui schématise la notion d'action d'un groupoïde $G$ de base $B$ sur un ensemble $E$ \emph{muni} d'une flèche $p$ vers $B$ (les fibres de $p$ sont verticales).

\begin{figure}[!htb]
\begin{center}
\fbox{
Action de \begin{diagram}[s=1.5em,p=.5em,inline]G\\ \dsTo\dsTo\\ B\\ \end{diagram} sur 
\begin{diagram}[s=1.5em,inline]E\\ \dTo_p\\ B\\ \end{diagram}}
\end{center}
\fbox{
\begin{diagram}[w=4em,h=2em,p=.5em,abut,tight,labelstyle=\scriptstyle,thick]
{}&\rLine&{}&\rLine&{}&\rLine&{}\\
&&\dLine&\overset{h=(g,x)}\curvearrowleft&\dLine&&\dLine\\
\dLine&&\scriptstyle{\overset{y={g\vdash x}={\beta_H}(h)}\bullet}&
&\scriptstyle{\overset{x={\alpha_H}(h)}\bullet}&&E\\
&&\dLine&&\dLine&&\dLine\\
{}&\rLine&{}&\rLine&{}&\rLine&{}\\
\underset{}G&&{}&\lTo_g&{}&&{}\\
\dsTo^{\beta_G}\dsTo_{\alpha_G}&&&&&&\dTo_p\\
{}&&&\underset{g}\curvearrowleft&&&{}\\
{}&\rLine&\scriptstyle{\overset{p(y)=\beta_G(g)}\bullet}&\rLine&\scriptstyle{\overset{p(x)=\alpha_G(g)}\bullet}&\rLine&B
\end{diagram}}
\begin{tabular}{||c||}
\hline
$H=G\underset{B}{\times}E\stackrel{\beta}{\longrightarrow}E$\\
$(g,x)\longmapsto g\vdash x=\beta(g,x)$\\
\begin{diagram}[h=1.3em,w=2em,tight,labelstyle=\scriptstyle,abut,midshaft]
&{}&&\rLine\til f&&{}&\\
H\ruLine(1,1)&&\riTo\til i&p^{\ast\ast}(G)&&\rTo^{\sim}_{\text{pr}}&G\rdTo(1,1)\\
\end{diagram}\\
$(g,x)\longmapsto(g;(g\vdash x,x))\longmapsto g$\\
\hline
\end{tabular}
\caption{\label{act}}
\end{figure}


Cette action est donnée par une application (définie sur le produit fibré de $\alpha_G$ et de $p$):
$$\beta:H=G\underset{B}{\times}E\longrightarrow E,\quad(g,x)\longmapsto g\vdash x=\beta(g,x)\,.$$

Le graphe de l'application $\beta$, composé avec une permutation des facteurs, dé\-ter\-mine une injection $i$ de $H$ dans le groupoïde $p^{\ast\ast}(G)=G{\underset{B\times B}\times}(E\times E)$ induit\,\footnote{Si l'on ne fait pas l'hypothèse $p\in D_s$, cette étape intermédiaire ne serait valable que dans le cadre ensembliste, mais la correspondance acteurs $\leftrightarrow$ lois d'action reste valable dans le cadre structuré.} par $G$ le long de $p$ (\ref{seq}).

Alors les axiomes des lois d'action de groupoïdes (associativité) se traduisent très exactement par le fait que l'image $i(H)$ est un \emph{sous-groupoïde} de $p^{\ast\ast}(G)$, ce qui munit $H$ d'une structure de groupoïde de base $E$ (connue sous le nom de \emph{groupoïde d'action}), pour laquelle $\beta_H=\beta$.

En composant $i$ avec la projection canonique vers $G$, on munit $H$ d'un morphisme vers $G$, qui est un \emph{acteur} au sens où nous venons de le définir.

\emph{Inversement} la donnée d'un acteur $f:H\rightarrow G$\,, ce qui signifie que l'on a un \textbf{D}-isomorphisme $H\rTo^u_{\approx}G\underset{B}{\times}E$, détermine une loi d'action $\beta=\beta_H\circ u^{-1}$.

Il faut noter que, si $h:H'\rightarrow H$ est un isomorphisme \emph{unifère} (\ref{casj}), les acteurs $f$ et $f'= f\circ h$, définissent la \emph{même} loi d'action. De façon plus précise on établit une équivalence entre la catégorie des morphismes équivariants entre lois d'action et celle des morphismes (\cad des carrés commutatifs) entre acteurs\,\footnote{En fait la suite montrera que le concept d'acteur se transfère au cas structuré beaucoup plus aisément que celui de loi d'action.}.

Une flèche $g$ de $G$ détermine une \emph{bijection entre les fibres} de $p$ au-dessus de sa source et de son but (voir le schéma de la figure \ref{act}).

\paragraph{\emph{Exemple 1}.} En prenant $H=\Delta G$\,\footnote{Considéré comme groupoïde principal, de base $G$, d'espace orbital $B$.}, $f=\delta_G$, $E=G$ et $p=\beta_G:G\rsTo B$, on aura: $$u^{-1}:\Gamma G\rightarrow\Delta G,\,(g,x)\mapsto(gx,x),$$ et l'on voit que:\\
\emph{le $s$-acteur principal canonique $\delta_G:\Delta G\rsTo G$\,\footnote{Dans le cas, évoqué dans l'introduction, où $G$ est un groupe, le groupoïde principal $\Delta G$ devient banal.}, $(y,x)\mapsto yx^{-1}$ décrit l'action $\lambda_G:(g,x)\rsTo gx$ de $G$ sur ($G\rsTo^{\beta_G}B$) par translation à gauche.}

\paragraph{\emph{Exemple 2.}} De même, en notant $\nabla G=\Delta(G^{\text{op}})$\,\footnote{Selon le même principe de notations qu'en \ref{not}.} le produit fibré de la flèche $\beta_G$ par elle-même, on vérifie que $$\overline{\delta}_G:\nabla G\rsTo G,\,(y,x)\mapsto y^{-1}x$$ (division à gauche, ou antidivision) détermine un $s$-acteur principal décrivant l'action (à gauche!) $$\overline{\lambda}_G:(g,x)\mapsto xg^{-1}$$ de $G$ sur ($G\rsTo^{\alpha_G}B$) par translation à droite. On a ici (en gardant toujours les mêmes notations) $$u^{-1}:\Delta G\rightarrow\nabla G,\,(g,x)\mapsto(xg^{-1},x)$$ (ce n'est pas un isomorphisme de groupoïdes principaux!).

Le couple $(\lambda_G,\overline{\lambda}_G)$ est un exemple de ce que nous appellerons plus loin actions principales conjuguées.

\smallskip

La symétrie $\varsigma_G:G\rightarrow G$ (vue ici comme flèche de $D$, échangeant $\alpha_G$ et $\beta_G$, et non pas comme un morphisme de groupoïdes) est équi\-variante par rapport aux deux actions $\lambda_G$ et $\overline{\lambda}_G$ de \textbf{G}. Ceci se traduit par le carré commutatif (dans \textbf{Grpd(D)}):
\begin{center}
\begin{diagram}[h=3em,w=4.5em,tight,labelstyle=\scriptstyle,inline]
\Delta G&\rTo\til {(\varsigma_G,\varsigma_G)}&\nabla G\\
\dsTo\til{\delta_G}&\text{Sym}(G)&\dsTo\til{\overline{\delta}_G}\\
G&\rEqual&G\\
\end{diagram}
\quad.
\end{center}

\paragraph{\emph{Remarque.}} 
Le carré précédent s'écrit, de façon plus fonctorielle, comme ci-dessous, 
\begin{diagram}[h=3em,w=5em,tight,labelstyle=\scriptstyle,textflow]
\Delta G&\rTo\til {\Delta(\varsigma_G)}&\Delta(G^{\text{op}})\\
\dsTo\til{\delta_G}&\text{Sym}(G)&\dsTo\til{\overline{\delta}_G}\\
G&\rEqual&G\\
\end{diagram}
et peut s'obtenir comme suit, avec les notations de \ref{notm}.

D'abord $\overline{\delta}_G$ s'interprète comme la diagonale du carré $\text{S}(\delta_{G^{\text{op}}})$, ou encore de son transposé:

\begin{center}
\begin{diagram}[h=3em,w=6.5em,tight,labelstyle=\scriptstyle,inline] 
\Delta({G^{\text{op}}})&\rTo\til{\varsigma_{(\Delta({G^{\text{op}}}))}}&(\Delta({G^{\text{op}}}))^{\text{op}}\\
\dsTo\til{\delta_{G^{\text{op}}}}
&^{\text t}({\text{S}(\delta_{G^{\text{op}}})})&\dsTo\til{(\delta_{G^{\text{op}}})^{\text{op}}}\\
G^{\text{op}}&\rTo\til{\varsigma_{G^{op}}}&G\\
\end{diagram}\quad.
\end{center}

\'Ecrivons alors la composition horizontale du carré transposé $^{\text t}({\text{S}(\delta_{G^{\text{op}}})})$ avec le carré $\text{D}(\varsigma_G)$, soit (en remplaçant $\Delta(G^{\text{op}})$ par $\nabla G$):

\begin{center}
\begin{diagram}[h=3em,w=5.5em,tight,labelstyle=\scriptstyle,inline]
\Delta G&\rTo\til {\Delta(\varsigma_G)}&\nabla G&\rTo\til{\varsigma_{(\nabla G)}}&(\nabla G)^{\text{op}}\\
\dsTo\til{\delta_G}&\text{D}(\varsigma_G)&\dsTo\til{\delta_{(G^{\text{op}})}}&\rdTo(2,2)\til{\overline{\delta}_G}&\dsTo\til{(\delta_{G^{\text{op}}})^{\text{op}}}\\
G&\rTo\til{\varsigma_G}&G^{\text{op}}&\rTo\til{\varsigma_{(G^{op}})}&G\\
\end{diagram}\quad.
\end{center}
Ceci redonne aussitôt le carré $\text{Sym}(G)$, puisque la flèche composée du bas donne l'identité.

On peut aisément vérifier le contenu ensembliste sous-jacent à ce diagramme commutatif, en observant que l'on a:
$$\varsigma_{(G^{op})}=\varsigma_G,\,\Delta({\varsigma_G})(y,x)=(y^{-1},x^{-1}),\,
\varsigma_{\Delta_G}(y,x)=(x,y),$$
$$(\delta_{G^{\text{op}}})^{\text{op}}=\delta_G,\,\delta_{(G^{\text{op}})}(y,x)=y{\underset{G^{\text{op}}}{.}}x^{-1}=x^{-1}{\underset{G}.}y=\overline{\delta}_G(x,y)\,.$$

\bigskip

Les actions (à gauche) du groupoïde opposé $\mathbf{G^{op}}$ sont aussi appelées \emph{actions à droite} de \textbf{G}. Nous les noterons $(x,g)\mapsto x\dashv g$. La symétrie $\varsigma_G$ associe à toute action à droite $x\dashv g$ une action à gauche $g\vdash x=x\dashv g^{-1}$. Les acteurs associés sont $f:H\rightarrow G$ et la diagonale du carré $\text{S}(f)$, à savoir $f^{\varsigma}=\varsigma_G\circ f:H\rightarrow G^\text{op}$.
En particulier les actions \emph{à droite} $\rho_G=(\lambda_G)^\varsigma$ et $(\overline\rho)_G=(\overline\lambda_G)^\varsigma$ associées à $\lambda_G$ et $\overline{\lambda}_G$ sont données par $x\dashv g=g^{-1}x$ et $x\dashv g=xg$.

Ceci est résumé par le tableau suivant, où l'on utilise partout la multiplication de $G$ (et non de $G^{\text{op}}$\,!).

\begin{center}
\begin{tabular}[]{||c||c|c||c|c||}
\hline
{}&\multicolumn{2}{|c||}{Actions de \textbf{G}}&\multicolumn{2}{|c||}{Actions de $\mathbf{G}^{\text{op}}$}\\
\hline\hline
\text{Acteurs:}&$\delta_G$&$\overline{\delta}_G$&$(\delta_G)^\varsigma$&
$(\overline{\delta}_G)^\varsigma$\\
\hline
$(y,x)\mapsto$&$yx^{-1}$&$y^{-1}x$&$xy^{-1}$&$x^{-1}y$\\
\hline
\hline
Actions à gauche:&$\lambda_G$&${\overline\lambda}_G$&&\\
\hline
$g\vdash x=$&$gx$&$xg^{-1}$&&\\
\hline\hline
Actions à droite:&&&
$\rho_G=(\lambda_G)^\varsigma$&${\overline{\rho}}_G=({\overline\lambda}_G)^\varsigma$\\
\hline
$x\dashv g=$&&&$g^{-1}x$&$xg$\\
\hline
\end{tabular}
\end{center}

\bigskip

Contrairement à un usage répandu dans la théorie des fibrés principaux, nous nous ramènerons toujours en principe à des actions à gauche\,\footnote{Vu la multiplicité des combinaisons possibles pour les choix des conventions concernant les sources et les buts, les actions à droite ou à gauche, et les constructions du groupoïde d'action associé à une loi d'action (combinaisons que les auteurs variés se sont efforcés d'exploiter toutes!), nous avons cru utile d'expliciter les formules ci-dessus concernant les quatre lois d'action à gauche ou à droite associées aux translations à gauche ou à droite dans \textbf{G} ou $\mathbf{G^{op}}$, et les groupoïdes d'action et acteurs associés.}.

Particulièrement important est le cas des \emph{acteurs principaux}, associés aux lois d'action \emph{principales}, \cad le cas où \textbf{H} est principal (\ref{degen}); nous le noterons alors souvent \textbf{R}. Nous consacrerons un paragraphe spécial à ce cas.

\subsubsection{Acteurs locaux}\label{aloc}
Considérons maintenant le \emph{flot} d'un champ de vecteurs sur une variété \emph{séparée} $B$\,\footnote{Nous prenons donc ici $\mbf{D}=\mbf{DifHaus}$.}, qui est défini dans un ouvert de ${\mbb R}\times B$:
$$(t,x)\mapsto f_t(x).$$
Bien qu'il s'agisse seulement d'une action \emph{locale} de $\mbb R$, l'application:
$$(t,x)\mapsto(t,(f_t(x),x)),\quad U\rightarrow{\mbb R}\times({B\times B})$$
définit encore un difféomorphisme de $U$ sur un sous-groupoïde de Lie $H$ plongé dans ${\mbb R}\times({B\times B})$, ceci traduisant la propriété d'associativité d'un flot.

Ici le morphime $f:H\rightarrow{\mbb R}$ définit un \emph{plongement ouvert} de $H$ dans le produit ${\mbb R}\times B$ , qui n'est bijectif que si le flot est complet. On dira que $f$ est un \emph{acteur local}.
\paragraph{\emph{Remarque.}} La théorie de R. Palais de la globalisation d'un loi d'action locale peut être étendue dans le cadre des diptyques en utilisant le diagramme du papillon indiqué plus loin et la notion d'actions conjuguées, mais ceci ne peut être développé ici.

\subsection{Exacteurs.}\label{exa}
La définition donnée ci-dessus est traduite par le diagramme suivant, qui factorise $\text{A}(f)$ (dans \textbf{D})\,\footnote{La condition ensembliste sous-jacente joue un rôle important. Dans la littérature catégorique, elle est généralement connue sous le nom de (op)fibrations ou foncteurs fibrants, et \emph{star surjections} chez Ronnie Brown, terminologies qui ne conviennent pas dans le cadre structuré.}.
\begin{center}
\begin{diagram}[size=1.5em,tight,labelstyle=\scriptstyle,textflow]
&{}&&&\rLine[abut]^f&&&{}&\\
H\ruLine(1,1)[abut]&&\rsTo^s&&p^{\ast}(G)\SEpbk&&\rTo^q&&G\rdTo(1,1)[abut]\\
&\rdsTo(4,4)\til{\alpha_H}&&&&&&&\\
&&&&\dsTo\til a&&&&\dsTo\til{\alpha_G}\\
&&&&&&&&\\
&&&&E&&\rTo^p&&B\\
\end{diagram}
\end{center}

Dans le cadre structuré, elle fournit la bonne condition pour l'existence de produits fibrés\,\footnote{Les produits fibrés de morphismes de groupoïdes existent toujours dans le cadre ensembliste, mais, dans le cadre des groupoïdes de Lie, contrairement à ce qui semble trop souvent admis implicitement à tort dans la littérature, la seule existence du produit fibré en tant que variété ne suffit pas à assurer que celle-ci-ci soit le support d'un groupoïde de Lie (on doit assurer la condition fondamentale $\alpha_G\in D_s$)!}. Mieux, on peut montrer que les $s$-exacteurs peuvent (entre autres choix possibles) être pris pour les bons épis d'une structure de diptyque (de Godement) sur la ca\-té\-go\-rie \textbf{Grpd(D)}, mais ceci ne peut être développé ici.

\smallskip

\paragraph{\emph{Remarques}.} 
\begin{enumerate}
	\item On a: $s$-extenseur $\Rightarrow$ $s$-exacteur.
	\item un exacteur est essentiellement $s$ (\ref{deq}) si et seulement si c'est un $s$-exacteur.
\end{enumerate}

\section{Conjugaison des actions principales.}
\subsection{Fibrés principaux et \textbf{D}-groupoïdes $s$-transitifs.}\label{fpgt}
Comme il est dit dans l'annexe historique A, la première motivation de C. Ehres\-mann pour l'introduction des groupoïdes différentiables en Géométrie fut l'étude des fibrés, et d'abord des fibrés principaux à groupe structural auxquels ils sont associés. Le \og groupoïde des isomorphismes de fibre sur fibre \fg\, est ce qui est souvent appelé aujourd'hui le \emph{groupoïde de jauge}. C'est un \textbf{D}-groupoïde $s$-transitif\,\footnote{Rappelons une fois de plus que le terme de groupoïdes de Lie est resté longtemps réservé au cas particulier des \textbf{D}-groupoïdes $s$-transitifs, avec $\mathbf D=\mathbf{DifHaus}$ (\ref{dip}).}, qui peut être décrit comme quotient d'un groupoïde banal $E\times E$ par l'action \emph{diagonale} du groupe structural (évidemment compatible avec cette loi triviale de groupoïde).

\emph{Inversement}, si \textbf{G} est un \textbf{D}-groupoïde $s$-transitif de base $B$, le \emph{choix} d'un point $e\in B$ détermine sur $E= G_e$ une structure de \emph{fibré principal} de base $B$, dont la projection est $p=\beta_{|E}$, admettant $\mathbf{G'}={_eG_e}$ pour \emph{groupe structural} et \textbf{G} pour groupoïde de jauge.

Les actions (principales) de \textbf{G} et $\mathbf{G'}$ sur $E$ sont les restrictions respectives à $E$ des actions (à gauche sur $G$) $\lambda_G$ et $\overline{\lambda}_G$, décrites dans \ref{ala} (Exemple 2), lesquelles \emph{commutent}.

La construction de l'équivalence ensembliste $G\rDotsTo^{\sim} G'$ décrite dans \ref{tect} ne subsiste pas dans le cadre structuré (en l'absence d'une section $\eta$ de $p$), mais nous esquisserons la construction d'une double $s$-équivalence:
$$\mathbf{G}\lsTo^{\sim}\mathbf{K}\rsTo^{\sim}\mathbf{G'}$$
qui sera un cas particulier d'équivalence de Morita.

Nous donnerons de cette construction une description diagrammatique. Celle-ci, non seulement transfère la construction dans le cadre des diptyques\,\footnote{La théorie classique se déroule dans le diptyque \textbf{DifHaus} (\ref{dip}), mais reste automatiquement valable dans le diptyque \textbf{Dif} (variétés non supposées séparées et plongements non supposés fermés).}, mais rétablit, dans un cadre plus large, la \emph{symétrie} des rôles joués par le groupe structural et le groupoïde de jauge (masquée par des dégénérescences très dissymétriques), et définit une remarquable correspondance involutive (à isomorphisme près), que nous appelons \emph{conjugaison}, entre acteurs principaux (donc entre lois d'action principales, ou entre classes d'isomorphie de fibrés principaux).

\subsection{Acteurs principaux.}
\subsubsection{Le diagramme du papillon.}\label{dipap}
Le paragraphe précédent fournit l'exemple fondamental de la double dégénérescence (dissymétrique!) de la notion générale de $s$-acteur principal (\ref{ala}):
$$\mathbf{R} \rsTo^r_\act \mathbf{G}.$$
Dans un cas, \textbf{G} dégénère en un \emph{groupe} $\mathbf{G'}$ (groupe structural), et le graphe d'équi\-va\-lence régulière $\mathbf{R'}$ définit la \emph{base} $B$ de l'espace fibré $E$, qui est aussi celle du groupoïde de jauge $s$-transitif \textbf{G}. Dans l'autre, c'est le groupoïde principal \textbf{R} qui dégénère en groupoïde \emph{banal}; la relation d'équivalence dégénère en l'équivalence \emph{\emph{grossière}}, dont le quotient (\emph{ponctuel}) définit la base (dégénérée) du groupe structural.

Il est très remarquable que, via le théorème de Godement, la description classique, relativement compliquée, d'une fibration principale (localement triviale), telle qu'on la trouve par exemple dans [VAR], soit équivalente à la donnée (fort simple) d'un $s$-acteur principal, et prenne ainsi automatiquement un sens dans le contexte de tout diptyque de Godement. En revanche, dans le contexte purement topologique, les conditions de locale trivialité n'ont pas de bonne traduction diagrammatique.

La situation générale, parfaitement symétrique\,\footnote{Cette situation est décrite dans le cadre topologique dans l'article d'A. Haefliger de [Ast. 116], selon une présentation attribuée à G. Skandalis. Il semble que, dans ce contexte, la stabilité des conditions de locale trivialité par composition des équivalences de Morita, et des morphismes généralisés (qui n'est d'ailleurs invoquée dans cet article que dans des cas particuliers) soulève des difficultés.}, est décrite par le diagramme suivant, 
\begin{diagram}[w=4em,h=2em,tight,labelstyle=\scriptstyle,textflow]
R&&&&R'\\
&\rdiTo\til{i}&&\ldiTo\til{i'}&\\
\dsTo^\act\til{r}&&K&&\dsTo_\act\til{r'}\\
&\ldsTo_\sim\til{q}&&\rdsTo_\sim\til{q'}&\\
G&&&&G'\\
\end{diagram}
qui est un cas particulier (symétrique) du \og diagramme du papillon \fg\,.

Dans ce diagramme (commutatif de \textbf{Grpd(D)}):
\begin{itemize}
	\item $R$ et $R'$ sont des \textbf{D}-groupoïdes principaux (\ref{god});
	\item $G$ et $G'$ sont des \textbf{D}-groupoïdes de bases respectives $B$ et $B'$;
	\item $K$ est un \textbf{D}-groupoïde de base $E$;
	\item $q$ et $q'$ sont des \emph{$s$-équivalences} (\ref{seq}), de noyaux respectifs $R'$ et $R$, de sorte que les \emph{diagonales} du diagramme sont des \emph{suites exactes} (courtes); les restrictions de $q$, $q'$ aux bases seront notées:
	$$p=q^{(0)}:E\rsTo B,\,p'={q'}^{(0)}:E\rsTo B';$$ par suite $B$ et $B'$ sont les quotients de $E$ par les relations d'équivalence définies respectivement par $R'$ et $R$;
	\item $r$ et $r'$ sont des \emph{$s$-acteurs}, qui seront dits \emph{conjugués}.
\end{itemize}

On dit parfois que $B'\lsTo^{p'} E\rsTo^p B$ est une \emph{bi-fibration principale}.

\smallskip La propriété remarquable du diagramme précédent est que la seule connaissance de l'acteur principal $r$ détermine $r'$ à isomorphisme près, d'où le terme de \emph{conjugaison}.

\bigskip

La dégénérescence dissymétrique décrite au paragraphe précédent \ref{fpgt} est celle qui suit,
\begin{diagram}[w= 3em,h=1.5em,tight,labelstyle=\scriptstyle,textflow]
E\underset{\{e\}}{\times}E&&&&E\underset{B}{\times}E\\
&\rdiTo^i&&\ldiTo^{i'}&\\
\dsTo^\delta&&K&&\dsTo_\delta'\\
&\ldsTo^\sim_s&&\rdsTo^\sim_{s'}&\\
G&(&\rDashTo^\sim&)&G'\\
\end{diagram}
 où $\delta$ et $\delta'$ sont des restrictions de $\delta_G$, ${\overline\delta}_G$.

La flèche en pointillés du bas rappelle l'équivalence purement ensembliste décrite en \ref{tect}. On lui donnera dans la suite un sens nouveau restant valable dans le contexte structuré (équivalence de Morita).

\subsubsection{Fibrés associés.}\label{fibasrem}
Nous ouvrons ici une simple parenthèse, car nous serions entraînés trop loin de notre sujet, pour indiquer comment se situerait la théorie des fibrés associés à un fibré principal.

Il convient de décrire l'action du groupe structural sur la fibre type par un acteur (qui n'est plus en général principal) et de former le \emph{carré cartésien} décrivant le produit fibré avec l'acteur principal, que nous appellerons \emph{carré principal}. Le diagramme s'étend alors en un diagramme tridimensionnel (on rajoute une épaisseur), qui permet d'étendre la conjugaison aux carrés principaux.

Ceci généralise la construction du fibré associé, laquelle conjugue l'action du groupe structural sur la fibre type avec l'action du groupoïde de jauge sur le fibré associé.

Entre autres applications de ce type de construction, on peut ainsi retrouver et généraliser considérablement la théorie de la globalisation des lois d'action de groupes de R. Palais\ldots

\subsection{Structure du c\oe ur du papillon.}\label{scp}
Le groupoïde \textbf{K} (de base $E$) possède une structure supplémentaire très riche et extrêmement remarquable (ici seulement esquissée), que l'on pourrait appeler structure d' \emph{actoïde bi-principal}, pour indiquer que, bien que non principale, sa structure de groupoïde peut être reconstituée, de différentes façons, à partir de groupoïdes principaux, \cad dégénérés, autrement dit à partir de relations d'équivalence (nous donnerons plus loin, dans la figure \ref{fpac} un schéma qui aidera à visualiser la situation algébrique):
\begin{itemize}
	\item il est induit (\ref{seq}) à la fois par $G$ le long de $q$ et par $G'$ le long de $q'$;
	\item il se décompose en produit (\og tordu \fg\,) des deux sous-groupoïdes (principaux) $R$ et $R'$, en ce sens que les \emph{restrictions} de $\delta_K$ et ${\overline\delta}_K$ induisent des \textbf{D}-isomorphismes:
$$\kappa:R\underset{\alpha}{\times}R'\rightarrow K,\,
\overline\kappa:R'\underset{\alpha}{\times}R\rightarrow K;$$
	\item cette décomposition permet de munir $K$ d'une structure de \emph{groupoïde double} (au sens d'Ehresmann) $\mathfrak K=(\mathfrak R,\mathfrak R')$, de bases respectives $R',R$.
\end{itemize}

Cette dernière propriété signifie que les flèches $\omega_{\mathfrak R},\alpha_{\mathfrak R},\delta_{\mathfrak R}$ définissant la structure de groupoïde $\mathfrak R$ sont des morphismes relativement à l'autre\,\footnote{On notera qu'en fait $\omega_{\mathfrak R}$ et $\delta_{\mathfrak R}$ seront donc des morphismes \emph{doubles.}} structure de groupoïde $\mathfrak R'$ (et de même en échangeant $\mathfrak R$ et $\mathfrak R'$).

Mais cette structure de groupoïde double possède ici deux propriétés sup\-plé\-men\-taires très particulières:
\begin{enumerate}
	\item $\mathfrak K$ est ce que nous appellerons un \emph{actoïde}, ce qui signifie que les morphismes $\omega_{\mathfrak R},\alpha_{\mathfrak R},\delta_{\mathfrak R}$, sont en fait ici des \emph{acteurs} (relativement aux structures de groupoïdes associées à $\mathfrak{R'}$), et symétriquement\,\footnote{On peut montrer que, de façon générale, un actoïde définit une \emph{troisième loi de groupoïde}, la \emph{loi mixte} (ou \emph{diagonale}), qui n'est autre ici que la loi de \textbf{K}. La construction de la loi mixte, dans le cas général, peut être interprétée comme une généralisation \emph{symétrique} du produit semi-direct. La situation algébrique sous-jacente a été rencontrée dans des contextes variés et divers cas particuliers, sous des noms différents, par de nombreux auteurs, entre autres Ronnie Brown, K. Mackenzie, A. Weinstein, P. Dazord, \emph{et alii}, et nous ne chercherons pas à faire l'historique de cette notion, qui, sous sa forme diagrammatique, acquiert des significations riches et variées dans les différents contextes structurés.};
	\item cet actoïde $\mathfrak K$ est \emph{bi-principal}, en ce sens que les groupoïdes $\mathfrak R$ et $\mathfrak R'$ sont principaux.
\end{enumerate}

Les lois de $\mathfrak R$ et $\mathfrak R'$ sont induites par deux lois banales sur $E^2$:
$$({\tau_{R'}}\times{\tau_{R'}})\circ{\tau_{\mathfrak R}}:{\mathfrak R}\riTo E^4$$
$$({\tau_{R}}\times{\tau_R})\circ{\tau_{\mathfrak R'}}:{\mathfrak R'}\riTo E^4\,.$$

La relation d'équivalence sur $R$ définie par le graphe $\mathfrak R'$ a pour quotient l'objet de \textbf{D} sous-jacent à \textbf{G}, ce que l'on peut traduire (\emph{cf.} \ref{seq}) en disant que $q$ détermine une $s$-équivalence $\mathfrak R'\rsTo^{\sim}|G|$, où la notation $|\quad|$ désignera le groupoïde \emph{nul} associé à l'objet de $D$ sous-jacent à un groupoïde. Mais par ailleurs $q$ détermine aussi un $s$-acteur: $\mathfrak R\rsTo^\act G$ (décrivant une action de $G$ sur $R'$).

Ceci conduit finalement à une \emph{nouvelle lecture} du diagramme du papillon (ayant le \emph{même} diagramme sous-jacent dans $D$)
\begin{diagram}[w=4.5em,h=2.5em,tight,labelstyle=\scriptstyle,textflow]
(R,|R|)&&&&(|R'|,R')\\
&\rdiTo_\act\til{\omega_{\mathfrak{R'}}}&&\ldiTo_\act\til{\omega_{\mathfrak{R}}}&\\
\dsTo^{\act}_{\act}\til{r}&&(\mathfrak R,\mathfrak{R'})&&\dsTo_{\act}^{\act}\til{r'}\\
&\ldsTo^{\act}_\sim\til{q}&&\rdsTo_\sim^\act\til{q'}&\\
(G,|G|)&&&&(|G'|,G')\\
\end{diagram}
\emph{dans la catégorie des morphimes doubles} (entre groupoïdes doubles\,\footnote{Qui sont en fait tous des \emph{actoïdes}(\ref{scp}), présentant diverses dégénérescences.}). Chacun de ces morphismes détermine \emph{deux} morphismes de groupoïdes, ayant le même support dans $D$, mais des propriétés dif\-fé\-rentes, indiquées par les deux étiquettes, \emph{dans l'ordre indiqué}.

Ici les diagonales ne sont plus des suites exactes. Le même couple de \textbf{D}-flèches détermine ici un couple d'acteurs principaux:
\begin{center}
$G\lsTo^q_\act\go R$, $\go R'\rsTo^{q'}_\act G'$
\end{center}
représentant les actions \emph{diagonales} (principales) associées aux actions définies par le couple
\begin{center}
$G\lsTo^r_\act R$, $R'\rsTo^{r'}_\act G'$
\end{center}
et non plus le couple d'équivalences:
$$\mathbf{G}\lsTo^q_{\sim}\mathbf{K}\rsTo^{q'}_{\sim}\mathbf{G'}.$$

\subsection{Le papillon canonique d'un \textbf{D}-groupoïde.}\label{papcan}
Nous donnons ici le diagramme canonique décrivant la conjugaison annoncée des lois $\lambda_G$ et ${\overline\lambda}_G$ (associées aux translations à droite et à gauche), \cad des acteurs $\delta_G$ et ${\overline\delta}_G$ (\ref{ala}). 
\begin{diagram}[w=4em,h=2em,tight,labelstyle=\scriptstyle,textflow]
\Delta G&&&&\nabla G\\
&\rdiTo\til{\iota_G}&&\ldiTo\til{{\overline\iota}_G}&\\
\dsTo^\act\til{\delta_G}&&\square G&&\dsTo_\act\til{{\overline\delta}_G}\\
&\ldsTo_\sim\til{\pi^\downarrow_G}&&\rdsTo_\sim\til{\pi^\uparrow_G}&\\
G&&&&G\\
\end{diagram}
Le c\oe ur du diagramme $\square G$ (notation d'Ehresmann) est ici le grou\-poïde des carrés commutatifs (ou quatuors) de $G$ (qu'il conviendra ici d'écrire avec des flèches dirigées vers la gauche et vers le bas), muni de sa loi de composition \emph{horizontale} (qu'il est bon d'écrire ici de droite à gauche). Les flèches $\pi^\downarrow_G,\pi^\uparrow_G$ sont, comme la notation le suggère, les projections sur les côtés du carré situés en bas ou en haut.

Les décompositions $\kappa$ et $\kappa'$ ci-dessus selon les deux sous-groupoïdes principaux prennent ici l'aspect (très géométrique!) des décompositions du carré $\square$ (à l'aide d'une diagonale) en deux triangles juxtaposés $\Delta\!\nabla$, ou $\nabla\!\Delta$, vus comme des carrés avec un côté dégénéré.

Les actions diagonales évoquées en \ref{scp} sont ici les actions de $G$ aux sommets du haut ou du bas des triangles $\Delta$ et $\nabla$.

On prendra bien garde que la structure de groupoïde double (ici actoïde principal) décrite ci-dessus n'a rien à voir avec la structure de groupoïde double bien connue définie par les lois horizontale et verticale (cette dernière n'intervient pas ici)!

La bi-fibration n'est autre que 
$$B\lsTo^{\alpha_G}G\rsTo^{\beta_G}B.$$

\subsection{Visualisation des trois lois de groupoïde sur le c\oe ur du papillon.}
\begin{figure}[!htb]
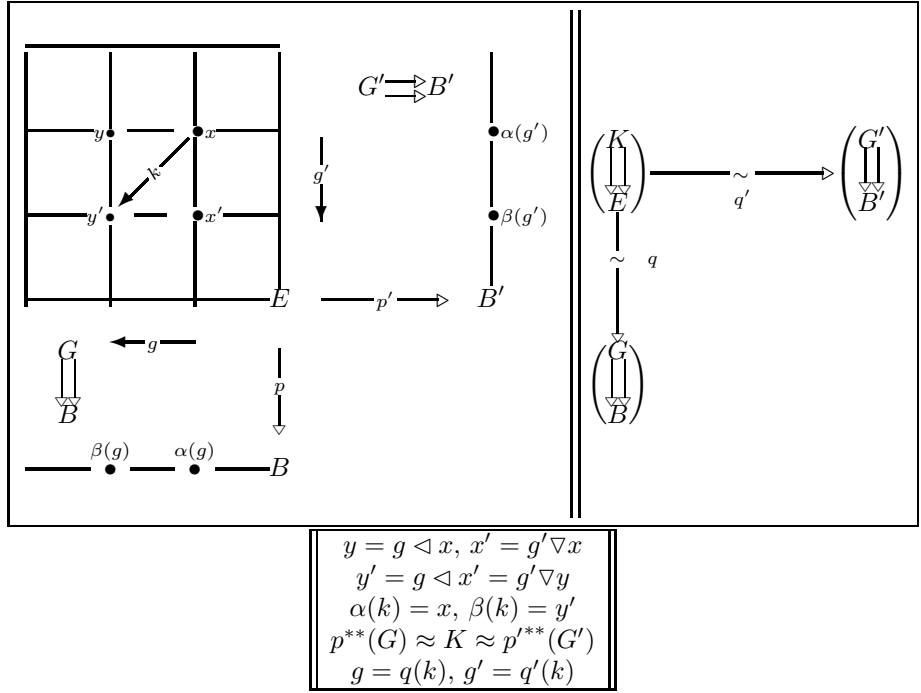

\def\KK{{
\begin{diagram}[size=1.2em,p=.5em,tight,thin]
K\\ \dsTo\dsTo\\E
\end{diagram}}}
\def\gg{{
\begin{diagram}[size=1.2em,p=.5em,tight,thin]
G'\\ \dsTo\dsTo\\B'
\end{diagram}}}
\def\gh{{\begin{diagram}[w=1.3em,p=0.5em,tight,labelstyle=\scriptstyle,thin]
G'&\pile{\rsTo[abut]\\ \rsTo[abut]}&B'\\
\end{diagram}}}
\def\GG{{
\begin{diagram}[size=1.2em,p=.5em,tight,thin]
G\\ \dsTo\dsTo\\B
\end{diagram}}}
\fbox{
\begin{diagram}[w=1.6em,h=1.6em,p=.3em,tight,labelstyle=\scriptstyle,abut,midshaft,inline,thick]
&&&&&&&&&&&&&&&&&&&&\\
{}&\rLine&{}&\rLine&{}&\rLine&{}&&&&&{}&&&&&&&&\\
\dLine&&\dLine&&\dLine&&\dLine&&&\gh&&\dLine&&&&&&&&&\\
{}&\rLine&\sst{y}\bullet\phantom{\sst{y}}&\rLine&\phantom{\sst{x}}\bullet\sst{x}&\rLine&{}&{}&&&&\phantom{\scriptstyle{\alpha(g')}}\bullet\scriptstyle{\alpha(g')}&&&&&&&&&\\
\dLine&&\dLine&\ldTo\til{k}&\dLine&&\dLine&\dTo\til{g'}&&&{}&\dLine&&\dLine\dLine&\left(\KK\right)&&&\rsTo_{q'}\til{\sim\,}&&&\left(\gg\right)\\
{}&\rLine&\sst{y'}\bullet\phantom{\sst{y'}}&\rLine&\phantom{\sst{x'}}\bullet\sst{x'}&\rLine&{}&{}&&&&\phantom{\scriptstyle{\beta(g')}}\bullet\scriptstyle{\beta(g')}&&&&&&&&&\\
\dLine&&\dLine&&\dLine&&\dLine&&&&&\dLine&&&\dsTo\til{\sim}_q&&&&&&\\
{}&\rLine&{}&\rLine&{}&\rLine&{E}&{}&&\rsTo\til{p'}&{}&B'&&&&&&&&&\\
&&{}&\lTo\til{g}&{}&&{}&&&&&&&&&&&&&&\\
&\GG&{}&&&&\dsTo\til{p}&&&&&&&&\left(\GG\right)&&&&&&\\
&&&&&&{}&&&&&&&&&&&&&&\\
{}&\rLine&{\overset{{\beta(g)}}\bullet}&\rLine&\overset{{\alpha(g)}}\bullet&\rLine&B&&&&&&&&&&&&&&\\
&&&&&&&&&&&&&&&&&&&&\\
\end{diagram}}
\begin{tabular}{||c||}
\hline
$y=g\lhd x,\,x'=g'\triangledown x$\\
$y'=g\lhd x'=g'\triangledown y$\\
$\alpha(k)=x,\,\beta(k)=y'$\\
$p^{\ast\ast}(G)\approx K\approx{p'}^{\ast\ast}(G')$\\
$g=q(k),\,g'=q'(k)$\\
\hline
\end{tabular}
\caption{\label{fpac}Conjugaison de deux fibrations principales.}
\end{figure}

Le schéma ci-joint (fig. \ref{fpac}) aide à visualiser, dans le cas général (symétrique), le contenu ensembliste ou algébrique, très élémentaire, sous-jacent au diagramme du papillon, et à comprendre la structure du c\oe ur du papillon.

La correspondance entre groupe structural et groupoïde de jauge pour un fibré principal classique est le cas particulier où $B'$ est un \emph{singleton}, la relation d'équi\-va\-lence correspondante étant \emph{grossière} (elle passe inaperçue!). Cette dissymétrie masque la simplicité de la description géométrique de la situation générale.

Sur la figure, les fibres horizontales/verticales sont les classes de l'équivalence définie par le graphe $R$/$R'$, associée à la projection $p'$/$p$.

Le groupoïde $G$/$G'$ est le quotient de $R$/$R'$ par la relation d'équivalence dont la projection est $q$/$q'$, et dont le graphe est $\go R'$/$\go R$, avec la compostion verticale/ hori\-zon\-tale; la deuxième structure de groupoïde de ce dernier, donnée par la loi horizontale/  ver\-ti\-cale $\go R$/$\go R'$, donne, par la loi quotient, la structure de groupoïde (non nul) de $G$/$G'$.

L'action de $G$ (sur les fibres verticales) est représentée horizontalement, et celle de $G'$ verticalement (sur les fibres horizontales). Ces actions commutent. Ceci est mémorisé en notant ces actions respectivement:
\begin{center}
$g\lhd x$ et $g'\triangledown x.$
\end{center}

Les éléments de $K$ sont figurés par les quatre sommets d'un carré, dont les sommets adjacents sont, horizontalement ou verticalement, équivalents.

Pour s'en tenir au point de vue purement ensembliste ou algébrique, la \emph{propriété fondamentale}, traduisant le fait que $K$ est un \emph{actoïde} (principal), (\ref{scp}) est que \emph{trois des sommets du carré décrivant $k$ déterminent uniquement le quatrième}, ce que l'on peut exprimer en disant que les relations d'équivalence $R$ et $R'$ sont \emph{transverses}.

La loi de composition horizontale/verticale de $\go R$/$\go R'$ se lit immédiatement sur la figure, d'où la structure de groupoïde double $\go K=(\go R,\go R')$. Rappelons que les flèches horizontales sont dans $R$, et les flèches verticales dans $R'$.

Quant à la loi (non dégénérée!) de $K$ (loi mixte ou diagonale de l'actoïde), elle se lit aisément sur le schéma ci-joint, 
\begin{diagram}[size=2em,tight,abut,labelstyle=\scriptstyle,textflow,thick]
(z)&\lTo&y&\lTo&x\\
\dTo&&\dTo&\ldTo\til{k}&\dTo\\
z'&\lTo&y'&\lTo&x'\\
\dTo&\ldTo\til{k'}&\dTo&&\dTo\\
z^{\prime\prime}&\lTo&y^{\prime\prime}&\lTo&(x^{\prime\prime})\\
\end{diagram}
compte tenu de la propriété fondamentale, qui permet de compléter le grand carré (sommets $x^{\prime\prime}$ et $z$), en complétant les deux petits carrés incomplets.

\bigskip
On prendra garde que, lorsqu'on applique la construction ci-dessus au cas particulier du papillon canonique d'un groupoïde \textbf{G} (\ref{papcan}) (avec $E=G$), les \emph{sommets} du carré ci-dessus, qui sont ici des flèches de $ G$, ne sont pas les quatre côtés du carré commutatif de $G$, mais les deux côtés verticaux et les deux diagonales, tandis que les \emph{côtés} du carré ci-dessus sont maintenant les quatre \emph{triangles} des deux décompositions du carré $\Delta\!\nabla$, et $\nabla\!\Delta$.

Lorsqu'on munit $\square G$ de la structure de groupoïde double ainsi définie (à ne pas confondre avec celle qui est définie par les compositions horizontale et verticale!), on pourra le noter $\boxtimes G$, pour évoquer la décomposition en quatre triangles admettant chacun une diagonale pour un de ses bords.

\section{Inversion des équivalences: équivalences de Morita.}\label{fract}
La référence est ici [P 1989], où figurent les énoncés précis (fascicule de résultats).
\subsection{Introduction.}
\label{intro}

L'équivalence de Morita, ou \textbf{D}-Morita-équivalence est ici introduite dans le cadre général d'un diptyque fixé \textbf{D}\,\footnote{Pour faciliter le langage, et bien que ce ne soit pas indispensable, nous supposons ici implicitement que \textbf{D} est muni d'un foncteur d'oubli vers \textbf{Ens}, et considérons les objets de \textbf{D} comme des structures possédant un support ensembliste, pour pouvoir utiliser le langage ensembliste. Nous emploierons souvent le terme de \emph{plongement} pour désigner les bons monos.}, en privilégiant les $s$-équi\-va\-lences.

Lorsque $\mbf{D}=\mbf{Dif}$, on peut dans certains cas particuliers remplacer les submersions surjectives par les submersions injectives (plongements ouverts); c'est cette version qui a été utilisée dans 
la conférence d'A. Weinstein (voir la note de bas de page de \ref{deq}, et ci-après \ref{eqmor}, où l'on esquisse le passage d'un point de vue à l'autre, avec renversement du sens des flèches d'équivalences).

Ces deux présentations correspondent à différents choix privilégiés de \emph{re\-pré\-sen\-tants} d'une même flèche de la \emph{catégorie de fractions} \merod\!, dans laquelle on rend \emph{inversibles} toutes les $s$-équivalences, et \emph{a posteriori} toutes les \textbf{D}-équivalences, ainsi transformées en \emph{isomorphismes}.

Ces flèches seront appelées \textbf{D}-\emph{méromorphismes}
\,\footnote{Dans le cadre topologique, ce sont essentiellement les morphismes généralisés introduits indépendamment, dans des langages et des cadres très différents, par W.T. van Est et A. Haefliger (s'inspirant de G. Skandalis). Toutefois les conditions de locale trivialité topologique imposées par ce dernier ne paraissent pas se prêter à une définition générale de la composition de ces morphismes topologiques.}
 pour rappeler qu'elles ont un dénominateur, et \textbf{D}-\emph{holomorphismes} lorsque ce dénominateur est trivial; il s'avère \emph{a posteriori} que l'holomorphisme $[f]$ défini par un morphisme de \textbf{D}-groupoïdes $f$ s'identifie à la classe des morphismes qui lui sont \emph{isomorphes} (au sens des transformations naturelles de foncteurs [McL])
\,\footnote{Appelée parfois \emph{classe d'homotopie} de $f$. Dans le cas particulier où $f$ est un automorphisme de \emph{groupe}, $[f]$ n'est autre que la classe de \emph{conjugaison} de $f$, parfois appelée \emph{automorphisme extérieur}; ceci suggère pour les holomorphismes une autre terminologie, à savoir \emph{exomorphismes}.} 
. Nous notons \holod la sous-catégorie des \textbf{D}-holomorphismes.

Il s'avère que, pour définir la composition de ces nouvelles flèches et les interpréter comme flèches d'une catégorie de fractions, il est bien préférable, plutôt que de chercher à travailler avec des représentants privilégiés (\cad définis à isomorphisme près), d'utiliser une classe de représentants, ce qui est l'analogue exact de la construction des rationnels utilisant l'équivalence engendrée par: $(p,q)\sim(pr,qr)$: le produit de deux fractions irréductibles ne fournit pas directement le représentant irréductible!

Cette construction est infiniment plus simple que celle du calcul général des fractions décrite par Gabriel-Zisman [GZ], d'ailleurs ici inapplicable, et où l'on doit considérer des équivalences de diagrammes plus compliquées. Un représentant privilégié, avec son \emph{interprétation géométrique} en termes de lois d'action, est alors l'analogue exact d'un \emph{représentant irréductible}. Toutefois une différence importante, exigeant quelque précaution, est qu'ici la fraction $p/1$ \emph{ne sera pas définie}, et l'image de $p$ est en fait $[p]$ et non $p/1$.

Les propriétés qu'il conviendra d'imposer au couple $(p,q)$, dont le dénominateur $q$ sera ensuite une $s$-équivalence, s'expriment très simplement dans le cadre du diagramme du papillon, apparu dans la section précédente dans un cas particulier.

Ainsi s'expliquent les définitions qui vont suivre, dans lesquelles la notion d'\emph{ex\-ac\-teur} (\ref{exa}) s'avère très utile (bien que d'interprétation géométrique moins simple). Rappelons qu'un exacteur possède un \emph{noyau} en tant qu'objet de \textbf{Grpd(D)}, défini par produit fibré avec la flèche unité $\omega_G$ (où $G$ est son but).

\subsection{Fractions et méromorphismes.}
\label{fem}
\subsubsection{Conditions de transversalité.}
\label{transv}

Nous avons déjà rencontré des cas particuliers de la situation suivante.

Soient $K$ un \textbf{D}-groupoïde de base $E$, et $M,\,N$ deux sous-groupoïdes (\cad munis de deux bons monos $i,\,j$ de $M,\,N$ vers $K$) unifères (\cad de même base $E$), et notons $L$ le produit fibré de $\alpha_M,\,\alpha_N$, qui est plongé dans $\Delta K$.

Nous dirons que $M$ et $N$ sont \emph{transversaux} (resp. \emph{transverses}) dans $K$ si la restriction de $\delta_K$ à $L$ est un bon épi (resp. un isomorphisme). Ceci sera noté $M\trl N$ (resp.$M\tre N$). Ceci implique que $S=M\cap N$ est un sous-groupoïde (plongé) de $M$ et de $N$, qui sera en particulier principal dès qu'il en est ainsi pour l'un d'eux. On pourra alors diviser par $S$ (nos définitions permettent de le faire dans le cadre structuré) pour ramener un couple transversal à un couple transverse.

\smallskip
\paragraph{\emph{Remarque.}} La donnée sur $K$ de deux sous-groupoïdes \emph{tranverses} munit $K$ d'une structure d'\emph{actoïde} (voir \ref{fibasrem}), en général non bi-principal ici, dont la loi de $K$ est la loi mixte.

\smallskip
Supposons maintenant que $N$ soit le noyau d'un exacteur $p:K\rTo^\exa G$. On a alors $M\trl N$ (resp. $M\tre N$) si et seulement si $u=pi$ est un exacteur (resp. un acteur). 

Par suite, si $M$ est aussi le noyau d'un exacteur $q:K\rTo^\exa H$, on aura (en posant $v=qj$):
\begin{center}
(u est un [ex]acteur) $\Leftrightarrow$ (v est un [ex]acteur).
\end{center}
On dira alors que les [ex]acteurs $p$ et $q$ sont \emph{cotransverses} [resp. \emph{cotransversaux}].

La situation est résumée par les diagrammes de papillon suivants:
\begin{center}
\begin{diagram}[w=2.5em,h=1.5em,tight,labelstyle=\scriptstyle,inline]
N&&&&M\\
&\rdiTo\til{j}&&\ldiTo\til{i}&\\
\dTo^\exa\til{v}&&K&&\dTo_\exa\til{u}\\
&\ldTo^\exa\til{q}&&\rdTo^\exa\til{p}&\\
H&&&&G\\
\end{diagram}\hfil
\begin{diagram}[w=2.5em,h=1.5em,tight,labelstyle=\scriptstyle,inline]
N&&&&M\\
&\rdiTo\til{j}&&\ldiTo\til{i}&\\
\dTo^\act\til{v}&&K&&\dTo_\act\til{u}\\
&\ldTo^\exa\til{q}&&\rdTo^\exa\til{p}&\\
H&&&&G\\
\end{diagram}
\end{center}
avec $N=\text{Ker}p,\,M=\text{Ker}q$.

\subsubsection{Fractions.}\label{fractions}
On considère maintenant des couples d'exacteurs $(p,\,q)$ de même source, de la forme:
$$H\lTo^\exa K\rTo^\exa G$$
sur lesquels on fait opérer à la source des exacteurs de la forme $k:K'\rTo^\exa K.$

On dira que les couples $(p,q)$ et $(p',q')$ sont \emph{isomorphes} s'il existe un \emph{isomorphisme} $k$ tel que l'on ait $p'=pk,\,q'=qk$. La classe d'\emph{isomorphie} sera dite la \emph{fraction} définie par le couple $(p,q)$, et notée \fbox{$p/q$}.

On dira qu'ils sont \emph{équivalents} s'il existe un couple de \emph{s-équivalences} $(k,k')$ telles que l'on ait: $pk=p'k',\,qk=q'k'$ (le fait que ce soit bien une relation d'équivalence résulte des propriétés de la catégorie des $s$-équivalences). La classe d'\emph{équivalence} de $(p,q)$ est notée \fbox{$pq^{-1}$}, et on dit qu'elle est \emph{représentée} par la \emph{fraction} $p/q$.

\subsubsection{Méromorphismes.}
La \emph{cotransversalité} du couple $(p,q)$ est invariante par \emph{équi\-va\-lence}; c'est donc une propriété de la classe $pq^{-1}$. La \emph{cotransversité} est invariante par \emph{isomorphie}, et traduit une propriété de la \emph{fraction} $p/q$, qui est alors dite \emph{irréductible}. Toute classe \emph{cotransversale} possède un \emph{représentant irréductible} $p/q$ (bien défini en tant que \emph{fraction}).

La propriété pour le numérateur $p$ ou le dénominateur $q$ d'être une $s$-équivalence est également invariante par équivalence.

\smallskip

Ceci permet de dire que la \emph{classe d'équivalence} $pq^{-1}$ est un \emph{méromorphisme} si elle est cotransversale et si ses dénominateurs sont des $s$-équivalences. Ses représentants peuvent donc être décrits par des papillons du type ci-dessous (qui décrivent un méromorphisme $pq^{-1}$ de$H$ vers $G$), celui de droite (défini à isomorphisme près) étant le \emph{représentant irréductible}:
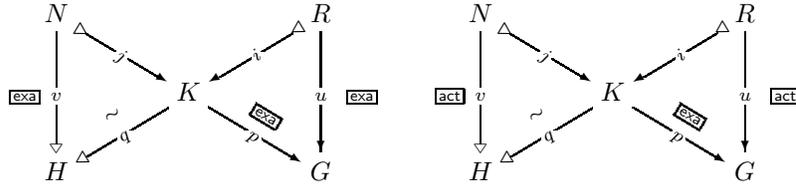
\begin{center}
\begin{diagram}[w=2.5em,h=1.5em,tight,labelstyle=\scriptstyle,inline]
N&&&&R\\
&\rdiTo\til{j}&&\ldiTo\til{i}&\\
\dsTo^\exa\til{v}&&K&&\dTo_\exa\til{u}\\
&\ldsTo\til{q}^\sim&&\rdTo^\exa\til{p}&\\
H&&&&G\\
\end{diagram}\hfil
\begin{diagram}[w=2.5em,h=1.5em,tight,labelstyle=\scriptstyle,inline]
N&&&&R\\
&\rdiTo\til{j}&&\ldiTo\til{i}&\\
\dsTo^\act\til{v}&&K&&\dTo_\act\til{u}\\
&\ldsTo\til{q}^\sim&&\rdTo^\exa\til{p}&\\
H&&&&G\\
\end{diagram}
\end{center}
avec $N=\text{Ker}p,\text{ et }R=\text{Ker}q$ \emph{principal}. La diagonale ascendante est une suite exacte. Dans le diagramme de droite, on a de plus $R\cap N$ \emph{nul}. Comme on l'a dit, on se ramène à ce cas en divisant par le sous-groupoïde principal $S=R\cap N$.

Les méromorphismes seront notés avec des tirets: $H\rDashTo G$.

Interprétant les acteurs $u,\,v$ d'un \emph{représentant irréductible} en termes de lois d'action (\ref{ala}), on retrouve la \emph{description géométrique} donnée par Haefliger dans le cadre topologique [Ast. 116] (noter qu'ici les conditions de locale trivialité sont des consé\-quen\-ces automatiques des données structurées): deux lois d'action de groupoïdes sur une bi-fibration au-dessus de leurs bases respectives, commutant entre elles, l'une d'elles étant principale.

\subsubsection{Composition.} 

\begin{figure}[!htb]
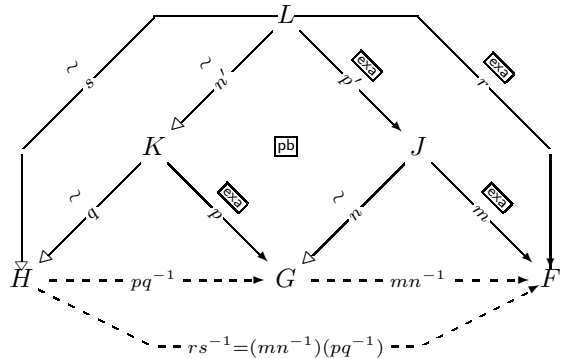

\begin{diagram}[size=2.5em,tight,labelstyle=\scriptstyle]
&&{}&\lLine[abut] &L&\rLine[abut]&{}&&\\
{}&\ldLine(2,2)[abut]^\sim\til s&&\ldsTo(2,2)^\sim\til{n'}&&\rdTo(2,2)^\exa\til{p'}&&\rdLine(2,2)[abut]^\exa\til r&{}\\
{}&&K&&\pb&&J&&{}\\
\dsTo[abut]&\ldsTo(2,2)^\sim\til q&&\rdTo(2,2)^\exa\til p&&\ldsTo(2,2)^\sim\til n&&\rdTo(2,2)^\exa\til m&\dTo[abut]\\
H&&\rDashTo\til{pq^{-1}}&&G&&\rDashTo\til{mn^{-1}}&&F\\
&\rdDashLine(2,1)[abut]&{}&&\rDashLine[abut]\til{rs^{-1}= (mn^{-1})(pq^{-1})}&&{}&\ruDashTo(2,1)[abut]&\\
\end{diagram}
\caption{\label{compomero}Composition des méromorphismes.}
\end{figure}

La composition des méromorphismes est beaucoup plus facile à définir si l'on ne s'astreint pas à travailler exclusivement avec des représentants irréductibles. Elle résulte du diagramme commutatif ci-joint (voir figure \ref{compomero}), que l'on construit en complétant les données $(m,n),(p,q)$ par produit fibré et composition des flèches.

Il faut noter que la composition de deux représentants irréductibles ne donne pas \emph{directement} un représentant irréductible du composé.

Une fois établies les propriétés des sous-catégories des $s$-équivalences et des exacteurs, c'est un simple jeu diagrammatique de s'assurer que la classe d'équivalence du composé ne dépend pas des choix des représentants. Il faut aussi vérifier (toujours de façon diagrammatique) la stabilité de la condition de cotransversalité (et non de cotransversité!) par composition.

On a ainsi construit la catégorie \merod.

\subsubsection{Holographes et holomorphismes.}
Comme nous l'avons annoncé, la définition du foncteur canonique (non fidèle! [McL]):
\begin{center}
$[\,]:\mbf{Grpd(D)}\rightarrow$\,\merod
\end{center}
 n'est pas $(f:H\rightarrow G)\mapsto f/\text{Id}_H$, car le couple $(f,\text{Id}_H)$ ne vérifie pas la condition de transversalité. La définition de la fraction $p/q$ associée à $f$ est donnée par la construction suivante 
\begin{diagram}[size=3em,tight,textflow,labelstyle=\scriptstyle]
&&G\\
&\ruTo^\exa\til p&\usTo_\sim\til{\pi^\uparrow_G}\\
K&\rTo&\square\,G\\
\dsTo_\sim\til q&\pb&\dsTo_\sim\til{\pi^\downarrow_G}\\
H&\rTo\til f&G\\
\end{diagram}
(voir les notations de \ref{papcan}), dont le contenu ensembliste sous-jacent semble connu dans divers contextes et cas particuliers. Le couple $(p,q)$ sera appelé l'\emph{holographe} de $f:H\rightarrow G$. On a $pq^{-1}=[f]$.

La $s$-équivalence $q$ est \emph{scindée} (\cad admet un morphisme section), et $p$ est un \emph{exacteur}. La base de $K$ est le produit fibré $W$ de $\alpha_G$ et de $f^{(0)}$ (notation de \ref{notm}). La fraction $p/q$ est irréductible et elle décrit les actions de $H$ et $G$ sur l'objet $W$, lequel se projette sur les bases de $H$ et de $G$. Elle représente l'holomorphisme défini par $f$, qui n'est autre que la classe d'isomorphie fonctorielle (ou classe d'homotopie) $[f]$ de $f$.

La propriété que le dénominateur soit scindé caractérise les holomorphismes parmi les méromorphismes.

\smallskip

L'holographe de l'identité de $G$ n'est autre que:
$$G\lsTo^{\pi^\downarrow_G}_\sim\square\,G\rsTo^{\pi^\uparrow_G}_\sim G$$
et le papillon associé n'est autre que le papillon canonique de $G$ (\ref{papcan}).
On a:
$$[\pi^\downarrow_G]=[\pi^\uparrow_G]\text{ et }\pi^\uparrow_G(\pi^\downarrow_G)^{-1}=[\text{Id}_G].$$

\smallskip

On montre que le foncteur:
\begin{center}
$[\,]:\mbf{Grpd(D)}\rightarrow$\,\merod
\end{center}
fournit la solution universelle du problème des fractions pour les $s$-équivalences, et \emph{a posteriori} pour toutes les \textbf{D}-équivalences.

La notation globale $pq^{-1}$ que nous avons introduite est \emph{a posteriori} justifiée par l'égalité, dans \merod:
$\boxed{pq^{-1}=[p][q]^{-1}}$.

\subsection{\'{E}quivalences de Morita.}\label{eqmor}

Les méromorphismes inversibles sont caractérisés par la propriété que le numérateur $p$ est aussi une $s$-équivalence. Ce sont les équivalences méromorphes, ou \textbf{D}-\emph{Morita-équivalences}\,\footnote{Lorsque le diptyque \textbf{D} reste fixé, on omettra le préfixe \textbf{D}.}.

Un représentant irréductible est alors décrit par un papillon \emph{symétrique}, du type de ceux rencontrés en \ref{dipap}, bien défini à isomorphisme près, qui fait intervenir deux $s$-équivalences.

Les équivalences holomorphes sont, de façon précise, les classes d'isomorphie des \textbf{D}-équivalences décrites en \ref{deq}.

Lorsqu'on prend pour \textbf{D} le diptyque \textbf{Dif} ou \textbf{Top} il peut exister d'autres re\-pré\-sen\-ta\-tions avec des équivalences \emph{ouvertes} (\ref{eqso}). Nous esquissons le passage d'une re\-pré\-sen\-ta\-tion à l'autre, pour faciliter la comparaison avec les constructions décrites par A. Weinstein au cours de ces \emph{Journées}.

Un représentant irréductible d'une équivalence de Morita est donc donné par un papillon symétrique 
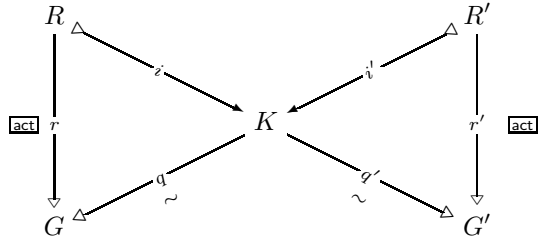
\begin{diagram}[w=4em,h=2em,tight,labelstyle=\scriptstyle]
R&&&&R'\\
&\rdiTo\til{i}&&\ldiTo\til{i'}&\\
\dsTo^\act\til{r}&&K&&\dsTo_\act\til{r'}\\
&\ldsTo_\sim\til{q}&&\rdsTo_\sim\til{q'}&\\
G&&&&G'\\
\end{diagram}
Ce diagramme décrit la bi-fibration principale:
$$B\lsTo^p E\rsTo^{p'}B'$$
où $E$ est la base de $K$.

Dans certains cas particuliers, on peut donner une autre description de l' équi\-va\-lence de Morita (qui existe toujours dans le cadre ensembliste, ou quand certaines flèches sont étales).

Nous supposons que $p$ et $p'$ (et par suite $q$ et $q'$) soient scindées, ce qui permet d'\emph{identifier} $B$ et $B'$ à des sous-variétés de $E$, que nous supposerons de plus \emph{disjointes}, et $G$, $G'$ à des sous -groupoïdes de $K$. Nous noterons $\widehat{B}$ la variété somme de $B$ et $B'$.

Supposons (ce qui est toujours le cas lorsque les sous-variétés $B$, $B'$ sont \emph{ouvertes}) que le \textbf{D}-groupoïde \emph{induit} par $K$ sur $\widehat{B}$ 
soit défini, et notons-le $\widehat{K}$. L'\emph{inclusion} de $\widehat{K}$ dans $K$ est une $i$-\emph{équivalence}.
Le groupoïde $\widehat{K}$ a alors une structure très simple. Il est somme disjointe des quatre ouverts $G$, $G'$, $S$, $S^{-1}$, où l'on a posé $S=\tau_K^{-1}(B'\times B)$, et les inclusions de $G$, $G'$, fournissant un couple de $i$-équivalences:
$$G\riTo^{i_1}_\sim \widehat{K}\liTo^{i_2}_\sim G'.$$

\subsection{Structures orbitales.}\label{storb}
Dans le cadre ensembliste ($\mbf D=\mbf{Ens}$), toutes les surjections sont scindées, ainsi que toutes les équivalences surjectives.

De ce fait la catégorie \merod se confond avec la catégorie \holod, et ses flèches sont simplement les classes d'isomorphie fonctorielle des morphismes (= foncteurs) de groupoïdes, ce qui n'apporte rien d'essentiellement nouveau.

Nous avons vu (\ref{tect} et \ref{doe}) que la classe d'équivalence (ou de Morita-équivalence) d'un groupoïde ensembliste est décrite par une somme de groupes, définis à isomorphisme près, mais que les morphimes ne se laissent nullement décrire en termes de morphismes de groupes.

Dans le cas structuré, on pourra penser les classes de Morita-équivalence de \textbf{D}-groupoïdes comme définissant une structure (en un sens généralisé) algébrico-différentielle très riche, dont deux exemples fondamentaux (\emph{ainsi unifiés}) sont:
\begin{itemize}
	\item une classe d'isomorphie de fibrés principaux, qui correspond à une classe de Morita-équivalence de groupoïdes $s$-transitifs, étudiés dans la section pré\-cé\-dente;
	\item la structure transverse d'un feuilletage régulier, que l'on peut identifier à la classe de \textbf{F}-équivalence (au sens introduit par P. Molino) d'un feuilletage régulier, que nous allons rencontrer dans la section suivante.
\end{itemize}

\section{Feuilletages.}

\subsection{Généralités et notations concernant la catégorie des variétés.}
\label{genot}
Toute variété [VAR] $X$ sera supposée de classe $\textbf{C}^\infty$, de dimension localement finie. L'espace topologique sous-jacent est localement connexe, de Baire, et accessible [TOP], mais non né\-ces\-sai\-rement séparé.

Il y a souvent lieu d'imposer en outre à $X$, en guise de substitut à la paracompacité, ici non disponible en général, une condition supplémentaire de \emph{dénombrabilité}, et parfois de travailler dans une sous-catégorie pleine.

Cette condition portera, tantôt sur $X$ tout entière, tantôt sur les \emph{composantes connexes} de $X$. On renvoie à l'annexe \ref{den} pour l'énoncé des deux choix les plus utiles.

Le plus souvent on considèrera des groupoïdes de Lie dont la \emph{base} est paracompacte; elle est alors séparée, et ses composantes connexes sont à base dénombrable au sens de l'annexe \ref{den}. Mais les groupoïdes d'holonomie seront en général non séparés. Toutefois ils seront \emph{à base dénombrable} dès que la base l'est.
On se place désormais sur la catégorie $\gr D= \textbf{Dif}$ des (morphismes de) variétés.

Elle possède entre autres trois sous-catégories importantes, dont les flèches seront figurées par des notations caractéristiques:
\begin{itemize}
	\item la sous-catégorie $\gr S$ des \emph{submersions}: \hfill$\rsubTo$;
	\item la sous-catégorie $\gr I$ des \emph{immersions}: \hfill$\rimTo$;
	\item la sous-catégorie $\gr E=\gr{S}\cap\gr I$ des \emph{étalements}: \hfill$\retaTo$.
\end{itemize}
Chacune de celles-ci possède à son tour diverses sous-catégories utiles (les dé\-fi\-ni\-tions nouvelles, accompagnées de commentaires, sont reportées dans l'\emph{annexe} \ref{itg}):
\begin{itemize}
	\item $\gr{S}_{\text{sur}}$, surmersions: \hfill$\rsTo$;
	\item $\overline{\gr{S}}_{\text{sur}}$, surmersions propres\,\footnote{D'après un théorème d'Ehresmann, ce sont des fibrations (localement triviales).}:\hfill$\rSTo$;
	\item $\gr{I}_\text{fid}$, \emph{immersions fidèles}:\hfill$\rfidTo$;
	\item $\gr{I}_\text{str}$, immersions \emph{strictes} :\hfill$\ristTo$;
		\item $\gr{I}_\text{str}$, \emph{plongements faibles}:\hfill$\rplfTo$;
	\item $\gr J$, plongements: \hfill$\riTo$;
	\item $\overline{\gr{J}}$, plongements propres: \hfill$\rITo$;
	\item $\gr{E}_{\text{inj}}$ plongements ouverts, ou étales (= étalements injectifs):\hfill$\rploTo$;
	\item $\overline{\gr{E}}_{\text{inj}}$, plongements ouverts fermés:\hfill$\rpofTo$;
	\item $\gr{E}_{\text{sur}}$, surétalements (ou surmersions étales): \hfill$\rsTo\til\etl$;
	\item $\overline{\gr{E}}_{\text{sur}}$, surétalements propres (= revêtements finis):\hfill$\rSTo\til\etl$.
\end{itemize}

\bigskip

\paragraph{\emph{Groupoïdes étales et pseudogroupes}.}\label{geps}
Un groupoïde de Lie \textbf{G} sera dit \emph{étale}\,\footnote{S-atlas dans la terminologie de W.T. van Est.} si $\alpha_G$ est étale. Ses $\alpha$-fibres sont alors discrètes, et le transiteur $\tau_G$ est un immersion. On prendra garde que, même si $\textbf{G}$ est (ensemblistement) transitif, il ne pourra pas être $s$-transitif. La base $\omega_G(B)$ est ouverte.

Tout \emph{pseudogroupe} de difféomorphismes locaux d'une variété $B$ (toujours supposé complet, quitte à le compléter, en tant que préfaisceau) sera identifié au faisceau de ses germes, et considéré ainsi comme un groupoïde étale. 

\subsection{Définition des feuilletages de \St\!.}\label{fst}
Par \emph{préfeuil\-letage} (dif\-fé\-ren\-tia\-ble) $(X,X')$ d'une variété $X$, nous entendrons ici la donnée sur l'ensemble support de $X$ d'une deuxième structure de variété $X'$ telle que l'application identique:
$$\iota_X:X'\longrightarrow X$$
soit une une immersion (nécessairement fidèle puisqu'injective).
On dit que $X'$ est la \emph{structure fine} du préfeuilletage.
On dit que le préfeuilletage $(X,X^{\prime\prime})$ est \emph{plus fin} que $(X,X')$ si $(X^{\prime},X^{\prime\prime})$ est un pré\-feuil\-le\-tage.

La \emph{feuille} de $a\in X$, notée $X'_a$, est la composante connexe de $a$ pour la structure fine $X'$.

On dira que le préfeuilletage est \emph{strict} si toutes ses feuilles sont faiblement plongées (annexe \ref{itg}). Ceci \emph{implique} que l'immersion $\iota_X$ est stricte. Inversement cette dernière condition impliquerait seulement que toutes les feuilles soient strictement immergées.

On note $(\dot{X},X)$ le préfeuilletage \emph{discret}; $(X,X)$ est le préfeuilletage \emph{grossier}.

On définit les \emph{opérations élémentaires} sur les préfeuilletages: produit, pré\-feuil\-le\-tage $U'$ induit sur un ouvert $U$ de $X$, et plus généralement le long d'un morphisme transversal [VAR] à $X'$.

Le préfeuilletage est dit \emph{élémentaire} si $(X,X')$ est isomorphe au produit d'un préfeuilletage discret $(Q,\dot{Q})$ par un préfeuilletage grossier $(F,F)$ (où $F$ est une variété connexe).

\smallskip

Un \emph{feuilletage régulier} est un préfeuilletage localement élémentaire (en un sens évident).

Toute surmersion $X\rsTo Q$ munit $X$ d'un feuilletage régulier $X'$ dit \emph{simple} induit par le feuilletage discret de $Q$. L'espace des feuilles est alors une variété $\tilde Q$ (en général non séparée, même quand $X$ et $Q$ le sont), étalée sur $Q$.

Pour tout groupoïde de Lie \textbf{G}, on notera $G^\alpha$ et $G^\beta$ les feuilletages simples ainsi associés aux surmersions $\alpha$ et $\beta$. Lorsque $G$ est étale, ces feuilletages sont discrets.

\smallskip

Un \emph{feuilletage de \St}est un préfeuilletage \og semi-continu inférieurement \fg\, au sens suivant:\\
pour tout $a\in X$ il existe un voisinage ouvert $U$ de $a$ dans $X$, et un feuilletage \emph{régulier} $(U,U^{\prime\prime})$ tel que:
\begin{enumerate}
	\item $(U,U^{\prime\prime})$ est \emph{plus fin} que le feuilletage \emph{induit} $(U,U')$;
	\item les feuilles $U'_a$ et $U^{\prime\prime}_a$ sont les mêmes.
\end{enumerate}

On peut choisir $(U,U^{\prime\prime})$ élémentaire, isomorphe à $(T,\dot{T})\times(F,F)$ (où $F=U'_a=U^{\prime\prime}_a$ est la feuille singulière) et on voit que le feuilletage induit sur $U$ est isomorphe au produit, par le feuilletage grossier de $F$, du feuilletage de \St $T'$ induit sur la sous-variété transverse $T\times\{a\}$ (dont la feuille singulière est réduite au point $a$).

C'est ce que suggère le schéma \ref{fs} ci-joint, dans lequel le feuilletage régulier, incomplètement dessiné, aurait ses feuilles verticales, parallèles à la feuille singulière.
\begin{figure}
\fbox{
\unitlength=.4em
\newsavebox{\Koupe}
\begin{picture}(45,55)
\thicklines
\savebox{\Koupe}(20,8)[h]
{\qbezier(-2,1)(0,1)(2,0)
\qbezier(2,0)(4,-1)(2,-1)
\qbezier(-4,2)(0,2)(4,0)
\qbezier(4,0)(8,-2)(4,-2)
\qbezier(-6,3)(0,3)(6,0)
\qbezier(6,0)(12,-3)(6,-3)
\qbezier(-8,4)(0,4)(8,0)
\qbezier(8,0)(16,-4)(8,-4)
\put(-.65,-.75){$\bullet$}}
\put(19,4){\usebox{\Koupe}}
\put(19,44){\usebox{\Koupe}}
\put(19,24){\usebox{\Koupe}}
\multiput(11,12)(2,-1){9}{\line(0,1){40}}
\put(31.2,5){\line(0,1){40}}
\put(42,23.5){\line(-2,1){10}}
\put(0,32.5){\line(1,0){10}}
\put(0,32.5){\line(2,-1){18}}
\put(18,23.5){\line(1,0){24}}
\qbezier[20](11,32.5)(17.5,32.5)(24,32.5)
\qbezier[12](24,32.5)(28,30.5)(31,29)
\end{picture}}
\caption{\label{fs}Feuilletages de $\check{\textrm{S}}$tefan.}
\end{figure}

\smallskip

Toute feuille d'un feuilletage de \St qui est à base dénombrable au sens de l'annexe \ref{den} est faiblement plongée (\ref{pfis}) (l'argument de C. Chevalley [LG] reste utilisable). Il en est toujours ainsi lorsque la variété ambiante est \emph{paracompacte} (ce qui implique séparée). Sous cette hypothèse on peut s'assurer que la formulation simple que nous donnons équivaut à l'ensemble des conditions, d'apparence plus technique, des définitions originales introduites par H. Sussmann, et, sous une forme plus précise, par P. \St.

\smallskip
Nous abordons maintenant l'étude des liens entre feuilletages de \St et grou\-poïdes de Lie.

\section{Feuilletages de \St attachés à un groupoïde de Lie.}\label{fagl}

Nous commencerons par une mise en garde contre un risque de confusion.

Bien qu'un pseudogroupe puisse être considéré comme un cas particulier de groupoïde (étale) (\ref{geps}), la construction rappelée (uniquement pour mémoire ) ci-après en \ref{faps} n'est aucunement un cas particulier de celle décrite ensuite en \ref{fogl} (la seule qui nous concerne ici), laquelle, appliquée à un pseudogroupe, donnerait toujours le feuilletage discret, et n'est donc pas intéressante dans ce cas.

\smallskip

\paragraph{\emph{Exemple.}} Soit $G$ un groupe de Lie opérant sur une variété $X$, et soit $H$ le grou\-poïde d'action correspondant (\ref{ala}). Notons $G'$ le groupe discret défini par $G$. Le grou\-poïde d'action $H'$ de $G'$ est un pseudogroupe (\ref{geps}) qui définit un feuilletage simple (et même élémentaire) de $H$.

Alors la construction de Sussmann-\St appliquée à $H'$ donnerait le même résultat que la construction \ref{fogl} appliquée à $H$, tandis que cette dernière construction appliquée à $H'$ donnerait la structure discrète de $H$.

\subsection{Construction de Sussmann-\St pour un pseudogroupe.}\label{faps}
Il faut rappeler qu'originellement la notion de feuilletage de \St a été introduite in\-dé\-pen\-dam\-ment par Sussmann et \St comme un feuilletage (singulier) de la base, attaché canoniquement à un pseudogroupe, tel que les orbites du pseudogroupe (composantes transitives de la base) soient des réunions de feuilles (c'est le moins fin satisfaisant à cette condition).

\smallskip

\paragraph{\textbf{Commentaire}.}Répétons qu'il s'agit d'une construction différente de celle du paragraphe suivant. Mais il pourrait être du plus haut intérêt d'examiner s'il serait possible, en étendant la construction de Sussmann-\St à un groupoïde de Lie \emph{non étale}, de traiter le problème universel (\emph{dual de celui considéré ci-après}) consistant à construire une flèche $G\rightarrow G^{\text{T}}$, universelle pour les morphismes de $G$ vers les groupoïdes de Lie localement $s$-transitifs\ldots

\subsection{Feuilletage de \St orbital pour un groupoïde de Lie.}\label{fogl}
 Nous renvoyons ici à [P 1985] pour des énoncés plus précis et des démonstrations.

Nous examinons ici ce qui subsiste, pour les groupoïdes de Lie, de la dé\-compo\-si\-tion orbitale ensembliste en composantes transitives, décrite en \ref{doe}.
Pour tout groupoïde de Lie \textbf{G}, on notera $G^\alpha$ et $G^\beta$ les feuilletages simples associés aux surmersions $\alpha_G$ et $\beta_G$. (Lorsque $G$ est étale, ces feuilletages sont discrets).

Si l'on ne fait pas d'hypothèses de régularité sur le transiteur $\tau_G$, les fibres de $\tau_G$ ne définissent pas en général un feuilletage régulier de $G$, mais on montre (\emph{loc. cit.}, Prop. B 6) qu'elles définissent un feuilletage \emph{régulier} $G^\tau$ raffinant à la fois  $G^\alpha$ et $G^\beta$ 

En sens inverse $G$ est muni canoniquement d'un feuilletage de \St $G^{\text t}$, engendré par les feuilletages réguliers $G^{\alpha}$ et $G^{\beta}$. On a ainsi la chaîne de raffinements de feuilletages de $G$ (réguliers ou de \St selon le cas):
$$G^{\tau}\prec G^\alpha(G^\beta)\prec G^{\text t}\prec G.$$

Chaque composante transitive de $G$ est un sous groupoïde ouvert fermé de $G^{\text t}$, qui est $s$-transitif, et dont la structure a donc été étudiée en \ref{fpgt}: pour tout $e\in B$, la $\alpha$-fibre $G_e$ possède une structure de \emph{fibré principal}, dont le \emph{groupe structural} est le groupe d'isotropie de $e$, et dont le \emph{groupoïde de jauge} est la composante transitive de $e$.

L'application identique
$$\iota^{\text t} :G^{\text t}\rightarrow G$$
définissant le feuilletage de \St canonique factorise universellement les morphismes des groupoïdes de Lie $s$-transitifs vers $G$.

Lorsque $B$ est paracompacte, ce feuilletage est strict (\ref{fst}).

\smallskip

\paragraph{\emph{Exemple (trivial)}.} Si l'on applique la construction au groupoïde décrivant l'action, sur la droite, du groupe multiplicatif $\mbb{R}^{+\ast}$, on obtient la décomposition de la droite en deux demi-droites ouvertes et l'origine, qui est bien un feuilletage de \St à trois feuilles. On peut répéter dans ce cas particulier ce qui a été dit à propos de l'exemple dans l'introduction de \ref{fagl}.

\smallskip

\paragraph{\emph{Composante neutre.}} Pour un groupoïde de Lie \textbf{G} (plus généralement un groupoïde topologique, pour lequel, d'après notre définition générale des groupoïdes \textbf{D}-struc\-tu\-rés, et du diptyque \textbf{Top}, on suppose toujours $\alpha_G$ \emph{ouverte}), le rôle de la composante neutre d'un groupe est ici joué par la réunion des composantes $\alpha$(ou $\beta$)-connexes des unités. C'est un sous-groupoïde \emph{ouvert}, mais qui n'est en général \emph{ni fermé ni invariant} dans $G$ (cette différence cruciale avec le cas des groupes est source de difficultés, mais de beaucoup de richesse). On n'a donc pas en général l'analogue du groupe discret quotient.

On note cependant que, lorsque $G$ est $s$-transitif (ceci s'étend au cas où le transiteur est régulier), on peut toujours diviser par la composante neutre du plurigroupe d'isotropie (qui existe alors en tant qu'objet structuré, et est invariant dans $G$). 
Le groupoïde quotient ainsi obtenu est alors \emph{de Galois}, ce qui signifie que le transiteur $\tau_G$ est \emph{étale surjectif}. Son isotropie est alors discrète, et les $\alpha$-fibres sont des \emph{revêtements galoisiens} de la base (d'où la terminologie).

\bigskip
On notera que, lorsque $G$ est $\alpha$-\emph{connexe} (\cad coïncide avec sa composante neutre), la décomposition orbitale en composantes transitives (lesquelles sont $s$-transitives) coïncide avec la décomposition en feuilles de $G^{\text t}$.

\section{Groupoïdes de Lie attachés à un feuilletage de \St.}
D'après G. Hector (communication privée, non publiée), tout feuilletage de \St d'une variété $B$ provient (par la construction qui précède) d'un groupoïde de Lie\,\footnote{Il conviendrait peut-être d'imposer ici à $B$ d'être paracompacte, ou tout au moins de satisfaire à une hypothèse de dénombrabilité.}, auquel nous imposerons (quitte à prendre la composante neutre) d'être $\alpha$-\emph{connexe}. On peut dire qu'un tel groupoïde \emph{intègre} le feuilletage de \St\!.

Mais il n'y a évidemment pas unicité, et ceci ouvre l'étude de la recherche d'un ou plusieurs choix canoniques, en imposant à ceux-ci des conditions supplémentaires, notamment des conditions universelles.

Il n'y a peut-être pas de réponse satisfaisante dans le cas le plus général (nous y reviendrons dans les pistes de réflexions ci-dessous), mais nous allons décrire des cas particuliers assez étendus, et très importants pour beaucoup d'applications, où l'on peut donner des réponses précises et très exhaustives.

\subsection{Cas des feuilletages simples.} (\ref{fst})
Il y a ici évidemment un choix canonique et minimal, donné par la composante neutre $G^{\text c}$ du graphe $G$ de la relation d'équivalence \emph{régulière} définie par la surmersion $B\rsTo Q$. C'est un groupoïde \emph{de Godement}.

On notera que, si $Q$ est \emph{séparée}\,, le transiteur $\tau_G$ est un plongement \emph{propre}, mais, comme nous l'avons déjà signalé, on peut seulement dire que $\tau_{G^{\text c}}$ est un \emph{plongement}. Le théorème de Godement reste valable dans ce cadre, si l'on accepte pour quotient une variété $\widetilde{Q}$ \emph{non séparée}. C'est bien sûr ce qui explique nos choix pour les structures de diptyques (de Godement) sur \textbf{Dif} et \textbf{DifHaus} (\ref{dip}).

\subsection{Cas des feuilletages réguliers: graphoïdes.} Nous référons ici à [P 1966] et [P 1984], ainsi qu'à l'annexe \ref{ght}.
\subsubsection{Holonomie d'un groupoïde microdifférentiable.}
Appliquée dans le cas d'un feuilletage différentiable \emph{régulier}, cas dans lequel la condition de locale simplicité topologique est bien vérifiée, on peut voir que la topologie définie par Ehresmann sur le groupoïde d'holonomie topologique (\ref{ght}) est en fait une topologie de variété.

Dans [P 1966] nous avons interprété cette situation dans le cadre plus général (nécessité par l'étude de la correspondance sous-algébroïde de Lie $\rightarrow$ groupoïde de Lie) d'un groupoïde algébrique $G$ (ici le graphe \emph{ensembliste} de la relation d'équi\-va\-ence des feuilles) \emph{localement différentiable} (ou morceau de groupoïde différentiable), \cad engendré par une partie munie d'une structure de variété $W$ sur laquelle la loi de composition partielle induite est compatible (en un sens précis) avec la structure de variété, le germe de cette structure différentiable étant bien défini (ceci résulte ici du fait que le feuilletage est \emph{localement élémentaire}\,\footnote{Condition beaucoup plus précise que la condition de \emph{locale simplicité topologique} d'Ehresmann, les deux conditions n'étant équivalentes que dans le cadre différentiable.}). La donnée d'un tel \emph{germe} de structure différentiable sur un groupoïde algébrique définit une structure de groupoïde \emph{microdifférentiable}.

Dans le cas classique des groupes (voir [LIE]), les translations (globales) permettent de globaliser la structure différentiable pour obtenir un groupe de Lie global.

Dans le cas général on dispose seulement des \emph{translations locales}\,\footnote{Associées aux bi-sections locales de $\alpha,\beta$, dont la composition a été introduite par Ehresmann. C. Albert les utilise sous le nom de \emph{glissements}.}, et il n'existera pas en général de structure de groupoïde de Lie global de support ensembliste $G$ admettant $W$ pour sous-variété ouverte.

On peut cependant, grâce à ces translations locales, en s'inspirant de la construction du groupoïde d'holonomie d'un feuilletage régulier, caractériser complètement la classe des groupoïdes de Lie\,\footnote{Ici la \emph{base} des groupoïdes est supposée \emph{paracompacte}.} (globaux) $K$ qui sont \emph{localement} isomorphes à un groupoïde localement différentiable (et, en passant aux germes, ce qui est possible grâce à la locale connexité, microdifférentiable) donné\,\footnote{Notre construction a inspiré à Ronnie Brown diverses généralisations dans le cadre topologique.}. Ce qui suit s'applique comme \emph{cas particulier} au graphe (ensembliste) de la relation d'équivalence associée à un feuilletage régulier, muni de sa structure localement différentiable, ou microdifférentiable définie à partir d'un recouvrement par des ouverts élémentaires.

Le problème admet une solution minimale universelle $H$, le \emph{groupoïde d'holonomie} du groupoïde microdifférentiable, et une solution maximale $M$ (que nous avions appelé le \emph{groupoïde de monodromie})\,\footnote{Dans le cas particulier des feuilletages, cette notion a été redécouverte beaucoup plus tard indépendamment par de nombreux auteurs, et appelée parfois aujourd'hui groupoïde d'homotopie.}, dont les $\alpha$-fibres sont les revêtements universels de celles du précédent, et qui résout un autre problème universel, celui de la globalisation des morphismes \emph{localement définis}\,\footnote{Le groupoïde d'holonomie d'un feuilletage différentiable régulier a été redécouvert, grâce à une tout autre construction, par H.E. Winkelnkemper, en 1982, sous le nom de graphe du feuilletage, et utilisé par A. Connes [NCG].}.

Toute autre solution $K$ donne lieu à une factorisation par deux morphismes étales:
$$M\retasTo K\retasTo H.$$

\bigskip
\subsubsection{Graphoïde universel.} La référence est ici [P 1984].

Le paragraphe précédent, bien que s'appliquant notamment aux graphes (ensemblistes) des feuilletages réguliers, n'utilise aucune aucune spécificité de ces derniers parmi les groupoïdes (localement différentiables).

Dans [P 1984] nous nous sommes au contraire préoccupé de caractériser la spécificité des groupoïdes d'holonomie des feuilletages \emph{réguliers} parmi les groupoïdes de Lie, ce qui nous a conduit à étudier plus généralement ceux pour lesquels le transiteur $\tau_G:G\rightarrow{B\times B}$ est une \emph{subimmersion}. Il s'avère en effet \emph{a posteriori }(ce qui n'était nullement évident \emph{a priori}), que la factorisation canonique de $\tau_G$ dans \textbf{Dif} (voir\ref{ifid})
$$G\rsTo H\rfidTo B\times B$$
définit en réalité une factorisation en morphismes de groupoïdes de Lie.

Cela résulte en fait de l'étude des groupoïdes de Lie $H$ dont le transiteur $\tau_H$ est une \emph{immersion fidèle}, pour lesquels nous avons introduit le terme, qui paraît très naturel, de \emph{graphoïdes}.

Nous avons montré que tout graphoïde de base $B$ fixée s'identifie à un sous-groupoïde ouvert d'un \emph{graphoïde universel} ${\mbf H}(B)$\,\footnote{Dans cette étape il faut évidemment renoncer à toute hypothèse de séparation et de dénombrabilité pour $H$, même lorsque $B$ vérifie de telles hypothèses. Toutefois, si $B$ est à base dénombrable, il en sera de même de la composante neutre $\mbf{H^c}(B)$, qui pourrait être appelée \emph{groupoïde d'holonomie universel de $B$}.}.

La base de $\mbf{H}(B)$ est la \emph{variété} $\Phi(B)$ \emph{des germes de feuilletages réguliers} de $B$, qui est étalée sur $B$.

$\mbf{H}(B)$ est lui-même un ouvert de la variété $\text{J}(B\times B)$ des germes de sous-variétés de $B\times B$. Le feuilletage orbital $\mbf{H'}(B)$ de $\mbf{H}(B)$ est \emph{régulier}.

Il y aurait lieu d'introduire aussi le groupoïde $\mbf{M}(B)$ dont les $\alpha$-fibres sont les revêtements universels de celles de $\mbf{H}(B)$.

La donnée d'un feuilletage \emph{régulier} de $B$ s'identifie à celle d'un \emph{section} de l'étalement canonique $\varphi_B$ de $\Phi(B)$ sur $B$. Son groupoïde d'holonomie est induit par le groupoïde d'holonomie universel.

\smallskip
\subsubsection{Caractérisations des pseudogroupes.} Si l'on abandonne ici l'$\alpha$-connexité, on peut caractériser les pseudogroupes parmi les groupoïdes étales par la condition d'être en outre des graphoïdes.

En fait on peut inversement caractériser les graphoïdes par la condition d'être \emph{Morita-équivalents} à des pseudogroupes. On retrouve ici les deux aspects classiques de l'holonomie (dont beaucoup d'auteurs privilégient l'un des deux):
\begin{itemize}
	\item groupoïde d'holonomie (qui est canoniquement attaché au feuilletage);
	\item pseudogroupe d'holonomie, attaché au choix d'un variété trasverse totale.
\end{itemize}
Tous ces choix appartiennent à la même classe de Morita-équivalence, dont on privilégie, selon le cas, les représentants $\alpha$-connexes ou les représentants étales.

\smallskip
\subsubsection{Groupoïdes de Barre.} Un cas particulier intéressant, \emph{intermédiaire} entre celui des groupoïdes de Godement (le transiteur est un plongement) et celui des graphoïdes (le transiteur est une immersion fidèle, voir\ref{ifid}), est celui où le transiteur $\tau_G$ est un \emph{plongement faible} (\ref{pfis}), propriété qui est encore invariante par la Morita-équivalence: $\tau_G:G\rplfTo B\times B$.

Nous dirons alors que le groupoïde est \emph{de Barre}. Elle caractérise en effet les groupoïdes de Lie dont l'espace orbital est une \textbf{Q}-\emph{variété de Barre}, comme cela résulte implicitement de la thèse de R. Barre [RB].

\bigskip
\subsection{Commentaires.} La propriété remarquable et fondamentale d'un graphoïde $G$ est que sa loi de composition est uniquement déterminée par la seule donnée de l'immersion fidèle $\tau_G:G\rfidTo B\times B$, ou encore du couple de surmersions $\alpha_G,\beta_G:G\rsTo B$, autrement dit par la seule donnée du \emph{graphe} (différentiable régulier) \emph{sous-jacent}. Elle est d'ailleurs induite par la composition des germes de graphes.

Bien entendu ces germes de graphes (réguliers) ne sont pas quelconques. On peut les interpréter comme des germes d'isomorphismes entre des germes de quotients locaux de $B$. On peut aussi les caractériser par une propriété algébrique qui traduit l'existence d'un 
germe de règle de trois $zy^{-1}x$. Ils sont appelés \emph{germes d'isonomie} dans [P 1985].

La théorie des graphoïdes est à notre sens un des moments cruciaux (annoncés dans l'introduction) où l'interaction de la structure (al\-gé\-bri\-que\-ment triviale) de \emph{graphe} (d'équivalence) avec celle de variété fait naître des groupoïdes non triviaux, et donc des \emph{groupes}.

C'est aussi un point où la théorie des feuilletages différentiables réguliers s'écarte fondamentalement de celle des feuilletages topologiques localement simples d'Ehresmann, manifestant une propriété cruciale du foncteur d'oubli $\mbf{Dif\rightarrow Top}$ (qui réside en fait dans le théorème d'inversion locale et ses conséquences, via le théorème de Godement).

Un autre point crucial est le suivant:
Nous venons de voir que la théorie des feuilletages réguliers est essentiellement équivalente à celle des graphoïdes, et que plus généralement, les groupoïdes de Lie dont le feuilletage de \St associé est régulier sont ceux pour lesquels le transiteur est une subimmersion (\cad est de rang localement constant).

Il en résulte que les singularités de ce feuilletage, \cad les sauts de dimensions des feuilles, sont reflétées par des sauts du rang du transiteur, ou de la dimension des groupes d'isotropie. En particulier les groupes d'isotropie des feuilles singulières ponctuelles seront non discrets, et seront susceptibles d'apporter une information sur le voisinage d'une telle singularité, alors que le groupe d'holonomie topologique d'une telle feuille est nul, et ne peut porter aucune information sur la nature de la singularité.

\subsection{Quasi-graphoïdes ou convecteurs.} La référence est ici [P 1985].

Nous avons montré, en collaboration avec notre étudiant B. Bigonnet, que la propriété fondamentale de la structure des graphoïdes (à savoir d'être uniquement déterminés par leur structure sous-jacente de graphe) est partagée par une classe plus large et importante de groupoïdes, appelés \emph{quasi-graphoïdes} ou encore \emph{convecteurs}.

Le transiteur va maintenant satisfaire à une forme très affaiblie d'injectivité locale, au niveau des germes.

\subsubsection{Définition des quasi-graphoïdes ou convecteurs.}
Nous noterons ici $S^\lambda G$ le groupoïde \emph{étale} des germes de bi-sections locales de $\alpha_G$, $\beta_G$, avec la loi de composition introduite par Ehresmann. Il est muni d'un morphisme canonique $s^\lambda_G:S^\lambda_G\rightarrow G$, qui est une \emph{immersion}. Celle-ci est toujours surjective lorsque $G$ est $\alpha$-connexe. Nous poserons, avec Ehresmann: $\Pi^\lambda(B)=S^\lambda(B\times B)$ (groupoïde des germes de difféomorphismes locaux de $B$).

Cette construction est fonctorielle vis-à-vis des $B$-morphismes, ce qui permet d'écrire le diagamme commutatif:
\begin{diagram}[w=4em,h=2em,\lab=\sst]
S^\lambda G&\retaTo^{S^{\lambda}(\tau_G)}&S^{\lambda}(B\times B)\\
\dimTo^{s^\lambda_G}&&\dimTo_{s^\lambda_{B\times B}}\\
G&\rTo^{\tau_G}&B\times B\\
\end{diagram}

Avec ces notations, on dira que $G$ est un \emph{quasi-graphoïde} (ou encore un \emph{convecteur}) si $S^{\lambda}(\tau_G)$
est \emph{injective}. Cela \emph{implique} que la restriction de $\tau_G$ à un ouvert dense de $G$ est une immersion. Cette condition peut aussi s'interpréter comme définissant un sous-objet de $B\times B$ dans une catégorie convenable. On notera cette propriété:
$$\tau_G:G\rqimTo B\times B$$
et on pourra dire que $\tau_G$ est une \emph{quasi-immersion}.

Bien qu'il n'existe plus maintenant de germes de quotients locaux, il demeure possible d'organiser les germes de convecteurs en un groupoïde de Lie $\widetilde{\mbf H}(B)$ (contenant le graphoïde universel ${\mbf H}(B)$ comme sous-groupoïde ouvert), qui est universel en ce sens que tous les convecteurs de base $B$ s'identifient à des sous-groupoïdes ouverts de celui-ci.
 La base $\widetilde{\Phi}(B)$ de ce groupoïde est la \emph{variété des germes de convections}, munie d'un étalement $\widetilde{\varphi}_B$ sur $B$. Elle est munie d'un feuilletage de \St qui sera dit \emph{quasi-régulier}. Elle contient comme sous-variété ouverte la variété $\Phi(B)$ des germes de feuilletage réguliers.

Une \emph{convection} sera une section de $\widetilde{\varphi}_B$.

Toutefois, alors que les feuilletages réguliers portent une convection canonique, il n'en est plus de même en général pour les feuilletages de \St quasi-réguliers, et il serait intéressant de caractériser une classe assez large de tels feuilletages pour lesquels on pourrait faire un choix canonique (minimal) d'un graphoïde associé.

\subsubsection{Exemples}
\label{sec:Exqg}
\begin{enumerate}
	\item Le groupoïde d'action associé au groupe des rotations du plan est un convecteur, associé au feuilletage par les cercles concentriques.
	\item Le \emph{groupoïde tangent} d'A. Connes [NCG], que nous notons $\widetilde{\text{T}}B$, est un convecteur. Les feuilles du feuilletage de \St qu'il définit sur $\mbb{R}\times B$ sont (en supposant $B$ connexe) $\{h\}\times B,\,h\neq0$ et les points de $\{0\}\times B$. Rappelons que les composantes transitives de ce groupoïde sont $\{h\}\times(B\times B),\,h\neq0$ et les fibres de $\text{T}B$ pour $h=0$.
	\item Les faces (des diverses dimensions) d'une \emph{variété à bords et à coins} munissent celle-ci d'un feuilletage de \St canonique qui peut être défini par un convecteur. A. Weinstein a attiré l'attention sur les applications de ce type de groupoïdes et des algébroïdes de Lie associés [W]. Ceux-ci ont été étudiés notamment par B. Monthubert.
\end{enumerate}
Dans tous ces exemples, il est intéressant d'étudier dans quelle mesure on peut faire des choix privilégiés pour les convecteurs définissant le feuilletage de \St\!, et quels sont les \emph{automorphismes} de ces convecteurs préservant le feuilletage.

\subsection{Quelques pistes à explorer.} Observons que le feuilletage de \St de $\mbb{R}^3$ par les sphères concentriques est associé au groupoïde d'action du groupe des rotations (ou bien encore de son revêtement universel $S_3$), mais celui-ci n'est plus un quasi-graphoïde. Le groupoïde qu'il induit sur la partie régulière (sur laquelle le feuilletage induit est simple) n'est pas un groupoïde principal, ni même un groupoïde trivial. \emph{Ses composantes transitives sont les groupoïdes de jauge d'une fibration principale de Hopf}.

Cet exemple n'est qu'un cas très particulier des feuilletages riemanniens (singuliers) étudiés par P. Molino.

Ceci montre que, dans l'étude des groupoïdes de Lie ($\alpha$-connexes) intégrant un feuilletage de \St donné, d'une part on ne doit pas imposer en général que le feuilletage induit dans un ouvert régulier coïncide avec le groupoïde d'holonomie, d'autre part cette étude est susceptible d'apporter une information géométrique extrêmement riche sur le voisinage des singularités.

On est ainsi conduit à élargir la notion d'équivalence de Morita, laquelle est engendrée par les $s$-équivalences, en considérant plus généralement l'équivalence \emph{engendrée par les $s$-extenseurs} (\ref{sex}), et éventuellement l'équivalence (au moins \emph{a priori}) plus restrictive obtenue en ne considérant que les $s$-extenseurs \emph{à noyau $\alpha$-connexe}.

Ceci amène à considérer la \emph{catégorie de fractions obtenue en inversant ces $s$-extenseurs}. Or il se trouve que la construction diagrammatique (reposant sur l'usage du papillon, dans le cadre des diptyques) qui est décrite dans [P 1989] subsiste sans changements fondamentaux lorsqu'on remplace les $s$-équivalences par les $s$-extenseurs.

Ceci définit des structures généralisées, que l'on peut considérer comme des structures transverses sur un feuilletage de \St (en général non uniquement déterminées par celui-ci: un tel choix comporte donc en général une donnée sup\-plé\-men\-taire), auxquelles on pourrait étendre le nom de convection.

Bien entendu ces remarques ne dispensent pas, bien au contraire, d'explorer de telles classes d'équivalence, d'en rechercher des représentants privilégiés, ainsi que la possibilité de choix canoniques\ldots

Par ailleurs l'exemple des feuilletages riemanniens singuliers et celui des convecteurs incite à la recherche de résultats généraux reliant les sauts des groupes d'isotropie avec les sauts corrélatifs des variétés transverses et des feuilles (compensation d'un saut vertical par un saut horizontal\ldots). La conjecture d'A. Weinstein serait le prototype de cette sorte de relations\ldots

\subsection{Lien avec l'intégration du non intégrable (A. Weinstein).}\label{lini}

Indiquons le \emph{lien très étroit avec le point de vue développé par A. Weinstein} lors de ces \emph{Journées}.

\subsubsection{Avec des $s$-équivalences (étales).}
Soit $(U_i)_{i\in I}$ un recouvrement ouvert d'une variété $B$. Posons $U_{ij} = U_i\cap U_j$, puis:
$$U=\coprod_{i\in I}U_i,\quad\quad R=\coprod_{(i,j)\in I\times I}U_{ij}.$$
$R$ n'est autre que le graphe (ici un groupoïde étale) de l'équivalence associée à la surmersion (étale) $r:U\rsTo B$, lequel est muni des projections vers $U$: $$U\letasTo^{\beta_R}_{p_1}R\retasTo^{\alpha_R}_{p_2}U.$$

La donnée d'un groupoïde de Lie $G$ de base $B$ induit le long de $r$ un groupoïde de Lie $H$ de base $U$, lequel induit, à son tour, le long de $p_1$ et $p_2$, un \emph{même} groupoïde de Lie $K$ de base $R$ (induit par $G$ le long de la projection $R\rsTo B$), d'où un couple de $s$-équivalences:
$$H\letasTo^\sim_{q_1}K\retasTo^\sim_{q_2}H,$$
qui sont étales et canoniquement scindées.
Lorsque la structure de $H$ provient de l'intégration \emph{locale} d'un algébroïde de Lie $\go{g}$ (de base $B$), on aura seulement \emph{a priori } un \emph{isomorphisme} entre les deux groupoïdes induits $k:K_1\rightarrow K_2$, traduisant les données de recollement, mais la \emph{classe de Morita} de $H$ est intrinsèquement attachée à $\go{g}$. On peut alors considérer que cette classe de Morita intègre $\go{g}$ en un sens généralisé. L'algébroïde de Lie sera intégrable au sens strict si cette classe possède un \emph{représentant de base} $B$. Dans le cas général, elle réalise, au sens où l'entendait A. Weinstein dans le titre de sa conférence, l'intégration d'un algébroïde de Lie non intégrable.

Ce que nous venons de suggérer, en nous inspirant de l'exemple ci-dessus du feuilletage de \St en sphères concentriques, est de remplacer les $s$-équivalences par des $s$-extenseurs à noyaux $\alpha$-connexes\ldots

\subsubsection{Avec des équivalences ouvertes (A. Weinstein).}
Compte tenu du fait que les flèches figurant dans le paragraphe précédent sont étales et scindées, on peut utiliser pour les équivalences de Morita la construction esquissée en \ref{eqmor}, et n'utiliser que des \emph{équivalence ouvertes}. La présentation donnée par A. Weinstein peut se lire dans ce cadre.

\section{\'{E}pilogue.}
Nous serions heureux d'avoir convaincu quelques auditeurs réticents de l'intérêt du concept de groupoïde de Lie d'Ehresmann pour comprendre la dualité graphe/ groupe via le différentiable.
\appendix
\section{Note historique.}
\label{nh}
En réponse à une remarque orale faite par Mlle Libermann concernant la date d'introduction des goupoïdes différentiables par Charles Ehresmann, nous croyons utile, à fin de clarification, d'apporter les quelques précisions qui vont suivre.

Pour plus de détails et beaucoup d'autres références, nous renvoyons au tome I des \OE uvres Complètes [EOC], et en particulier au paragraphe 2 (page 588) du Synopsis présenté par Andrée Charles Ehresmann, ainsi qu'aux Commentaires de l'auteur pp. 526-539.
\par\smallskip

La \emph{définition générale} des groupoïdes différentiables n'apparaît, comme nous l'avions indiqué dans notre présentation lors des \emph{Journées}, qu'à l'occasion du Colloque de Géométrie Différentielle Globale tenu à Bruxelles en 1958. Les Actes de ce Colloque, publiés en 1959 [EOC/50/]\,\footnote{Les /numéros/ réfèrent à la numérotation des publications de C. Ehresmann figurant en tête de ses \OE uvres complètes [EOC]} sans tenir compte des corrections d'épreuves, comportaient une erreur (corrigée dans les \OE uvres Complètes, partie I, pp. 237-250) concernant la propriété de rang localement constant, exigée des applications source et but $\alpha\,,\beta$. La formulation plus simple en termes de surmersions sur la variété base, adoptée ici, qui s'est imposée par la suite comme beaucoup plus maniable et se prêtant aux bonnes généralisations, n'est explicitée que dans les Actes du Congrès de Géométrie tenu à Bologne en 1967 [EOC/103/] (la variété base restant toutefois toujours \emph{identifiée} à une sous-variété fermée\,\footnote{Nous n'imposons pas ici cette dernière condition, laquelle équivaut, lorsqu'on suppose la variété base séparée, à la condition de séparation de la variété support du groupoïde, afin de pouvoir inclure sans restriction l'exemple du groupoïde d'holonomie.} du groupoïde, C.\til Ehresmann s'étant toujours systématiquement refusé à exploiter la désidentification unités/objets, pourtant souvent fort éclairante, et même indispensable).
\smallskip

Cependant les \emph{exemples} et cas particuliers qui ont motivé cette définition géné\-rale, et lui ont servi de modèles, sont, bien entendu, apparus beaucoup plus tôt:
\begin{itemize}
	\item le groupoïde différentiable $PP^{-1}$ des isomorphismes de fibre sur fibre, attaché à un fibré principal et souvent appelé aujourd'hui \emph{groupoïde de jauge}, lequel joue un rôle clef dans la définition géométrique des connexions infinitésimales telles que C. Ehresmann les concevait, a été introduit dès 1950 au Colloque de Topologie de Bruxelles [EOC/28/];
	\item les \emph{groupoïdes de jets} sont introduits dans une Note aux C.R.A.S. de 1952 (\textbf{234}, pp. 1424-25) [EOC/38/].
\end{itemize}
\par\bigskip

Il est à noter que, dans ces deux cas particuliers, qui ont motivé la définition générale ultérieure, le transiteur $\tau=(\beta,\alpha)$ est en fait une surmersion, ce qui, bien sûr, implique \emph{a fortiori} qu'il en soit de même pour $\alpha$ et $\beta$. Mais C. Ehresmann a bien vu (motivé vraisemblablement par la considération d'autres exemples de catégories structurées) que seule cette dernière condition, plus faible, devait être retenue pour le cas général, malgré l'intérêt spécial qu'il portait à ce cas privilégié fondateur.

Notons aussi qu'à la suite de Ngo van Que, de nombreux auteurs, notamment A. Kumpera et K. Mackenzie, ont adopté le terme de \emph{groupoïde de Lie} pour qualifier le cas restrictif précédent, et cette terminologie s'est maintenue jusque vers la fin des années 80. On doit à P. Dazord et A. Weinstein l'initiative (qui semble assez largement acceptée aujourd'hui, et que nous suivons ici) d'avoir étendu cette terminologie au cas général, comme substitut à celle de \emph{groupoïde différentiable}. Cette terminologie (qu'Ehresmann n'avait jamais voulu adopter, malgré ma suggestion) a l'avantage de rétablir la cohérence avec celle d'\emph{algébroïde de Lie}, que j'avais adoptée en 1967 pour qualifier l'objet infinitésimal général correspondant (connu antérieurement seulement dans les deux cas particuliers sus-mentionnés: déplacements infinitésimaux des fibres), et de conforter le terme de \emph{foncteur de Lie} pour dénommer la correspondance groupoïde de Lie $\rightarrow$ algébroïde de Lie, comme étendant la correspondance groupe de Lie $\rightarrow$ algèbre de Lie.
\smallskip

Le thème principal du présent exposé est précisément d'insister sur le cas particulier \emph{complémentaire du précédent}: la condition de submersion pour $\tau$, autrement dit de \og\emph{surjectivité locale}\fg\,, liée aux \emph{fibrations}, fait alors place à celle d' \og\emph{injectivité locale}\fg\,, (c'est-à-dire diverses formes affaiblies de plongements ou d'immersions), liée aux \emph{feuilletages}.

Observons que, bien que la définition générale de C. Ehresmann englobe, nous venons de le souligner, ce dernier cas, et qu'il ait introduit, au Congrès de Montréal de 1961 [EOC/54/], la notion de \emph{groupoïde d'holonomie} d'un feuilletage topologique, en tant que groupoïde \emph{topologique}, il ne semble pas s'être intéressé à cet exemple, qu'il ne mentionne jamais dans ses présentations ultérieures des groupoïdes différentiables, estimant (de façon certainement abusive, comme nous le montrons) que la différentiabilité est une hypothèse superflue dont il convient de s'affranchir dans l'étude de l'holonomie des feuilletages.

Rappelons que le groupoïde d'holonomie différentiable d'un feuilletage régulier (que nous avions utilisé dans une Note aux C.R.A.S. de 1966 pour illustrer les obstacles nouveaux rencontrés dans l'étude de la correspondance sous-algébroïde de Lie $\rightarrow$ sous-groupoïde de Lie) a été redécouvert en 1983 par Winkelnkempen sous le nom de \og graphe d'un feuilletage\fg\, et simultanément exploité par A. Connes [NCG].

\bigskip
Nous voudrions également profiter de ces remarques historiques pour souligner la difficulté que les groupoïdes de Lie ont eu et ont encore à s'imposer comme source de la Géométrie, selon la conception visionnaire de C. Ehresmann.

Bien que cette situation ait évolué, très lentement, en près d'un demi-siècle, et un peu plus rapidement à partir des années 80, sous l'impulsion (entre autres) des travaux (que nous venons d'évoquer) d'Alain Connes, lequel utilise les groupoïdes de Lie comme source de Géométrie Non Commutative [NCG], et d'Alan Weinstein [WDC] (parallèlement à M. V. Karasev [KM]), qui les exploite comme sources de réalisations symplectiques des variétés de Poisson, et aussi grâce à la publication du livre de référence de K. Mackenzie [McK]( dont une nouvelle version est sous presse), cette notion continue de souffrir de la double malédiction des fées malicieuses qui se sont penchées sur son berceau à l'époque de sa naissance, et qui ont contribué à en éloigner, pour des raisons opposées (tenant surtout à des malentendus sémantiques), tant les Géomètres purs que les Catégoriciens purs, alors que les groupoïdes différen\-tiables sont à la fois des variétés et des catégories:
\smallskip

-- d'une part, les fondateurs de la Théorie des Catégories avaient à l'origine principalement en vue de grandes catégories de morphismes. Dans ce contexte, les grou\-poïdes ne peuvent décrire que des isomorphismes. En outre, en tant qu'objets algébriques, les groupoïdes, étant équivalents à des groupes, paraissent à première vue des objets triviaux, donc inintéressants, et, bien que les foncteurs entre grou\-poïdes ne soient en fait pas du tout, nous le verrons, des objets triviaux, peu d'exemples élaborés de foncteurs entre groupoïdes se présentent de façon naturelle dans ce contexte des grandes catégories. De plus la hantise des paradoxes ensemblistes, avait développé chez les Catégoriciens purs la religion des ensembles \og Hom\fg\, et le tabou de la considération globale des ensembles (ou plutôt des classes) d'objets ou de flèches, de sorte que la définition même des groupoïdes différentiables leur paraissait au départ suspecte d'hérésie;
\smallskip

-- d'autre part, la plupart des Géomètres (et des \og Working Mathematicians\fg\, [McL] en général) avaient la plus grande méfiance envers tout ce qui paraissait se rapprocher de l'\og abstract nonsense\fg\,, et une forte tendance à fuir à la seule audition du mot \og catégorie\fg\,. Malheureusement C. Ehresmann, loin de chercher à les rassurer en faisant observer que, dans la plupart des applications géométriques des groupoïdes différentiables, on peut se contenter de petites catégories (et en fait quasiment ignorer la théorie des grandes catégories), s'est au contraire ingénié --- tout en adoptant de fait, vis-à-vis des paradoxes logiques, une attitude pragmatique et \og naïve\fg\,, reposant sur une théorie des classes non véritablement explicitée --- à transporter systématiquement le langage des grandes catégories même là où les petites suffisaient, allant jusqu'à proposer d'appeler \og foncteurs\fg\, les homomorphismes de groupes (dans le présent exposé nous adopterons au contraire le terme de \og morphismes\fg\, pour les foncteurs entre petits groupoïdes, conformément au vocabulaire général des structures algébriques, dont les petits groupoïdes sont des cas particuliers, au même titre que les groupes, et nous substituerons le terme d'\og acteur\fg\, à celui de \og foncteur d'hypermorphisme\fg\, d'Ehresmann).

Il faut dire aussi qu'Ehresmann s'est assez rapidement détourné des applications géométriques, pourtant très prometteuses dans ses premières Notes, des groupoïdes de Lie, se trouvant entraîné dans l'étude générale des \og catégories structurées\fg\,. Or nous pensons que le cadre formel dans lequel il a développé cette étude (à savoir les catégories structurées par des catégories complètes), bien que conduisant à des résultats très généraux et profonds, laisse échapper, lorsqu'on cherche à les appliquer au cas particulier des groupoïdes de Lie, la plupart des propriétés les plus intéressantes et les plus utiles, notamment en ce qui concerne les quotients, que l'on rencontre aussi bien dans l'étude des fibrés que dans celle des feuilletages et des espaces de feuilles. Nous suggérons dans le présent exposé quelques pistes qui, privilégiant notamment le \og théorème de Godement\fg\,, nous paraissent plus adaptées à l'étude des applications géométriques.

Nous souhaiterions que le présent exposé puisse modestement contribuer à dissiper ces malentendus.

\section{Axiomes des diptyques.}
\label{ad}
(Les notations sont celles de \ref{dip}).
\subsection{Liste des axiomes.}
\begin{enumerate}
	\item $D_i\cap D_s=D_\ast$ ;
	\item $(f,f'\in D_i/D_s)\Longrightarrow f\times f'\in D_i/D_s$ (stabilité par produit);
	\item 
\begin{enumerate}
	\item les flèches de $D_i$ sont des monomorphismes; 
	\item les flèches de $D_s$ sont des épimorphismes \emph{stricts} [McL]\footnote{La légère dissymétrie des axiomes est intentionnelle, au vu de la plupart des exemples fondamentaux. La catégorie duale ne sera pas toujours un diptyque.};

\end{enumerate}
	\item 
\begin{enumerate}
	\item $h=gf\in D_i\Longrightarrow f\in D_i$ ($\alpha$-stabilité forte)\,\footnote{Dans les schémas qui suivent, on convient que les points d'interrogation espagnols \IE ? encadrent des objets ou des propriétés dont l'existence figure dans les conclusions des axiomes; de même l'existence des flèches en tirets résultera de l'énoncé des axiomes, par contraste avec les données et les hypothèses en traits pleins.};
	\begin{diagram}[h=1.5em,w=2em,labelstyle=\scriptstyle]
&{}&&&\rLine[abut]\til{h=gf}&&&{}&\\
\ruiLine(1,1)[abut]\ \ \ A\ \ \ \IE&\riLine&
?&\rTo\til f&B&&\rTo\til g&&C\rdTo(1,1)[abut]\\
\end{diagram}
	\item $(gf\in D_s\,\&\,f\in D_s)\Longrightarrow g\in D_s$ ($\beta$-stabilité faible);
	\begin{diagram}[h=1.5em,w=2em,labelstyle=\scriptstyle]
&{}&&&\rLine[abut]\til{h=gf}&&&{}&\\
A\ruLine(1,1)[abut]&&\rsTo\til f&&B&\rLine\til g&\IE&\rsTo&            {?\ \ \ C\ \ \ \ }\rdsTo(1,1)[abut]\\
\end{diagram}

\end{enumerate}
\item
(\emph{transversalité} : $D_s\pitchfork D_i$) :
         \begin{enumerate}
         \item(\emph{transfert parallèle}) :\\étant donné un couple $(A\rsTo^s B)\in D_s\,,(B'\riTo^i B)\in D_i$ , celui-ci se complète en un carré cartésien comme 
\begin{diagram}[h=3.4em,w=4em,tight,labelstyle=\scriptstyle,textflow]
        {\IE A'?}     &\rDashiTo\til{\IE i'?}            &         A   \\
\dDashsTo\til{\IE s'?}         &            {\IE\pb?}    & \dsTo\til s       \\
        B'     &\riTo\til i            &         B   \\
  \end{diagram}
indiqué, avec de plus $s'\in D_s$ et $i'\in D_i$ ;
          \item [(b)]inversement (\emph{descente}, ou \emph{transfert inverse}) :\\ étant donné un carré cartésien comme ci-dessous
\begin{diagram}[w=2em,h=3.4em,tight,labelstyle=\scriptstyle,textflow]
        A' &&    \riTo\til{i'} &           &         A   \\
\dsTo\til{s'}  &&            {\pb} &   & \dsTo\til s       \\
{\ \ B'\ \IE}&\riLine[abut]&{?}   &\rTo[abut]\til {i}            &         B   \\
  \end{diagram}
(où l'on suppose donc, selon nos conventions, $s,\!s'\in D_s$, $i'\in D_i, i\in D$), on a $i\in D_i$.
            \end{enumerate}
\end{enumerate}

\subsection{Quelques conséquences des axiomes.}
\label{qca}
\begin{itemize}
	\item Pour toute flèche $f:A\rightarrow B$ de $D$, son graphe $(\text{id}_A,f):A\riTo A\times B$ est un bon mono;
	\item Avec les notations de l'axiome (5), on a aussi:
\begin{itemize}
	\item $D_s\pitchfork D$;
	\item $D_s\pitchfork D_s$; les carrés cartésiens dont les quatre côtés sont dans $D_s$ seront appelés \emph{parfaits};
\end{itemize}
\item Tout carré parfait est aussi cocartésien (\emph{push out square})\,\footnote{Cette propriété importante est en fait équivalente à l'axiome (3)(b) en présence des autres axiomes.}.
\end{itemize}

\section{Conditions de dénombrabilité.}
\label{den}
 Les notations sont ici celles du paragraphe \ref{genot}, auquel la présente annexe se réfère.

Même pour une variété séparée les deux conditions suivantes ne sont pas équi\-va\-lentes. La première, qui dans le cas séparé équivaut à la paracompacité (jointe à la dé\-nom\-bra\-bi\-lité de l'ensemble des composantes connexes quand on l'impose à $X$ tout entière), est la plus satisfaisante, et définit des sous-catégories pleines qui déterminent encore des diptyques de Godement. Mais on est parfois contraint techniquement à travailler avec la seconde.

$X$ sera dite \emph{à base dénombrable} si elle vérifie l'une des quatre conditions équi\-va\-lentes suivantes:
\begin{itemize}
	\item [a)] $X$ possède un atlas dénombrable;
	\item [b)]$X$ est réunion dénombrable de compacts;
	\item [c)] $X$ admet une base dénombrable d'ouverts;
	\item [d)] $X$ est un espace de Lindelöf [TOP].
\end{itemize}

$X$ est dite \emph{séparable}\,\footnote{Cette terminologie désastreuse (surtout en français, à cause de l'emploi concurrent de \emph{séparé}, qui n'a rien à voir), attend toujours un substitut acceptable et accepté.} si elle vérifie l'une des deux conditions (équivalentes) suivantes:
\begin{itemize}
	\item [(i)] $X$ admet un sous-ensemble dense dénombrable;
	\item [(ii)] $X$ admet une sous-variété ouverte dense à base dénombrable.
\end{itemize}

\section{Immersions en tout genre.}\label{itg}
 Les notations sont ici encore celles du paragraphe \ref{genot}, auquel la présente annexe se réfère.

Nous exploitons maintenant des propriétés spécifiques à la catégorie \textbf{Dif}, fournies par la donnée de la sous-catégorie des \emph{immersions}, laquelle renferme tout une hiérarchie de notions d'injectivité, mêlant étroitement le local et le global, dont nous précisons les définitions, souvent trop mal connues.

Transférées au transiteur (\ref{am}), ces notions créent l'interaction (annoncée dans l'introduction) de la structure algébrique triviale (localement!) de graphe et de la structure de variété, pour faire naître des groupes non triviaux.

\subsection{Plongements.}

Dans le cadre des diptyques (de Godement), la notion qui joue dans \textbf{Dif} le rôle des injections (et qui est la plus restrictive parmi les immersions) est celle de \emph{plongement} [VAR] (\emph{propre} ou fermé quand on se place dans la sous-catégorie pleine \textbf{DifHaus}), laquelle est associée, via le théorème de Godement, à la caractérisation (par leur graphe) des équivalences \emph{régulières}, donc à la notion de feuilletage \emph{simple}.

Encore plus restrictives sont les notions de plongements ouverts et de plongements ouverts propres, qui sont des étalements injectifs.

On a besoin aussi de notions \emph{plus faibles}, qui sont fondamentales pour l'étude des feuilletages. Notons que celle d'immersion injective est la moins intéressante (bien que certains auteurs utilisent encore improprement le terme de sous-variétés pour les variétés immergées), et nous n'aurons pas besoin d'introduire de notation spéciale pour de telles flèches, qui définissent la sous-catégorie $\gr{I}_{\text{inj}}$.

\subsection{Plongements faibles et immersions strictes}\label{pfis}
Il s'agit en fait de la traduction, dans le cas particulier des variétés, de la notion générale de $p$-sous-structure au sens d'Ehresmann, où $p$ va être:
\begin{itemize}
	\item soit le foncteur d'oubli: $\mbf{Dif}\rightarrow\mbf{Ens}$;
	\item soit le foncteur d'oubli: $\mbf{Dif}\rightarrow\mbf{Top}$.
\end{itemize}

Nous dirons que la flèche $f:A\rightarrow B$ est un \emph{plongement faible} (resp. une \emph{immersion stricte}) si c'est une \emph{immersion}, et si, pour toute variété $Z$ et toute \emph{application ensembliste} (resp. \emph{continue}) $g:Z\rDotTo A$, la condition que l'application composée $h=f\circ g:Z\rightarrow B$ soit un morphisme (de variétés) implique qu'il en est de même pour $f$\,\footnote{La première de ces deux notions (la plus forte) a été introduite, sous une forme équivalente, mais d'aspect moins catégorique (ou du moins intensément utilisée, sans pourtant y recevoir de nom!), dans la thèse de R. Barre [RB], où elle joue un rôle crucial dans la définition des \textbf{Q}-\emph{variétés}, ainsi nommées parce qu'elles décrivent un cas très particulier (impliquant notamment la nullité de l'holonomie) d'espace \emph{quotient} d'un feuilletage (espace des feuilles). Nous l'avions par la suite utilisée, sous le nom d'immersion stricte (repris ensuite par P. Molino), et ultérieurement de plongement propre, comme propriété de l'immersion d'une feuille. D'autres auteurs utilisent celui de structure initiale. Désormais nous distinguons les deux termes, en affaiblissant le sens du terme \og immersion stricte \fg\, pour pouvoir l'utiliser comme propriété de l'application identique de la structure fine d'un feuilletage. Voir les exemples ci-après.}.

Notations: $\rplfTo$ et$\ristTo$.

\subsection{Immersions fidèles.}\label{ifid}
\'{A} toute variété $X$ est associée (voir [VAR]) la variété (\emph{non séparée} et non séparable!) $\text{J}(X)$ de ses \emph{germes de sous-variétés}, munie d'une \emph{immersion canonique} (surjective!):
$$\text{j}_X:\text{J}(X)\rimTo X$$
(laquelle possède d'ailleurs deux sections canoniques, associées aux germes discrets et grossiers).
Celle-ci est \emph{universelle} en ce sens que toute immersion $Y\rimTo X$ admet une factorisation unique de la forme:
$$f=\text{j}_X\circ u\circ s$$
où $s:Y\retasTo U$ est un étalement surjectif, et $u:U\rploTo X$ l'injection canonique d'un ouvert de $\text{j}(X)$ (plongement ouvert, ou étalement injectif)\,\footnote{Plus généralement toute \emph{subimmersion} [VAR] se factorise uniquement par une surmersion sur un ouvert de $\text{J}(X)$.}. Ehresmann dit alors que $Y$ est une \emph{variété extraite} de $X$.

On dira que l'immersion $f$ est \emph{fidèle} si la factorisation $s$ ci-dessus est injective (donc un isomorphisme).

Notation: $\rfidTo$.

La donnée d'un feuilletage de $X$ équivaut à celle d'un ouvert de $\text{J}(X)$ tel que la restriction de $\text{j}_X$ soit une bijection sur $X$.

On a les implications:
\begin{center}
\begin{tabular}[h]{|c|}
\hline
$\text{plongement fermé (ou ouvert)}\Longrightarrow\text{plongement}\Longrightarrow\text{plongement faible}\Longrightarrow$\\
$\text{immersion stricte}\Rightarrow\text{immersion injective}\Rightarrow\text{immersion fidèle}\Rightarrow\text{immersion.}$\\
\hline
\end{tabular}
\end{center}

\subsection{Exemples.}\label{exim}

Les exemples très simples qui suivent vont illustrer la hiérarchie stricte de ces diverses notions (on va des notions les plus fortes vers les moins fortes):
\begin{enumerate}
	\item un cercle privé d'un point, ou une spirale s'enroulant autour d'un point asymptotique, sont des sous-variétés plongées non proprement dans le plan;
	\item une géodésique dense du tore est faiblement plongée, non plongée;
	\item il en est de même de l'injection canonique de $\mbb Q$ dans $\mbb R$, ou de $\mbb{Q}\times\mbb{R}$ dans $\mbb{R}\times\mbb{R}$;
	\item l'application identique de $\dot{\mbb R}$ vers $\mbb R$, ou de $\dot{\mbb R}\times\mbb{R}$ vers $\mbb{R}\times\mbb{R}$ (où $\dot{\mbb R}$ désigne $\mbb R$ muni de la structure discrète), est une immersion stricte, mais non un plongement faible;
	\item la même situation se produit, de façon plus subtile, dans l'exemple connu (voir M. Spivak [SPI]) d'un feuilletage à une seule feuille d'une variété (séparée mais non paracompacte) à l'aide de la droite transfinie; c'est un feuilletage régulier, mais non strict au sens de \ref{fst}.
	\item le chiffre 6 est une variété fidèlement (car injectivement), mais non strictement immergée;
	\item la lettre $\alpha$ est une variété fidèlement (mais non injectivement) immergée;
	\item le cercle est l'image de la droite par une immersion non fidèle;
	\item le chiffre 3 n'est pas l'image d'une immersion.
\end{enumerate}

\subsubsection{Remarque.}\label{exre}
\begin{enumerate}
	\item 
Munissons le plan ${\mbb R}\times{\mbb R}$ de la structure de variété définie par:
$$({\dot{\mbb R}^-}\times{\mbb R})\cup({{\mbb R}^{+\ast}}\times{\dot{\mbb R}}).$$
Alors l'application identique détermine un préfeuilletage strict, qui n'est pas un feuilletage de \St\!.
	\item
On peut aussi aisément incorporer le chiffre 6 dans un préfeuilletage (non strict et non de \St\!, contenant une ou plusieurs copies du 6) en entourant le 6 par des courbes parallèles, de diverses façons.
\end{enumerate}

\smallskip

\section{Groupoïde d'holonomie topologique (C. Ehresmann).}\label{ght}
\subsection{Feuilletages topologiques}
La notion de \emph{groupoïde d'holonomie} d'un feuilletage a été introduite par C. Ehresmann en 1961 dans [EOC/54/] dans le cadre très faible des \emph{feuilletages topologiques}.

Dans ce cadre, $X$ est un espace topologique, et $\text{J}(X)$ est l'espace topologique des \emph{germes de sous-espaces} de $X$. La topologie fine peut être identifiée à un ouvert $X'$ de $\text{J}(X)$ tel que la restriction de $\text{j}_X$ soit bijective\,\footnote{Il s'agit donc de l'analogue topologique de ce que nous appelons ici immersion fidèle. Ehresmann utilise le terme de \emph{plongement}, en contradiction avec l'usage de ce terme qui s'est imposé à la suite de [VAR] dans le cadre des variétés, usage que nous respectons ici, et qui nous contraint à adopter une autre terminologie.}. On la suppose de plus \emph{localement connexe}.

\'{A} chaque ouvert $U$ de $X$, on peut ici associer l'\emph{espace topologique} $Q(U)$ quotient de $U$ par la relation d'équivalence définie par les composantes connexes (ou feuilles) de la topologie $U'$ \emph{induite} par $X'$ sur $U$. Cette notion est fonctorielle pour les inclusions d'ouverts, \cad qu'une telle inclusion $\rho:U\rightarrow V$ détermine une application \emph{continue} $Q(\rho):Q(U)\rightarrow Q(V)$.

La \emph{simplicité} (topologique) au sens d'Ehresmann signifie que la topologie $X$ admet une base d'ouverts $\go U$ telle que, pour tout $U\in\go U$, l'application canonique $Q(U)\rightarrow Q(X)$ soit un plongement ouvert (\cad un homéomorphisme sur un ouvert).

La \emph{locale simplicité} signifie que tout point de $X$ possède un voisinage ouvert $U$ dans $X$ tel que le feuilletage \emph{induit} par $(X,X')$ soit simple.

Ehresmann effectue la construction du groupoïde d'holonomie topologique sous l'hypothèse de locale simplicité.

\smallskip
\subsection{Remarques.}\label{ghtr}
\begin{enumerate}
	\item Le groupe d'holonomie \emph{topologique} d'un feuille sigulière ponctuelle est toujours nul. L'holonomie topologique ne peut donc apporter aucune information sur le voisinage d'un telle feuille.
	\item Les exemples de la remarque \ref{exre} sont des feuilletages topologiquement localement simples, qui ne sont pas des feuilletages de \St. La construction du groupoïde topologique d'Ehresmann s'applique, mais n'éclaire guère sur la nature de la singularité de ce feuilletage.
	\item Nous donnons maintenant un exemple simple de feuilletage de \St du plan \emph{non topologiquement localement simple}, auquel la construction d'Ehresmann ne s'applique donc pas.
	
	Nous prenons l'origine comme feuille singulière, et définissons en dehors de l'origine un feuilletage régulier dont les feuilles sont les cercles tangents à l'axe des abscisses (privés de l'origine), et les deux demi-droites. On a ainsi construit un feuilletage de \St sur le plan.
	
	Il est facile d'expliciter un groupoïde de Lie dont il est le feuilletage orbital. En effet ce feuilletage peut être défini par un champ de vecteurs s'annulant à l'origine, lequel détermine un flot (action locale), donc, par la construction décrite en \ref{aloc}, un groupoïde de Lie associé au feuilletage.
\end{enumerate}

\providecommand{\bysame}{\leavevmode\hbox to3em{\hrulefill}\thinspace}


\begin{thebibliography}{99}
\bibitem[Ast. 116]{Ast} Astérisque 116, \emph{Structure transverse des feuilletages}, Société Mathématique de France, Paris, 1984.

					\textsc{A. Haefliger}, Groupoïdes d'holonomie et classifiants, 70-97.

					\textsc{W. T. van Est}, Rapport sur les S-atlas, 235-292.

\bibitem[RB]{RB}\textsc{R. Barre}, Quelques propriétés des Q-variétés, \emph{Ann. Inst. Fourier}, Grenoble, 23, (1973), 227-312.
\bibitem[VAR]{B} \textsc{N. Bourbaki}, \emph{Variétés différentielles et analytiques}, Hermann, Paris, 1971.
\bibitem[LIE]{LIE} \bysame, \emph{Groupes et algèbres de Lie}, Hermann, Paris, 1972
\bibitem[TOP]{TOP} \bysame, \emph{Topologie générale}, Hermann, Paris, 1974.
\bibitem[NCG]{NCG} \textsc{A. Connes}, \emph{Noncommutative Geometry}, Acad. Press, 1994.
\bibitem[EOC]{OC} \textsc{C. Ehresmann}, \emph{\OE uvres Complètes}, Amiens, 1984.\\
					/28/, Les connexions infinitésimales dans un espace fibré différentiable, I, 179-205, et \emph{Coll. de Topo.} Bruxelles, C.B.R.M. (1950), 29-55.\\
					/38/, Les prolongements d'une variété différentiable. V. Covariants différentiels et prolongements d'une structure infinitésimale, I, 358-360, et \emph{C.R.A.S., Paris}, 234, (1952) , 1424-25.\\
					/50/, Catégories topologiques et catégories différentiables, I, 237-250, et \emph{Coll. Géom. Diff. Globale,}, Bruxelles, C.B.R.M. (1959), 137-150.\\
					/54/, Structures feuilletées, II, 2, 563-626, et \emph{Proceedings of the Fifth Canadian Mathematical Congress, Montréal} (1961), 109-172.\\
					/103/ Sur les catégories différentiables, I, 261-270, et \emph{Atti Conv. Int. Geom. Diff.} Bologna (1967, 31-40).
\bibitem[GZ]{GZ} \textsc{P. Gabriel, M. Zisman}, \emph{Calculus of fractions and homotopy theory}, Ergebn. Math. 35, Springer, (1965).
\bibitem[KM]{KM} \textsc{M. V. Karasev, V. P. Maslov}, \emph{Nonlinear Poisson brackets, Geometry and Quantization}, Transl. of Math. Monographs, vol. 119, (1993), A.M.S., Rhode Island.
\bibitem [McK]{McK} \textsc{K. Mackenzie}, \emph{Lie groupoids and Lie algebroids in Differential Geometry}, Cambridge Univ. Press, Cambridge, (1987) (Nouvelle version sous presse).
\bibitem[McL]{ML} \textsc{S. Mac Lane}, \emph{Categories for the Working Mathematicians}, Springer-Verlag, New York, 1971.
\bibitem[P 1966]{P66} \textsc{J.Pradines}, Théorie de Lie pour les groupoïdes différentiables. Relations entre propriétés locales et globales, \emph{C.R.A.S., Paris}, 263, (1966), 907-910.
\bibitem[P 1984]{P84} \bysame, Feuilletages: holonomie et graphes locaux, \emph{C.R.A.S., Paris}, 298, 13 (1984), 297-300.
\bibitem[P 1985]{P85} \bysame, How to define the graph of a singular foliation, \emph{Cahiers Top. et Géom. Diff.}, XVI-4 (1985), 339-380.
\bibitem[P 1986]{P86} \bysame, Quotients de groupoïdes différentiables, \emph{C.R.A.S., Paris}, 303, 16 (1986), 817-820.
\bibitem[P 1989]{P89} \bysame, Morphisms between spaces of leaves viewed as fractions, \emph{Cahiers Top. et Géom. Diff.}, XXX-3 (1989), 229-246.
\bibitem[P 2003]{P03} \bysame, Lie groupoids as generalized atlases, \emph{preprint}, à paraître dans \emph{Central European Journal of Mathematics}.
\bibitem[LALG]{SLG} \textsc{J.-P. Serre}, \emph{Lie Algebras and Lie Groups}, W. A. Benjamin Inc., New York, 1965.
\bibitem[SPI]{SPI} \textsc{M. Spivak}, \emph{A comprehensive introduction to Differential Geometry}, vol. I, Brandeis University, (1970).
\bibitem[WDC]{WDC} \textsc{A. Weinstein, P. Dazord, A. Coste}, Groupoïdes symplectiques, \emph{Publications du département de Mathématiques de l'Université de Lyon 1}, 2/A (1987).
\bibitem[W]{W} \textsc{A. Weinstein}, Groupoids: unifying internal and external symmetry, \emph{Notices of the Am. Math. Soc.}, vol. 43, 7, (1996), 744-752.



\end{thebibliography}
\end{document}